\documentclass[11pt,a4paper]{amsart}
\usepackage{amsfonts}
\usepackage{mathrsfs}
\usepackage{amsthm,amsxtra}
\usepackage{amssymb}
\usepackage{amsmath}
\usepackage{amscd}
\usepackage{hyperref}
\usepackage{bookmark}
\hypersetup{
colorlinks=true,
linkcolor=blue,
citecolor=gray,
}
\usepackage[utf8]{inputenc}
\usepackage{t1enc}
\usepackage[mathscr]{eucal}
\usepackage{indentfirst}
\usepackage{graphicx}
\graphicspath{ {images/} }
\usepackage{subcaption}
\usepackage{enumerate}
\usepackage{graphics}
\usepackage{pict2e}
\usepackage{epic}
\numberwithin{equation}{section}
\usepackage[margin=2.9cm]{geometry}
\usepackage{epstopdf} 
\usepackage{longtable}
\usepackage{enumerate}
\usepackage[inline]{enumitem}
\usepackage{tikz-cd}
\usepackage{tikz}
\usepackage{color}
\DeclareMathOperator{\MCG}{MCG}

\DeclareMathOperator{\Out}{Out}
\DeclareMathOperator{\Aut}{Aut}
\DeclareMathOperator{\Inn}{Inn}
\DeclareMathOperator{\Stab}{Stab}

\newcommand{\FF}{\mathbb{F}}
\newcommand{\MM}{\mathbb{M}}
\newcommand{\ZZ}{\mathbb{Z}}

\theoremstyle{plain}
\newtheorem{theorem}{Theorem}[section]
\newtheorem{proposition}[theorem]{Proposition}
\newtheorem{lemma}[theorem]{Lemma}
\newtheorem{corollary}[theorem]{Corollary}

\theoremstyle{definition}

\theoremstyle{remark}

\newtheorem*{remark*}{Remark}

\begin{document}

\title[Dynamics of Dehn twists in \texorpdfstring{$\operatorname{Out}(\FF)$}{Out(F)}]{Dynamics of Dehn twists in the outer automorphism group of a free group}
	
\author{Donggyun Seo}
\address{Institute of Mathematical Sciences \\ Chungnam National University \\ Daejeon, Korea}
\email{seodonggyun@cnu.ac.kr \\ seodonggyun7@gmail.com}
\urladdr{https://www.seodonggyun.com}


\begin{abstract}
	We study Dehn twists in the outer automorphism group of a free group through the doubled handlebody, a connected sum of copies of the product of a circle and a two-sphere, in which the Dehn twist along a primitive torus is induced by twisting a regular neighborhood of the torus.
	Our main result states that sufficiently large powers of the Dehn twists along a compatible collection of primitive tori generate a right-angled Artin group.
	The proof is a ping-pong argument in the universal cover, the complement of a Cantor set in the three-sphere, where intersection numbers become counts of orbits of lifts.
	Along the way we show that commuting Dehn twists along primitive tori have tori whose lifts never cross, and that a Dehn twist determines the oriented cyclic splitting dual to its torus.
	Finally, we construct a compact space of primitive curves on which multitwists act parabolically and fully irreducible outer automorphisms act with north--south dynamics.
\end{abstract}

\maketitle

\section{Introduction}

Let $\FF = \FF_g$ be a free group of rank $g \ge 2$, and let $\Out(\FF) = \Aut(\FF) / \Inn(\FF)$ be its outer automorphism group.
Many questions about $\Out(\FF)$ have a topological counterpart in the doubled handlebody
\[
	\MM = \MM_g = \#_g (\mathbb{S}^1 \times \mathbb{S}^2) ,
\]
the double of the handlebody of genus $g$, whose fundamental group is $\FF$.
Whitehead \cite{MR1575455, MR1503309} used $\MM$ to solve the automorphic equivalence problem for $\FF$, and proved that the natural map $\MCG(\MM) \to \Out(\FF)$ from the mapping class group is surjective.
Laudenbach \cite{MR356056} showed that its kernel is $(\ZZ / 2 \ZZ)^g$, and Brendle, Broaddus and Putman \cite{MR4557874} that the resulting short exact sequence splits.
So $\Out(\FF)$ is a subgroup of $\MCG(\MM)$, and an outer automorphism may be studied through a homeomorphism of $\MM$.

A \emph{Dehn twist} in $\Out(\FF)$ is the outer automorphism obtained by twisting along a one-edge $\ZZ$-splitting of $\FF$; see Section~\ref{subsec:dehn_prelim}.
The conjugacy class of the edge group of the splitting is the \emph{core} of the twist.
An essential torus in $\MM$ determines such a splitting, its dual tree of Section~\ref{subsec:curve-lifts}, whose core is the homotopy class of a core curve of the torus.
We say that the torus is \emph{primitive} if its core is represented by an element of a free basis of $\FF$.

For a primitive torus, the Dehn twist along its dual tree is the image of the homeomorphism obtained by twisting a regular neighborhood of the torus, as for Dehn twists on surfaces, by Section~\ref{sec:dehn-twist}.
We do not claim that every one-edge $\ZZ$-splitting of $\FF$ is dual to a torus.
Dehn twists along tori and their dynamics have been studied by Clay and Pettet \cite{MR2720227} and G\"ultepe \cite{MR3677931}, among others.

This paper studies the subgroups of $\Out(\FF)$ generated by powers of Dehn twists along primitive tori, and the dynamics of these twists on a compact space of curves.
Throughout, a primitive torus enters the statements only through its dual tree, the one-edge $\ZZ$-splitting of $\FF$ it determines.
Two primitive tori are in the same \emph{class} if their dual trees are $\FF$-equivariantly isomorphic; homotopic tori are in the same class, by Section~\ref{subsec:classes}.
The Dehn twist $\delta_{\vec T} \in \Out(\FF)$ depends only on the class $T$ together with the orientations of the torus and of its core.
Two classes of essential spheres or of primitive tori are \emph{core-independent} if, when both are classes of tori, their cores are not freely homotopic.

\subsection{The main theorem}

For a finite simplicial graph $\Gamma$ with vertex set $V$, the \emph{right-angled Artin group} on $\Gamma$ is
\[
	A(\Gamma) := \langle \, v \in V \mid [u, v] = 1 \text{ whenever } u, v \text{ are adjacent} \, \rangle .
\]
For a homotopy class $\gamma$ of curves and a homotopy class $S$ of essential spheres or tori, the \emph{geometric intersection number} $i(\gamma, S)$ is the minimal number of points of $c \cap s$ over transverse representatives $c \in \gamma$ and $s \in S$.
See Section~\ref{subsec:primitive}.
For a class $T$ of primitive tori we put $i(\gamma, T) := i(\gamma, [t])$, where $[t]$ is the homotopy class of some $t \in T$; this does not depend on $t$, by Section~\ref{subsec:classes}.

For two core-independent classes $T$, $T'$ of primitive tori, we use instead the \emph{crossing number} $N(T, T')$.
It is the number of orbits of pairs of lifts of representatives to the universal cover that cross in the sense of Scott \cite{MR1608688}.
It is modeled on the intersection number of two splittings of Scott \cite{MR1608688}, and it vanishes when $T$ and $T'$ have disjoint representatives.

Let $V = \{ T_1, \dots, T_k \}$ be pairwise core-independent classes of primitive tori, with cores $\gamma_1, \dots, \gamma_k$.
The \emph{coincidence graph} $\Gamma$ of $V$ has vertex set $V$, with distinct $T_a$ and $T_b$ adjacent when they have disjoint representatives.
Twists along disjoint tori commute, so the assignment $T_a \mapsto \delta_{\vec T_a}^{\,n_a}$ defines a homomorphism $A(\Gamma) \to \Out(\FF)$ for all integers $n_a$.
We say that $V$ is \emph{compatible} if, whenever distinct $T_a$ and $T_b$ are not adjacent, the core of each has positive intersection number with the other.
Put $M(\Gamma) := \max_{a,b} i(\gamma_a, T_b)$, $I(\Gamma) := \max_{a \ne b} N(T_a, T_b)$, with $I(\Gamma) := 0$ when $k = 1$, and $D(\Gamma) := 2 (M(\Gamma) + I(\Gamma))^2$.

\begin{theorem}[Theorem~\ref{thm:main}] \label{thm:intro-main}
	Let $V = \{ T_1, \dots, T_k \}$ be a compatible collection with coincidence graph $\Gamma$.
	Then for all integers $n_1, \dots, n_k$ with $\lvert n_a \rvert \ge (2k + 7) D(\Gamma) + 4$, the homomorphism
	\[
		\Phi \colon A(\Gamma) \to \Out(\FF) , \qquad T_a \mapsto \delta_{\vec T_a}^{\,n_a} ,
	\]
	is injective.
\end{theorem}

The problem of constructing free subgroups of $\Out(\FF)$ has a long history.
Bestvina, Feighn and Handel \cite{MR1765705} proved the Tits alternative for $\Out(\FF)$.
Kapovich and Lustig \cite{MR2652906} showed that large powers of two hyperbolic fully irreducible elements generate a virtually cyclic or a free group of rank two, and Uyanik \cite{MR3954286} obtained a related result for hyperbolic extensions.
Clay and Pettet \cite{MR2720227} proved that if two one-edge $\ZZ$-splittings fill, then large powers of the associated Dehn twists generate a free group of rank two.
They also proved that its nontrivial elements, apart from conjugates of the twists, are atoroidal and fully irreducible.
G\"ultepe \cite{MR3677931} gave another criterion for two Dehn twists to generate a free group.

Bering \cite{MR3946718} considered two Dehn twists along one-edge $\ZZ$-splittings of $\FF$ in \emph{hyperbolic--hyperbolic} position, the edge generator of each splitting acting hyperbolically on the tree of the other.
He proved that the powers $\delta^n$ and $\delta'^n$ generate a free group of rank two for all $\lvert n \rvert$ at least a constant depending only on the rank of $\FF$.
By Proposition~\ref{prop:orbit-count-translation}, the condition $i(\gamma_a, T_b) > 0$ says that the core of $T_a$ acts hyperbolically on the dual tree of $T_b$.
So the compatibility of a pair of classes of primitive tori is Bering's hyperbolic--hyperbolic condition, and for a pair, Theorem~\ref{thm:intro-main} is his result with a threshold depending on the crossing number instead of the rank.

What Theorem~\ref{thm:intro-main} adds is the finite collection with its prescribed commutation graph and independently chosen exponents.
The passage from a pair to a collection is where Proposition~\ref{prop:link-twist} is needed.
The twists in the link of a vertex must preserve enough of the lifts penetrated by a curve for the next step of the ping-pong to apply.

Taylor \cite{MR3343354} proved that every right-angled Artin group embeds quasi-isometrically into some $\Out(\FF)$, using powers of fully irreducible elements supported on free factors.
Theorem~\ref{thm:intro-main} is a counterpart of Taylor's theorem in which the generators are Dehn twists, and it follows the ping-pong scheme that Koberda \cite{MR3000498} used for powers of Dehn twists on surfaces.

The compatibility hypothesis does not cover every right-angled Artin group generated by Dehn twists.
For instance, $\Out(\FF_g)$ contains a direct product of $2g-4$ free groups generated by Dehn twists, by Horbez and Wade \cite{MR4069231} and Bridson and Wade \cite{MR4792836}.
It also contains a free abelian group of rank $2g-3$ generated by Dehn twists, by Feighn and Handel \cite{MR2496054}.
In both examples several of the generating twists have freely homotopic cores, such as the twists $a_i \mapsto a_i a_1$ for different $i$.
So their tori are not pairwise core-independent, and the hypotheses of Theorem~\ref{thm:intro-main} do not hold for these collections.

Koberda's argument needs two inputs that are standard for curves on a surface.
First, Dehn twists along disjoint curves commute, while the twists along two curves that intersect are detected by their action on a third curve.
Second, intersection numbers can be read off from representatives in minimal position.
For primitive tori in $\MM$, the first input is replaced by Theorem~\ref{thm:intro-commuting} below together with the compatibility hypothesis, which makes the cores of non-adjacent members penetrate lifts of each other.
The second is replaced by counting orbits of crossing lifts, which needs no minimal position at all.
Theorems~\ref{thm:intro-commuting-crossing} and~\ref{thm:intro-coincide} are independent results about the twists themselves, and the proof of Theorem~\ref{thm:intro-main} does not use them.

\subsection{Dehn twists}

If $T$ and $T'$ have disjoint representatives, the twists have disjoint supports and commute.
Commuting twists have crossing number zero; whether crossing number zero conversely forces disjoint representatives is not known to us, see Section~\ref{subsec:commuting}.

\begin{theorem}[Theorem~\ref{thm:commuting}] \label{thm:intro-commuting}
	Let $T$ and $T'$ be core-independent classes of primitive tori, with any orientations.
	If $T$ and $T'$ have disjoint representatives, then $\delta_{\vec T}$ and $\delta_{\vec T'}$ commute.
\end{theorem}

\begin{theorem}[Theorem~\ref{thm:commuting-crossing}] \label{thm:intro-commuting-crossing}
	Let $T$ and $T'$ be core-independent classes of primitive tori, with any orientations.
	If $\delta_{\vec T}$ and $\delta_{\vec T'}$ commute, then $N(T, T') = 0$.
\end{theorem}

The proof of Theorem~\ref{thm:intro-commuting-crossing} rests on the fact that a Dehn twist remembers its splitting, which we deduce from a rigidity theorem of Cohen and Lustig \cite{MR1691946}.

\begin{theorem}[Proposition~\ref{prop:coincide}] \label{thm:intro-coincide}
	Let $\vec T$ and $\vec T'$ be oriented classes of primitive tori.
	Then $\delta_{\vec T} = \delta_{\vec T'}$ if and only if $\vec T = \vec T'$.
\end{theorem}

The second theme is the dynamics of Dehn twists.
Let $\mathcal{C}$ be the set of homotopy classes of unoriented primitive curves, that is, of simple closed curves representing elements of a free basis, and $\mathcal{S}$ the set of isotopy classes of essential spheres.
The intersection numbers $i(\gamma, S)$ define a map $i_\ast \colon \mathcal{C} \to \mathbb{R}_+^{\mathcal{S}}$, and its projectivization $p \circ i_\ast$ is injective by Proposition~\ref{prop:injective}, and the closure $\overline{\mathcal{C}}$ of its image is compact by Proposition~\ref{prop:C_compact}.
This is an analogue of the embedding of curves into the projective space $\mathcal{PML}(\Sigma)$ of measured laminations on a surface $\Sigma$.
The action of $\Out(\FF)$ on $\overline{\mathcal{C}}$ exhibits the two model behaviors of the action of a mapping class group on $\mathcal{PML}(\Sigma)$: multitwists act parabolically and fully irreducible elements act with north--south dynamics.

Let $t_1, \dots, t_r$ be pairwise disjoint primitive tori with oriented cores, whose classes $T_1, \dots, T_r$ are pairwise core-independent, with cores $\gamma_1, \dots, \gamma_r$.
Their Dehn twists commute, and we call a product $\Delta = \delta_{\vec T_1}^{\,m_1} \dots \delta_{\vec T_r}^{\,m_r}$ a \emph{multitwist}.

\begin{theorem}[Corollaries~\ref{cor:multitwist-parabolic} and~\ref{cor:multitwist-parabolic-converge}] \label{thm:intro-parabolic}
	Let $\Delta = \delta_{\vec T_1}^{\,m_1} \dots \delta_{\vec T_r}^{\,m_r}$ be a multitwist along pairwise disjoint primitive tori with pairwise core-independent classes $T_1, \dots, T_r$ and cores $\gamma_1, \dots, \gamma_r$, and let $\alpha \in \mathcal{C}$.
	If $i(\alpha, T_v) = 0$ for every $v$ with $m_v \ne 0$, then $\Delta^n(\alpha) = \alpha$ for all $n$.
	Otherwise $\Delta^n(\alpha)$ converges in $\overline{\mathcal{C}}$, as $\lvert n \rvert \to \infty$, to the projective class of
	\[
		\sum_{v=1}^{r} \lvert m_v \rvert \, i(\alpha, T_v) \, i_\ast(\gamma_v) .
	\]
\end{theorem}

\begin{theorem}[Corollary~\ref{cor:multitwist-parabolic-fixed}] \label{thm:intro-parabolic-fixed}
	Let $\Delta = \delta_{\vec T_1}^{\,m_1} \dots \delta_{\vec T_r}^{\,m_r}$ be a multitwist along pairwise disjoint primitive tori with pairwise core-independent classes $T_1, \dots, T_r$ and cores $\gamma_1, \dots, \gamma_r$.
	Then every projective class $\bigl[ \sum_v s_v \, i_\ast(\gamma_v) \bigr]$ with $s_v \ge 0$ not all zero is fixed by $\Delta$.
\end{theorem}

For $r = 1$ Theorem~\ref{thm:intro-parabolic} is the parabolic dynamics of a single Dehn twist: it fixes the core of its torus and pushes every primitive curve $\alpha$ with $i(\alpha, T) > 0$ toward this core.
For $r > 1$ the limit is a weighted combination of the cores.
For a fixed multitwist and a fixed $\alpha$ it is one weight vector, with coefficients proportional to $\lvert m_v \rvert \, i(\alpha, T_v)$, whose ratios are rational.
Every projective weighting on the cores of the tori met by $\alpha$ occurs as a limit along a sequence of multitwists.

\begin{theorem}[Corollary~\ref{cor:core_simplex}] \label{thm:intro-simplex}
	Suppose that some $\alpha \in \mathcal{C}$ satisfies $i(\alpha, T_v) > 0$ for all $v$.
	Then $\overline{\mathcal{C}}$ contains the projectivized cone $\bigl\{ \bigl[ \sum_v s_v \, i_\ast(\gamma_v) \bigr] \mid s_v \ge 0 \text{ not all zero} \bigr\}$.
\end{theorem}

By Theorem~\ref{thm:intro-parabolic-fixed}, this cone is fixed pointwise by every multitwist along $t_1, \dots, t_r$.
So $\overline{\mathcal{C}}$ contains, besides the primitive curves, the projectivized cones of weighted cores attached to the systems of disjoint primitive tori $t_1, \dots, t_r$ for which some $\alpha \in \mathcal{C}$ has $i(\alpha, T_v) > 0$ for all $v$.
These are the counterpart of the simplices of weighted multicurves in $\mathcal{PML}(\Sigma)$.
We do not know whether $i_\ast(\gamma_1), \dots, i_\ast(\gamma_r)$ are linearly independent, that is, whether these cones are simplices.
Theorems~\ref{thm:intro-parabolic}, \ref{thm:intro-parabolic-fixed} and~\ref{thm:intro-simplex} rest on an estimate for the growth of intersection numbers under multitwists, Proposition~\ref{prop:twist-growth}.
From this estimate, Corollary~\ref{cor:disjoint_twists_free_abelian} deduces that Dehn twists along disjoint core-independent primitive tori generate a free abelian group of the expected rank; Proposition~\ref{prop:pairwise-free-abelian} extends this to tori that are only pairwise disjoint up to homotopy, by results of Guirardel and Levitt on compatible trees and of Levitt on groups of twists, and this is the form used in the proof of Theorem~\ref{thm:intro-main}.

The counterpart of a pseudo-Anosov mapping class is a fully irreducible outer automorphism, which acts with north--south dynamics on the minimal set of projective geodesic currents by Martin \cite{MR2693216} and Uyanik \cite{MR3273532, MR3370027}.
Section~\ref{subsec:iwip} supplies one step of Uyanik's argument for nonorientable surfaces.
Through the intersection form of Kapovich and Lustig \cite{MR2496058}, this dynamics descends to $\overline{\mathcal{C}}$.

\begin{theorem}[Theorem~\ref{thm:loxodromic}] \label{thm:intro-loxodromic}
	A fully irreducible $\phi \in \Out(\FF)$ acts on $\overline{\mathcal{C}}$ with north--south dynamics.
\end{theorem}

\subsection{Acknowledgments}

The author thanks Sang-hyun Kim for many helpful and inspiring conversations during the preparation of this work.
During the revision process, discussions with Michael Kopreski, Thomas Hill, George Shaji and Brian Udall highlighted an error in a previous version of Corollary~\ref{cor:dehn_disj}, leading to its improvement.
An error in an earlier version of this work was corrected thanks to comments of Sanghoon Kwak, Richard Wade and the referees.
The author is also grateful to Thomas Koberda for valuable comments that strengthened the manuscript, and to Ignat Soroko for pointing out Whitehead's work.

This paper is dedicated to my son, Jian, who arrived as this manuscript was taking shape.

\section{Preliminaries} \label{sec:prelim}

This section fixes notation and collects the standard facts used in the paper.
These are one-edge $\ZZ$-splittings and Dehn twists in $\Out(\FF)$, the doubled handlebody $\MM$ and its mapping class group, and the results from $3$-manifold topology and from the theory of trees that we cite by name.
The section ends with the primitive curves and Lemma~\ref{lem:primitive-basis}.

\subsection{Free groups and one-edge \texorpdfstring{$\ZZ$}{Z}-splittings} \label{subsec:free_group}

This subsection recalls standard facts on free groups and trees and defines one-edge $\ZZ$-splittings and their orientations.

Let $g \ge 2$, and let $\FF = \FF_g$ be the free group of rank $g$.

\begin{theorem} \label{thm:free}
	Let $a, b \in \FF \setminus \{ 1 \}$.
	\begin{enumerate}
		\item \label{item:free:commute}
		The centralizer of $a$ is the maximal cyclic subgroup containing $a$.
		So $a$ and $b$ commute if and only if they are powers of a common element, if and only if they have a common nonzero power.
		\item \label{item:free:basis}
		If $a$ and $b$ do not commute, then they freely generate a free subgroup of rank $2$.
		\item \label{item:free:primitive}
		An element of a free basis of $\FF$ is not a proper power, and its centralizer is the cyclic group it generates.
	\end{enumerate}
\end{theorem}

\begin{theorem} \label{thm:trees}
	Let a group $G$ act on a tree $X$ without edge inversions, and let $h \in G$.
	\begin{enumerate}
		\item \label{item:trees:dichotomy}
		Either $h$ fixes a vertex, and its fixed point set is a subtree, or $h$ is \emph{hyperbolic}.
		In the latter case there is a unique $h$-invariant line, the \emph{axis} of $h$, on which $h$ translates by $\lVert h \rVert := \min_x d(x, h x) > 0$.
		Moreover $d(x, h^n x) = n \lVert h \rVert + 2 \, d(x, \mathrm{axis})$ for every vertex $x$ and $n \ge 1$.
		\item \label{item:trees:power}
		For $n \ne 0$, $h^n$ is hyperbolic if and only if $h$ is, with the same axis.
		\item \label{item:trees:fg}
		If $G$ is finitely generated and every element of $G$ fixes a vertex, then $G$ fixes a vertex.
		\item \label{item:trees:ends}
		A hyperbolic $h$ fixes the two ends of its axis and acts on the space of ends of $X$ with north--south dynamics.
		That is, for every compact set $K$ of ends not containing the repelling end and every neighborhood $U$ of the attracting end, $h^n K \subset U$ for all large $n$.
	\end{enumerate}
\end{theorem}

An \emph{$\FF$-tree} is a simplicial tree with a simplicial action of $\FF$ without edge inversions.
An $\FF$-tree $T$ is a \emph{one-edge $\ZZ$-splitting} if $\FF$ acts transitively on the edges of $T$ and the stabilizer of an edge is infinite cyclic.
The stabilizers of the edges are then conjugate, and the conjugacy class of a generator, up to inversion, is the \emph{core} of $T$.
We do not require $T$ to be minimal; the splittings to which this section is applied are the dual trees of primitive tori, which have no vertex of valence one by Lemma~\ref{lem:no-leaf}, so that they are minimal.

Let $\vec E T$ be the set of oriented edges of $T$, and for $\vec e \in \vec E T$ let $\ast \vec e$ be the same edge with the opposite orientation.
An \emph{orientation} of $T$ is a map $o \colon \vec E T \to \FF$ such that
\[
	\langle o(\vec e) \rangle = \Stab(e) , \qquad
	o(\ast \vec e) = o(\vec e)^{-1} , \qquad
	o(w \vec e) = w \, o(\vec e) \, w^{-1}
\]
for all $\vec e \in \vec E T$ and $w \in \FF$.
An orientation is determined by its value on one oriented edge, and there are exactly two orientations.
We write $\vec T := (T, o)$ for $T$ with an orientation.

The Cayley tree $\mathrm{Cay}(\FF)$ of $\FF$ with respect to a free basis is an $\FF$-tree on which every nontrivial element is hyperbolic.
The ends of $\mathrm{Cay}(\FF)$ form the Gromov boundary of $\FF$, which we identify with the Cantor set of the universal cover in Section~\ref{subsec:ends}.

\subsection{Dehn twists in \texorpdfstring{$\operatorname{Out}(\FF)$}{Out(F)}} \label{subsec:dehn_prelim}

This subsection defines the Dehn twist $\delta_{\vec T} \in \Out(\FF)$ along an oriented one-edge $\ZZ$-splitting $\vec T$.

Let $\vec T = (T, o)$ be an oriented one-edge $\ZZ$-splitting, and fix a vertex $v_0$ of $T$.
For a vertex $v$, let $\vec e_1, \dots, \vec e_k$ be the oriented edges of the geodesic from $v_0$ to $v$, and put
\[
	\delta_{\vec T, v_0}(v) := o(\vec e_1) \dots o(\vec e_k) \, v .
\]

\begin{lemma} \label{lem:twist-path}
	Let $\vec e_1, \dots, \vec e_k$ be the oriented edges of any edge path from $v_0$ to $v$.
	Then $\delta_{\vec T, v_0}(v) = o(\vec e_1) \dots o(\vec e_k) \, v$.
\end{lemma}

\begin{proof}
	If the path is not a geodesic, then $\vec e_{i+1} = \ast \vec e_i$ for some $i$, and $o(\vec e_i) \, o(\vec e_{i+1}) = 1$.
	Removing $\vec e_i$ and $\vec e_{i+1}$ gives a shorter path with the same product.
	Repeating this ends at the geodesic.
\end{proof}

\begin{lemma} \label{lem:twist-automorphism}
	The map $\delta_{\vec T, v_0}$ extends to a simplicial automorphism of $T$, whose inverse is given on vertices by $v \mapsto o(\vec e_1)^{-1} \dots o(\vec e_k)^{-1} \, v$ for the geodesic $\vec e_1, \dots, \vec e_k$ from $v_0$ to $v$.
\end{lemma}

\begin{proof}
	Let $v$ and $v'$ be adjacent, with $v'$ closer to $v_0$, and let $\vec e_1, \dots, \vec e_k$ be the geodesic from $v_0$ to $v$.
	Then $\vec e_1, \dots, \vec e_{k-1}$ is the geodesic from $v_0$ to $v'$, and $o(\vec e_k)$ fixes $v'$ since it stabilizes $e_k$.
	So $\delta_{\vec T, v_0}(v) = o(\vec e_1) \dots o(\vec e_k) \, v$ and $\delta_{\vec T, v_0}(v') = o(\vec e_1) \dots o(\vec e_k) \, v'$ are adjacent, and $\delta_{\vec T, v_0}$ extends to a simplicial map.

	Let $\epsilon$ denote the map in the statement.
	The $j$-th oriented edge of the geodesic from $v_0$ to $\delta_{\vec T, v_0}(v)$ is $\vec f_j := o(\vec e_1) \dots o(\vec e_{j-1}) \, \vec e_j$: these edges form a reduced path since $o(\vec e_j)$ fixes $e_j$, so $\vec f_j$ ends where $\vec f_{j+1}$ begins, and the path ends at $\delta_{\vec T, v_0}(v)$.
	Since $o(w \vec e) = w \, o(\vec e) \, w^{-1}$, we have $o(\vec f_j)^{-1} = o(\vec e_1) \dots o(\vec e_{j-1}) \, o(\vec e_j)^{-1} \, o(\vec e_{j-1})^{-1} \dots o(\vec e_1)^{-1}$.
	So $\epsilon(\delta_{\vec T, v_0}(v)) = o(\vec f_1)^{-1} \dots o(\vec f_k)^{-1} \, \delta_{\vec T, v_0}(v)$, where the product telescopes to $(o(\vec e_1) \dots o(\vec e_k))^{-1}$, which gives $v$.
	Likewise $\delta_{\vec T, v_0}(\epsilon(v)) = v$.
	Hence $\delta_{\vec T, v_0}$ is a bijection on vertices, and a simplicial automorphism.
\end{proof}

\begin{lemma} \label{lem:twist-conjugate}
	For $w \in \FF$, we have $w \, \delta_{\vec T, v_0} \, w^{-1} = \delta_{\vec T, w v_0}$.
\end{lemma}

\begin{proof}
	Let $v$ be a vertex and $\vec f_1, \dots, \vec f_m$ the geodesic from $w v_0$ to $v$.
	Then $w^{-1} \vec f_1, \dots, w^{-1} \vec f_m$ is the geodesic from $v_0$ to $w^{-1} v$, and
	\[
		w \, \delta_{\vec T, v_0} \, w^{-1} v
		= w \, o(w^{-1} \vec f_1) \dots o(w^{-1} \vec f_m) \, w^{-1} v
		= o(\vec f_1) \dots o(\vec f_m) \, v
		= \delta_{\vec T, w v_0}(v) .
	\]
\end{proof}

\begin{lemma} \label{lem:twist-conjugate-action}
	For $w \in \FF$, we have $\delta_{\vec T, v_0} \, w \, \delta_{\vec T, v_0}^{-1} = o(\vec e_1) \dots o(\vec e_k) \, w$, where $\vec e_1, \dots, \vec e_k$ is the geodesic from $v_0$ to $w v_0$.
\end{lemma}

\begin{proof}
	Let $v$ be a vertex and $\vec f_1, \dots, \vec f_m$ the geodesic from $w v_0$ to $v$.
	By Lemma~\ref{lem:twist-conjugate}, $\delta_{\vec T, v_0} \, w \, \delta_{\vec T, v_0}^{-1} = \delta_{\vec T, v_0} \, \delta_{\vec T, w v_0}^{-1} \, w$.
	The path $\vec e_1, \dots, \vec e_k, \vec f_1, \dots, \vec f_m$ runs from $v_0$ to $v$, so Lemma~\ref{lem:twist-path} gives $\delta_{\vec T, v_0}(v) = o(\vec e_1) \dots o(\vec e_k) \, o(\vec f_1) \dots o(\vec f_m) \, v = o(\vec e_1) \dots o(\vec e_k) \, \delta_{\vec T, w v_0}(v)$.
	Applying this to $\delta_{\vec T, w v_0}^{-1}(w v)$ in place of $v$ gives $\delta_{\vec T, v_0} \, \delta_{\vec T, w v_0}^{-1}(w v) = o(\vec e_1) \dots o(\vec e_k) \, w v$.
\end{proof}

The group $\FF$ acts faithfully on $T$.
Indeed, the kernel of the action fixes every edge, so it is a normal subgroup of $\FF$ contained in the infinite cyclic group $\Stab(e)$.
If it were nontrivial, with generator $a$, then every $b \in \FF$ would satisfy $b a b^{-1} = a^{\pm 1}$, so $b^2$ would commute with $a$.
By Theorem~\ref{thm:free}\eqref{item:free:commute}, the squares of two elements of a free basis would then be powers of a common element, which is false since $g \ge 2$.
So we regard $\FF$ as a subgroup of the automorphism group of $T$.

By Lemma~\ref{lem:twist-conjugate-action}, conjugation by $\delta_{\vec T, v_0}$ preserves $\FF$ in the automorphism group of $T$, so it defines an automorphism $\varphi_{v_0}$ of $\FF$.
For vertices $v_0, v_1$ with geodesic $\vec e_1, \dots, \vec e_k$ from $v_0$ to $v_1$, Lemma~\ref{lem:twist-path} gives $\delta_{\vec T, v_0} = o(\vec e_1) \dots o(\vec e_k) \, \delta_{\vec T, v_1}$, so $\varphi_{v_0}$ and $\varphi_{v_1}$ differ by an inner automorphism.
So the outer class of $\varphi_{v_0}$ does not depend on $v_0$; we call it the \emph{Dehn twist} along $\vec T$ and write $\delta_{\vec T} \in \Out(\FF)$.
Reversing the orientation replaces $\delta_{\vec T}$ by its inverse, by Lemma~\ref{lem:twist-automorphism}.
In Section~\ref{sec:dehn-twist} we show that, when $\vec T$ is the oriented dual tree of a primitive torus, $\delta_{\vec T}$ is the image in $\Out(\FF)$ of a homeomorphism of $\MM$ twisting a regular neighborhood of that torus.

\subsection{Fundamental results on three-manifolds} \label{subsec:three-manifolds}

This subsection collects the results from $3$-manifold topology that we cite by name; the reader may skip this subsection and come back when a theorem is invoked.
They concern spheres and balls in $3$-manifolds, and are used to build the model of the universal cover in Section~\ref{sec:univ}.
The results specific to the doubled handlebody are collected in Section~\ref{subsec:doubled_handlebody}.

The universal cover of the doubled handlebody is the complement of a Cantor set in $\mathbb{S}^3$, and much of the paper takes place there.
The first theorem is what makes spheres in $\mathbb{S}^3$ tractable.
Throughout, embedded curves, spheres and tori in $\MM$, and their lifts, are smooth, hence locally flat; homeomorphisms are arbitrary.

\begin{theorem}[Alexander \cite{Ale24}; see Hatcher's note \cite{HatcherNotes}] \label{thm:alexander}
	Every locally flat embedded $2$-sphere in $\mathbb{S}^3$ bounds a closed $3$-ball on each side.
\end{theorem}

We use Theorem~\ref{thm:alexander} in the proof of Lemma~\ref{lem:sphere-lift}.

The following statement about self-maps of the $2$-sphere is used only in the proof of Lemma~\ref{lem:extend}, to build the model in Section~\ref{sec:univ}.

\begin{theorem}[Smale \cite{MR112149}, Kneser \cite{MR1544816}; see Martelli's note \cite{Mar22}] \label{thm:isotopy}
	Every orientation-preserving diffeomorphism of $\mathbb{S}^2$ is isotopic to the identity through diffeomorphisms, and every orientation-preserving homeomorphism of $\mathbb{S}^2$ is isotopic to the identity through homeomorphisms.
\end{theorem}

When we say \emph{by general position in dimension $3$}, we use that a map of a circle or an arc into a $3$-manifold can be perturbed to an embedding.
The perturbation can be taken supported in a small neighborhood of the locus where the map fails to be injective.
When this locus misses a given surface to which the map is transverse, the perturbation can be taken away from the surface, so that the transverse intersection points with the surface are unchanged.

\subsection{Doubled handlebodies} \label{subsec:doubled_handlebody}

Let $\MM = \MM_g$ be the connected sum of $g$ copies of $\mathbb{S}^1 \times \mathbb{S}^2$, the double of the handlebody of genus $g$.
Its fundamental group is free of rank $g$, and we identify it with $\FF$ once and for all.
An embedded sphere in $\MM$ is \emph{essential} if it bounds no $3$-ball.
Essential tori and primitive tori are defined in Section~\ref{subsec:tori}, using the universal cover.
Throughout the paper, $\MM$ and every torus in $\MM$ are oriented, and the \emph{positive side} of a torus is the side into which its positive normal points.
For a nonseparating torus this is a local notion, a choice of one boundary component of a regular neighborhood, and it becomes a genuine side for the lifts of the torus to the universal cover.

Let $\MCG(\MM)$ be the group of isotopy classes of orientation-preserving homeomorphisms of $\MM$.
A homeomorphism induces an automorphism of $\FF$ up to inner automorphisms, which gives a homomorphism $\MCG(\MM) \to \Out(\FF)$.
Whitehead \cite{MR1575455, MR1503309} showed that it is surjective, and used it to solve the automorphic equivalence problem for $\FF$.
Laudenbach \cite{MR356056} showed that its kernel is $(\ZZ / 2 \ZZ)^g$, generated by twists along essential spheres, and Brendle, Broaddus and Putman \cite{MR4557874} showed that the resulting extension splits.
A sphere twist acts trivially on $\pi_1(\MM)$ and on $\pi_2(\MM)$ \cite[Lemma~2.1]{MR5014056}, so it preserves the homotopy class of every embedded sphere, hence its isotopy class by Theorem~\ref{thm:laudenbach}\eqref{item:laudenbach:spheres} below.
We fix a section
\[
	\Out(\FF) \hookrightarrow \MCG(\MM)
\]
once and for all, and regard outer automorphisms of $\FF$ as mapping classes of $\MM$ through this section.

The next theorem is the reason that spheres in $\MM$ can be treated up to homotopy, and the surjectivity of $\MCG(\MM) \to \Out(\FF)$ in a sharper form.
Its first item has no analogue for primitive tori, and Section~\ref{subsec:classes} therefore works with classes of tori defined through their dual trees.

\begin{theorem}[Laudenbach \cite{MR356056}] \label{thm:laudenbach}
	The following hold.
	\begin{enumerate}
		\item \label{item:laudenbach:spheres}
		Homotopic embedded $2$-spheres in $\MM$ are isotopic.
		\item \label{item:laudenbach:mcg}
		Every automorphism of $\FF$ is induced by a homeomorphism of $\MM$.
	\end{enumerate}
\end{theorem}

Maximal sphere systems cut $\MM$ into pieces with trivial fundamental group.
This gives the intersection bound of Lemma~\ref{lem:bounded} in Section~\ref{sec:curve-space}.

A \emph{sphere system} is a finite disjoint union of pairwise non-isotopic essential spheres in $\MM$; it is \emph{maximal} if no further sphere can be added.

\begin{theorem}[Laudenbach \cite{MR356056}, Hatcher \cite{MR1314940}] \label{thm:systems}
	Every sphere system is contained in a maximal one.
	If $\Sigma$ is maximal, every component of $\MM$ cut along $\Sigma$ is a holed sphere with three boundary spheres.
\end{theorem}

\subsection{Intersection numbers and primitive curves} \label{subsec:primitive}

A simple closed curve in $\MM$ is \emph{essential} if it is not null-homotopic.
Embedded submanifolds $a$ and $a'$ of $\MM$ are \emph{homotopic} if there is a homotopy from the inclusion of $a$ to a homeomorphism onto $a'$.

Let $A$ and $B$ be homotopy classes of embedded submanifolds of $\MM$, such as simple closed curves, spheres and tori.
Their \emph{geometric intersection number} is
\[
	i(A, B) := \min \bigl\{ \lvert \pi_0(a \cap b) \rvert \;\big|\; a \in A , \ b \in B \text{ transverse} \bigr\} ,
\]
the minimal number of components of $a \cap b$.
For two curves it is zero.
For a curve and a surface it is the minimal number of intersection points, and for two surfaces the minimal number of circles.

An essential simple closed curve $c$ is \emph{primitive} if $i([c], S) = 1$ for some homotopy class $S$ of essential spheres.

\begin{lemma} \label{lem:primitive-basis}
	Let $\gamma$ be the homotopy class of an essential simple closed curve, represented by $w \in \FF$.
	Then $\gamma$ is represented by a primitive curve if and only if $w$ belongs to a free basis of $\FF$.
\end{lemma}

\begin{proof}
	Suppose first that $w$ belongs to a free basis $\{ w, b_2, \dots, b_g \}$.
	Fix a free basis $\{ a_1, \dots, a_g \}$, an essential sphere $s_1$ and a simple closed curve $c_1$ representing $a_1$ that meets $s_1$ in one point.
	Let $\phi$ be the automorphism of $\FF$ with $\phi(a_1) = w$ and $\phi(a_i) = b_i$ for $i \ge 2$.
	By Theorem~\ref{thm:laudenbach}\eqref{item:laudenbach:mcg}, $\phi$ is induced by a homeomorphism $h$ of $\MM$.

	The curve $h(c_1)$ represents $\gamma$ and meets the essential sphere $h(s_1)$ in one point, so $i(\gamma, [h(s_1)]) \le 1$.
	The algebraic intersection number of $h(c_1)$ with $h(s_1)$ is $\pm 1$ and depends only on the classes, so $i(\gamma, [h(s_1)]) \ge 1$.
	Hence $\gamma$ is primitive.

	Conversely, let $c$ be a primitive curve representing $\gamma$, meeting an essential sphere $s$ in one point.
	A closed curve meets a separating sphere in an even number of points, so $s$ is nonseparating.
	Let $\MM'$ be $\MM$ cut along $s$, with the two copies $s^\pm$ of $s$ in its boundary.
	Then $\FF = \pi_1(\MM') \ast \langle a \rangle$, where $a$ is represented by a loop crossing $s$ once.
	The curve $c$ cut at its point on $s$ is an arc in $\MM'$ from $s^+$ to $s^-$, so $w$ is conjugate to $a u$ for some $u \in \pi_1(\MM')$, after replacing $a$ by $a^{-1}$ if necessary.

	As a free factor of $\FF$, the group $\pi_1(\MM')$ is free, of rank $g - 1$ by Grushko's theorem.
	Let $\{ b_2, \dots, b_g \}$ be a free basis of $\pi_1(\MM')$, so that $\{ a, b_2, \dots, b_g \}$ is a free basis of $\FF$.
	The map $a \mapsto a u$, $b_i \mapsto b_i$ is an automorphism of $\FF$, so $\{ a u, b_2, \dots, b_g \}$ is a free basis containing a conjugate of $w$.
	Hence $w$ belongs to a free basis.
\end{proof}

\section{The universal cover} \label{sec:univ}

The main goal of this section is Proposition~\ref{prop:cantor}: the universal cover of $\MM$ is the complement of a tame Cantor set $C$ in $\mathbb{S}^3$, built from a tree of pieces indexed by $\FF$.
On this model we describe the lifts of essential spheres and tori, their halfspaces, and when two lifts cross; that crossing depends only on homotopy classes is shown in Section~\ref{subsec:homotopic-lifts}.

Let $\rho: \widetilde\MM \to \MM$ be the universal cover of $\MM$.

\subsection{A model in the sphere}

This subsection constructs the pieces $P'_w$, the balls $D_w$ and the Cantor set $C$, and proves Proposition~\ref{prop:cantor}.

We regard $\mathbb{S}^3$ as the one-point compactification of $\mathbb{R}^3$.
A Cantor set $C$ in $\mathbb{S}^3$ is said to be \emph{tame} if there exists a homeomorphism of $\mathbb{S}^3$ that moves $C$ into a line segment.
We prove the following.

\begin{proposition} \label{prop:cantor}
	$\widetilde{\MM}$ is homeomorphic to $\mathbb{S}^3$ minus a tame Cantor set $C$.
	The homeomorphism can be chosen so that $\mathbb{S}^3 \setminus C$ is the union of the holed spheres $P'_w$, $w \in \FF$, constructed below, with disjoint interiors, and each deck transformation $g \in \FF$ carries $P'_w$ onto $P'_{gw}$.
\end{proposition}

We begin the proof with a sequence of lemmas.
A \emph{holed sphere} is $\mathbb{S}^3$ minus finitely many disjoint open balls bounded by round spheres.
Cutting $\MM$ along the spheres $\{\mathrm{pt}\} \times \mathbb{S}^2$ of the $g$ summands gives a manifold $P$ homeomorphic to a holed sphere with $2g$ boundary spheres, and we fix such a homeomorphism.

Write the boundary spheres of $P$ as $s_i^\pm$ for each $i = 1, \dots, g$.
Let $\varphi_i \colon s_i^+ \to s_i^-$ be the gluing maps, so that $\MM$ is obtained from $P$ by gluing $s_i^+$ to $s_i^-$ via $\varphi_i$.
Let $\FF$ be the free group with a free basis $\beta = \{ b_1, \dots, b_g \}$, and identify its elements with reduced words in $\beta^{\pm 1}$.
Let $\mathrm{Cay}(\FF) = \mathrm{Cay}(\FF, \beta)$ be its Cayley graph, a tree in which $w$ is joined to $w b_i$ by an edge for each $i$.
For each $w \in \FF$, let $P_w$ be a copy of $P$.
Let $U_0$ denote the space constructed from $\bigsqcup_w P_w$ by gluing $s_i^+ \subset P_w$ to $s_i^- \subset P_{w b_i}$ via $\varphi_i$ for every edge of $\mathrm{Cay}(\FF)$.

\begin{lemma} \label{lem:tree}
	$U_0$ is the universal cover of $\MM$.
\end{lemma}

\begin{proof}
	The copies $P_w$ project to $\MM$, and the projections agree on the glued spheres.
	The resulting map from $U_0$ to $\MM$ is a covering map.
	Indeed, a point off the cutting spheres has a neighborhood inside $P$, whose preimage is one copy in each $P_w$.
	A point on the $i$-th cutting sphere has a neighborhood made of a collar of $s_i^+$ and a collar of $s_i^-$ glued by $\varphi_i$.
	Its preimage is one such pair of collars for each edge of $\mathrm{Cay}(\FF)$ labeled $b_i$.

	We claim that $U_0$ is simply connected.
	The union of the pieces over a finite subtree, thickened by collars, is built by gluing simply connected open sets along connected simply connected intersections, one piece at a time.
	So it is simply connected by van Kampen.
	Every loop in $U_0$ lies in such a finite union.
\end{proof}

Let us say a family of closed round balls $\{ D_w \subseteq \mathbb{S}^3 \}_{w \neq 1}$ is \emph{nested} when $D_{w'} \subset \operatorname{int} D_w$ if $w$ is a proper prefix of $w'$, and $D_w \cap D_{w'} = \varnothing$ if neither is a prefix of the other.
Given a nested family, let $P'_1$ be $\mathbb{S}^3$ minus the interiors of the balls of one-letter words.
For $w \ne 1$, let $P'_w$ be $D_w$ minus the interiors of the balls of words one letter longer than $w$.
By nesting, $P'_w$ and $P'_{w'}$ intersect only if one word is the other with one more letter, and then in the boundary of the ball of the longer word.
So $P'_w \cap P'_{w'}$ is one round sphere if $w$, $w'$ are adjacent in $\mathrm{Cay}(\FF)$, and empty otherwise.
Each $P'_w$ is a holed sphere with $2g$ boundary spheres, since $\mathbb{S}^3$ contains $2g$ balls of one-letter words and $D_w$ contains $2g - 1$ balls of words one letter longer.

\begin{lemma} \label{lem:model}
	There are a line $L$ and a nested family $\{D_w\}$ such that every $D_w$ is centered on $L$ and has radius at most $2^{1 - \lvert w \rvert}$.
\end{lemma}

\begin{proof}
	Let $L$ be the $x$-axis of $\mathbb{R}^3 \subset \mathbb{S}^3$.
	Let $D_{b_1^{\pm 1}}, \dots, D_{b_g^{\pm 1}}$ be disjoint closed balls of radius $1$ centered on $L$.
	Suppose that $D_w$ is defined for all reduced words $w$ of length $n$, with center $c_w \in L$ and radius $r_w$.
	Let such $w$ end in the letter $x$, and order the $2g - 1$ letters $y \ne x^{-1}$.
	For the $j$-th letter $y$, let $D_{wy}$ be the ball of radius $r_w / 8g$ centered at the point of $L$ at signed distance $(2j - 1) r_w / 4g - r_w / 2$ from $c_w$.
	These balls are disjoint, lie in $\operatorname{int} D_w$, and have radius at most $r_w / 2$.

	So $r_w \le 2^{1 - \lvert w \rvert}$ by induction, and $D_{w'} \subset \operatorname{int} D_w$ when $w$ is a proper prefix of $w'$.
	If neither of $w$, $w'$ is a prefix of the other, their shortest differing prefixes $u$, $u'$ have equal length, so $D_u$ and $D_{u'}$ are disjoint by construction, and $D_w \subset D_u$, $D_{w'} \subset D_{u'}$.
\end{proof}

From now on, we fix a line $L$ and a nested family $\{D_w\}$ given by Lemma~\ref{lem:model}, so that every $D_w$ is centered on $L$ and has radius at most $2^{1 - \lvert w \rvert}$, and we let $C$ denote the complement in $\mathbb{S}^3$ of the union of the sets $P'_w$.

\begin{lemma} \label{lem:cantor}
	$C$ is a tame Cantor set.
\end{lemma}

\begin{proof}
	We claim that $C$ is the set of points lying in balls $D_w$ with $\lvert w \rvert$ arbitrarily large.
	A point outside every $P'_w$ lies in a one-letter ball, since it is outside $P'_1$; it lies in the interior of a child of that ball, since it is outside the corresponding $P'_w$; and so on.
	Conversely, a point lying in balls of arbitrarily long words lies in $\operatorname{int} D_{uy}$ for every $u$ and the letter $y$ with $u y$ a prefix of a longer such word, so it lies in no $P'_u$.

	So $C$ is an intersection of nonempty compact sets, hence nonempty and compact.
	Every point of $C$ lies in balls $D_w$ of arbitrarily small radius, and $D_w \cap C$ is open and closed in $C$; so $C$ is totally disconnected.
	Every $D_w$ contains at least two balls of the next generation, as $2g - 1 \ge 2$; so $C$ has no isolated points.
	The centers lie on $L$ and the radii tend to zero, so $C$ lies on $L$, hence in a segment of $L$.
\end{proof}

\begin{lemma} \label{lem:extend}
	Let $Q$ and $Q'$ be holed spheres, and $\psi$ a bijection from the boundary spheres of $Q$ to the boundary spheres of $Q'$.
	Suppose that the centers of the boundary spheres of $Q'$ lie in a plane.
	Then every homeomorphism from a boundary sphere $\Sigma$ of $Q$ to $\psi(\Sigma)$ extends to a homeomorphism from $Q$ to $Q'$ carrying each boundary sphere to its image under $\psi$.
\end{lemma}

\begin{proof}
	Let $f$ be the given map.
	Any two configurations of the same number of disjoint closed round balls in $\mathbb{S}^3$ are ambient isotopic by an isotopy realizing any given bijection between the two configurations.
	Indeed, we may shrink the balls about their centers, and then move the centers, since the configuration space of finitely many distinct points of $\mathbb{S}^3$ is connected.
	Here a round ball is a metric ball of the round metric on $\mathbb{S}^3$.
	Stereographic projection carries round spheres to round spheres, so the complementary balls of a holed sphere are round in this sense, and we shrink each about its spherical center.
	Applied to the complementary balls of $Q$ and $Q'$, this gives a homeomorphism $k$ from $Q$ to $Q'$ carrying each boundary sphere to its image under $\psi$.

	Then $k^{-1} \circ f$ is a self-homeomorphism of $\Sigma$.
	If it reverses orientation, replace $k$ by $r \circ k$, where $r$ is the reflection of $\mathbb{S}^3$ in a plane containing the centers of the boundary spheres of $Q'$; it preserves $Q'$ and each of its boundary spheres.
	Now $k^{-1} \circ f$ preserves orientation, so it is isotopic to the identity by Theorem~\ref{thm:isotopy}.
	Realize the isotopy on a collar of $\Sigma$ in $Q$, and the identity outside the collar.
	This is a self-homeomorphism $F$ of $Q$ restricting to $k^{-1} \circ f$ on $\Sigma$ and preserving each boundary sphere.
	So $k \circ F$ extends $f$.
\end{proof}

Lemma~\ref{lem:extend} applies to $Q = P_w$ through the fixed homeomorphism of $P$ with a holed sphere, and to $Q' = P'_w$, whose boundary spheres are centered on $L$.
For each $w$ and each letter $y$, let $\sigma_{w,y}$ be the boundary sphere of $P'_w$ shared with $P'_{wy}$; it is $\partial D_{wy}$ if $wy$ is longer than $w$, and $\partial D_w$ otherwise.
Let $\psi_w$ be the bijection sending $s_i^+$ to $\sigma_{w, b_i}$ and $s_i^-$ to $\sigma_{w, b_i^{-1}}$.

\begin{lemma} \label{lem:homeo}
	There are homeomorphisms $h_w \colon P_w \to P'_w$ for all $w \in \FF$, each carrying boundary spheres as $\psi_w$ does, such that $h_{w b_i} \circ \varphi_i = h_w$ on $s_i^+$ for all $w \in \FF$ and $i$.
\end{lemma}

\begin{proof}
	Let $h_1$ be any homeomorphism carrying boundary spheres as $\psi_1$ does; it exists by Lemma~\ref{lem:extend}.
	Suppose that $h_w$ is defined for all words of length at most $n$.
	Let $w'$ have length $n + 1$ and let $w$ be its prefix of length $n$.
	The edge between $w$ and $w'$ glues a boundary sphere $\Sigma$ of $P_{w'}$ to a boundary sphere of $P_w$, by $\varphi_i$ or $\varphi_i^{-1}$, and $h_w$ carries the latter onto $\partial D_{w'} = \psi_{w'}(\Sigma)$.
	So $h_w$ and the gluing map prescribe $h_{w'}$ on $\Sigma$, as a homeomorphism onto $\psi_{w'}(\Sigma)$.
	By Lemma~\ref{lem:extend}, this prescription extends to a homeomorphism $h_{w'}$ carrying boundary spheres as $\psi_{w'}$ does.

	Every edge of $\mathrm{Cay}(\FF)$ joins a word to a word one letter longer, so the required identity holds for every $w$ and $i$.
\end{proof}

\begin{proof}[Proof of Proposition~\ref{prop:cantor}]
	Let $U_1$ be the union of the sets $P'_w$, so that $U_1 = \mathbb{S}^3 \setminus C$ with $C$ a tame Cantor set by Lemma~\ref{lem:cantor}.
	By the identity of Lemma~\ref{lem:homeo}, the maps $h_w$ agree on the spheres glued in $U_0$, so they descend to a map $h$ from $U_0$ to $U_1$.
	It is a bijection, since the sets $P'_w$ cover $U_1$ and intersect only along the spheres of the edges.
	It is continuous, since $U_0$ carries the quotient topology and $h$ is continuous on each $P_w$.

	The family $\{P'_w\}$ is locally finite in $U_1$: a point has positive distance from $C$, while $P'_w$ lies within distance $2^{2 - \lvert w \rvert}$ of $C$.
	So $h^{-1}$ is continuous on each member of a locally finite closed cover of $U_1$, hence continuous.
	Therefore $h$ is a homeomorphism, and $\widetilde{\MM} \cong U_0 \cong U_1$ by Lemma~\ref{lem:tree}.

	The deck transformation $g \in \FF$ acts on $U_0$ by sending $P_w$ onto $P_{gw}$, and $h$ sends $P_w$ onto $P'_w$.
	So under $h$ it sends $P'_w$ onto $P'_{gw}$.
\end{proof}

\subsection{Essential tori} \label{subsec:tori}

This subsection characterizes essential spheres by their lifts, and defines essential tori, their cores and primitive tori.

Fix the identification $\widetilde\MM = \mathbb{S}^3 \setminus C$ constructed in the proof of Proposition~\ref{prop:cantor}, with its pieces $P'_w$ and balls $D_w$.
The set $C$ is compact and totally disconnected, and each ball $D_u$ meets $\widetilde\MM$ in the connected union of the pieces $P'_w$ over the subtree of words with prefix $u$.
So the points of $C$ are exactly the ends of $\widetilde\MM$ and $\mathbb{S}^3$ is its end compactification \cite{MR1545233}.
Hence the definitions below would not change under another identification of $\widetilde\MM$ with $\mathbb{S}^3$ minus a tame Cantor set given by Proposition~\ref{prop:cantor}.

A closed subset of $\mathbb{S}^3$ is said to \emph{split} $C$ if each component of its complement contains a point of $C$.
A $2$-sphere in $\mathbb{S}^3$ is said to be \emph{essential} if it is disjoint from $C$ and splits $C$.
A $2$-sphere in $\mathbb{S}^3$ is said to be \emph{pinched} if it intersects $C$ in exactly two points and splits $C$.

Recall that an embedded sphere in $\MM$ is essential if it bounds no ball.

\begin{lemma} \label{lem:sphere-lift}
	An embedded sphere $s$ in $\MM$ is essential if and only if a lift of $s$ is an essential $2$-sphere, if and only if every lift of $s$ is.
\end{lemma}

\begin{proof}
	Let $\tilde s$ be a lift of $s$.
	Since $s$ is simply connected, $\tilde s$ is an embedded $2$-sphere in $\widetilde\MM$, disjoint from $C$, with trivial stabilizer.
	It suffices to show that $s$ bounds a ball if and only if $\tilde s$ does not split $C$.

	Suppose first that $s$ bounds a ball $B$.
	Since $B$ is simply connected, the lift of $B$ through a point of $\tilde s$ is a ball $\tilde B \subset \widetilde\MM$ with boundary $\tilde s$.
	Its interior is a component of $\mathbb{S}^3 \setminus \tilde s$ disjoint from $C$, so $\tilde s$ does not split $C$.

	Conversely, suppose that $\tilde s$ does not split $C$.
	By Theorem~\ref{thm:alexander}, $\tilde s$ bounds a closed ball $\tilde B$ on the side missing $C$, so $\tilde B \subset \widetilde\MM$.
	We claim that $\rho$ is injective on $\tilde B$.
	Otherwise $g \tilde B \cap \tilde B \ne \varnothing$ for some $g \in \FF \setminus \{ 1 \}$.
	The spheres $g \tilde s$ and $\tilde s$ are distinct lifts of the embedded sphere $s$, so they are disjoint.
	So one of $g \tilde B$ and $\tilde B$ contains the other, or $g \tilde B \cup \tilde B = \mathbb{S}^3$.

	The last case is impossible, since both balls miss the nonempty set $C$.
	In the other cases, $g$ or $g^{-1}$ maps the ball $\tilde B$ into $\tilde B$, and has a fixed point by Brouwer's theorem, which is impossible for a nontrivial deck transformation.
	So $\rho$ embeds $\tilde B$, and $s$ bounds the ball $\rho(\tilde B)$.
\end{proof}

An embedded torus $t$ in $\MM$ is said to be \emph{essential} if the closure in $\mathbb{S}^3$ of some lift of $t$ is a pinched $2$-sphere.
This does not depend on the lift, since deck transformations extend to homeomorphisms of $\mathbb{S}^3$ preserving $C$, as shown in Section~\ref{subsec:ends}.

\begin{lemma} \label{lem:essential_pi1}
	For an essential torus $t$, the image of a map $\pi_1(t) \to \FF$ induced from the inclusion $t \hookrightarrow \MM$ is an infinite cyclic group.
\end{lemma}

\begin{proof}
	Let $H$ be the image, and let $\tilde t$ be a lift of $t$.
	The restriction $\tilde t \to t$ of the universal cover is the connected covering of $t$ corresponding to the kernel of $\pi_1(t) \to \FF$, so its deck group is $H$ and $H \cong \pi_1(t) / \pi_1(\tilde t)$.
	Since $t$ is essential, $\tilde t$ is an open annulus, so $\pi_1(\tilde t) \cong \ZZ$.
	Hence $H$ is a quotient of $\ZZ^2$ by an infinite cyclic subgroup, and in particular $H$ is abelian and infinite.
	A subgroup of a free group is free, and an abelian free group is trivial or infinite cyclic.
	So $H$ is infinite cyclic.
\end{proof}

Let $t$ be an essential torus.
By Lemma~\ref{lem:essential_pi1}, the image of $\pi_1(t) \to \FF$ is infinite cyclic; it is well-defined up to conjugation, and it depends only on the homotopy class of $t$.
Every element of $\FF$ is represented by a simple closed curve, since $\dim \MM = 3$.
So there is a homotopy class $\gamma$ of an essential simple closed curve whose representatives generate this image, and $\gamma$ is unique, since a generator of an infinite cyclic group is determined up to inversion.
We call $\gamma$ the \emph{core} of $t$; it depends only on the homotopy class of $t$.
We say that $t$ is \emph{primitive} if its core is a homotopy class of a primitive curve.

\subsection{Halfspaces and crossing} \label{subsec:lifts}

This subsection defines the halfspaces of a lift and crossing of two lifts; that both depend only on homotopy classes is shown in Section~\ref{subsec:homotopic-lifts}.

Let $u$ be an essential sphere or torus and $\tilde u$ a lift.
If $u$ is a sphere, then $\tilde u$ is a $2$-sphere with trivial stabilizer.
If $u$ is a torus, then $\tilde u$ is an open annulus, and its stabilizer $\langle y \rangle$ acts cocompactly on $\tilde u$, where $y$ represents the core of $u$ by Lemma~\ref{lem:essential_pi1}.
In both cases $\widetilde{\MM} \setminus \tilde u$ has exactly two components, since $H_1(\widetilde{\MM}; \ZZ / 2) = 0$.
We call them the \emph{halfspaces} of $\tilde u$.
For an annulus $\tilde u$, both halfspaces are $\langle y \rangle$-invariant, since $y$ preserves the orientations of $\tilde u$ and of $\widetilde{\MM}$.

Let $\tilde u$ and $\tilde v$ be lifts of essential spheres or tori.
A \emph{corner} of $(\tilde u, \tilde v)$ is an intersection $Y \cap Z_0$ of a halfspace $Y$ of $\tilde u$ with a halfspace $Z_0$ of $\tilde v$.

Fix a Riemannian metric on $\MM$ and lift it to $\widetilde{\MM}$.
Then $\FF$ acts on $\widetilde{\MM}$ by isometries, and bounded subsets of $\widetilde{\MM}$ have compact closure.
For $A \subset \widetilde{\MM}$ and $R \ge 0$, let $N_R(A)$ be the closed $R$-neighborhood of $A$.
A subset $A$ is \emph{coarsely contained} in $B$ if $A \subset N_R(B)$ for some $R \ge 0$.
Following Scott \cite{MR1608688}, we say that $\tilde u$ and $\tilde v$ \emph{cross} if no corner of $(\tilde u, \tilde v)$ is coarsely contained in $\tilde u \cup \tilde v$.

Disjoint lifts do not cross.
Indeed, if $\tilde u \cap \tilde v = \varnothing$, then $\tilde v$ lies in one halfspace of $\tilde u$.
The other halfspace $Y'$ is connected and misses $\tilde v$, so it lies in one halfspace of $\tilde v$.
Hence $Y'$ misses the other halfspace $Z'$ of $\tilde v$, and the corner $Y' \cap Z'$ is empty.

We now use the model $\widetilde{\MM} = \mathbb{S}^3 \setminus C$, with the pieces $P'_w$ and the balls $D_w$ of the identification fixed in Section~\ref{subsec:tori}, indexed by the vertices $w$ of the Cayley tree $\mathrm{Cay}(\FF)$.
By Proposition~\ref{prop:cantor}, the deck transformation $g$ carries $P'_w$ onto $P'_{gw}$.
Let $y \in \FF \setminus \{1\}$.
The element $y$ acts on $\mathrm{Cay}(\FF)$ with an axis, a line on which $y$ translates.
The vertices $y^n w$ converge, as $n \to \pm\infty$, to the two ends of the axis.
For each end, the balls $D_u$, for the prefixes $u$ of the end, decrease to a point of $C$.
Let $y^{+\infty}$ be the point for the end $\lim_{n \to +\infty} y^n w$ and $y^{-\infty}$ the point for the end $\lim_{n \to -\infty} y^n w$.
We keep this notation for Lemmas~\ref{lem:closure} and~\ref{lem:closure-limit}.

\begin{lemma} \label{lem:closure}
	Let $y \in \FF \setminus \{1\}$, let $y^{+\infty}$ and $y^{-\infty}$ be the points of $C$ for the two ends of the axis of $y$ in $\mathrm{Cay}(\FF)$, and let $A \subset \widetilde{\MM}$ be a $\langle y \rangle$-invariant subset with compact quotient.
	Then the closure of $A$ in $\mathbb{S}^3$ is $A \cup \{ y^{+\infty}, y^{-\infty} \}$.
\end{lemma}

\begin{proof}
	We first show that the limit points of $A$ in $C$ lie in $\{ y^{+\infty}, y^{-\infty} \}$.
	The set $A$ is closed in $\widetilde\MM$: its image in the covering space $\widetilde\MM / \langle y \rangle$ of $\MM$, a Hausdorff space, is compact and hence closed, and $A$ is the full preimage of its image.
	So the closure of $A$ in $\mathbb{S}^3$ is $A$ together with its limit points in $C$.
	Let $W$ be the set of $w$ with $P'_w \cap A \ne \varnothing$.
	Finitely many pieces meet a compact set.
	So $W$ consists of finitely many $\langle y \rangle$-orbits, and lies within distance $d$ of the axis for some $d$.

	Let $a_n \in A$ leave every compact subset of $\widetilde{\MM}$, and choose $w_n \in W$ with $a_n \in P'_{w_n}$.
	Finitely many pieces have words of bounded length, so $\lvert w_n \rvert \to \infty$.
	Let $v_n$ be a vertex of the axis within distance $d$ of $w_n$.
	Then $\lvert v_n \rvert \to \infty$, and after passing to a subsequence, $v_n$ converges to one end of the axis, say the end of $y^{+\infty}$.
	Two words at distance at most $d$ share a prefix of length at least the smaller length minus $d$.
	Since $v_n$ converges to the end of $y^{+\infty}$, the common prefix of $v_n$ and that end has length tending to infinity.

	So $w_n$, $v_n$ and the end share a prefix $u_n$ of length tending to infinity.
	Hence $a_n \in P'_{w_n} \subset D_{u_n}$, and $D_{u_n}$ converges to $y^{+\infty}$.
	So every limit point of $A$ in $C$ is $y^{+\infty}$ or $y^{-\infty}$.

	It remains to show that both points are limit points of $A$.
	Let $a \in A$ lie in the piece $P'_w$, and fix a vertex $v$ of the axis.
	Then $y^n a \in A$ lies in $P'_{y^n w}$, and $y^n w$ is within distance $d(w, v)$ of the vertex $y^n v$ of the axis.
	The vertices $y^n v$ converge to the end of $y^{+\infty}$ as $n \to +\infty$, so the common prefix of $y^n v$ and that end has length tending to infinity.
	Two words at distance at most $d(w, v)$ share a prefix of length at least the smaller length minus $d(w, v)$, so $y^n w$ and that end share a prefix $u_n$ of length tending to infinity.
	Hence $y^n a \in D_{u_n}$ and $y^n a \to y^{+\infty}$ as $n \to +\infty$.
	Likewise $y^n a \to y^{-\infty}$ as $n \to -\infty$.
\end{proof}

\begin{lemma} \label{lem:closure-limit}
	Let $y \in \FF \setminus \{1\}$, let $y^{+\infty}$ and $y^{-\infty}$ be the points of $C$ for the two ends of the axis of $y$ in $\mathrm{Cay}(\FF)$, and let $K \subset \widetilde{\MM}$ be compact.
	Then the sets $y^n K$ converge to $y^{\pm\infty}$ as $n \to \pm\infty$: every neighborhood of $y^{+\infty}$ contains $y^n K$ for all large $n$, and likewise for $y^{-\infty}$ and $n \to -\infty$.
\end{lemma}

\begin{proof}
	Let $a_n \in y^{m_n} K$ with $m_n \to +\infty$.
	Fix a vertex $v$ of the axis.
	Finitely many pieces meet $K$, so $a_n$ lies in a piece $P'_{w_n}$ with $w_n$ within bounded distance of $y^{m_n} v$.
	The vertices $y^{m_n} v$ lie on the axis and converge to the end of $y^{+\infty}$.
	Two words at distance at most $d$ share a prefix of length at least the smaller length minus $d$.
	Since $w_n$ lies within bounded distance of these vertices, $w_n$ and that end share a prefix $u_n$ of length tending to infinity, so $a_n \in D_{u_n}$ and $a_n \to y^{+\infty}$.
	The case $m_n \to -\infty$ is the same.
\end{proof}

Let $\tilde u$ be a lift of an essential torus with stabilizer $\langle y \rangle$.
By Lemma~\ref{lem:closure}, the closure of $\tilde u$ in $\mathbb{S}^3$ is the pinched sphere $\hat u := \tilde u \cup \{ y^{\pm\infty} \}$.
We call $y^{+\infty}$ and $y^{-\infty}$ the \emph{ends} of $\tilde u$.
For a lift $\tilde u$ of an essential sphere, we put $\hat u := \tilde u$.

\begin{lemma} \label{lem:axes}
	Let $A, B \subset \widetilde{\MM}$ be invariant under $\langle a \rangle$ and $\langle b \rangle$ respectively, with compact quotients.
	If $A \cap N_R(B)$ is unbounded for some $R$, then $a$ and $b$ have a common nonzero power.
\end{lemma}

\begin{proof}
	Choose compact sets $K_A$, $K_B$ with $A \subset \langle a \rangle K_A$ and $B \subset \langle b \rangle K_B$.
	There are integers $n_k$, taking infinitely many values, and $m_k$ such that $a^{n_k} K_A$ and $b^{m_k} K_B$ are within distance $R$.
	Since $\FF$ acts properly discontinuously, the elements $a^{-n_k} b^{m_k}$ lie in a finite set.
	So $a^{-n_k} b^{m_k} = a^{-n_l} b^{m_l}$ for some $k$, $l$ with $n_k \ne n_l$, that is, $a^{n_l - n_k} = b^{m_l - m_k}$.
	This element is nontrivial since $\FF$ is torsion-free.
\end{proof}

\begin{lemma} \label{lem:halfspace-unbounded}
	No halfspace of a lift of an essential sphere or torus is coarsely contained in the lift.
\end{lemma}

\begin{proof}
	Let $\tilde u$ be a lift of an essential sphere or torus $u$.
	The set $\hat u$ splits $C$, by Lemma~\ref{lem:sphere-lift} if $u$ is a sphere and by definition if $u$ is a torus.
	Suppose that a halfspace $Y$ of $\tilde u$ satisfies $Y \subset N_R(\tilde u)$.
	The closure of $N_R(\tilde u)$ in $\mathbb{S}^3$ meets $C$ only in points of $\hat u$.
	Indeed, if $u$ is a torus, it meets $C$ only in the ends of $\tilde u$ by Lemma~\ref{lem:closure}, and if $u$ is a sphere, then $N_R(\tilde u)$ is compact and misses $C$.

	Let $W$ be the component of $\mathbb{S}^3 \setminus \hat u$ containing $Y$.
	Since $C$ is closed and totally disconnected, $W \setminus C$ is connected: a closed subset of dimension $0$ does not separate a connected manifold of dimension at least $2$ \cite[Chapter~IV, Corollary~1, p.~48]{MR6493}.
	It contains $Y$ and misses $\tilde u$, so $W \setminus C = Y$.
	Hence $W \cap C$ lies in the closure of $Y$, hence in $\hat u$.
	But $W \cap \hat u = \varnothing$.
	So $W$ contains no point of $C$, contradicting that $\hat u$ splits $C$.
\end{proof}

\begin{lemma} \label{lem:boundary}
	Let $V \subset \widetilde{\MM}$ be a connected closed $3$-dimensional submanifold whose boundary is a locally finite disjoint union of lifts $\sigma$ of essential spheres or tori.
	For each $\sigma$, let $U_\sigma$ be the closure of the halfspace of $\sigma$ disjoint from $\operatorname{int} V$.
	Then $\widetilde{\MM} = V \cup \bigcup_\sigma U_\sigma$, where $U_\sigma \cap V = \sigma$ and the sets $U_\sigma$ are pairwise disjoint.
\end{lemma}

\begin{proof}
	The connected set $V \setminus \sigma$ lies in the halfspace of $\sigma$ containing $\operatorname{int} V$, and so do the other boundary surfaces.
	Let $W$ be a component of $\widetilde{\MM} \setminus V$; its closure intersects some boundary surface $\sigma$.
	Near $\sigma$, the set $W$ lies in the other halfspace of $\sigma$, and $W$ misses $\sigma$, so $W \subset U_\sigma \setminus \sigma$.
	Since $U_\sigma \setminus \sigma$ is connected and misses $V$, it equals $W$.
	In particular $\bar W$ intersects no other boundary surface.
\end{proof}

\section{Ends, sides and lifts of curves} \label{sec:ends}

The goal of this section is Proposition~\ref{prop:orbit-count}: the intersection number of a curve with an essential sphere or a primitive torus is the number of orbits of lifts of the surface penetrated by a lift of the curve.
To state it, we record a lift of a curve by its two ends in the Cantor set $C$, and a lift of a surface by the two sides into which it cuts $C$.
We also define penetration and the dual tree of a surface.
On the way we show, in Lemma~\ref{lem:representatives}, that the ends and sides of the lifts of a torus depend only on its homotopy class, and we express crossing through sides in Lemma~\ref{lem:crossing-criterion}.

\subsection{Ends and sides} \label{subsec:ends}

This subsection collects the facts about the Cantor set $C$ that the rest of the paper uses.
They concern how sequences converge to points of $C$, how $\FF$ acts on $C$, and how a lift of a sphere or torus cuts $C$ into two sides.

We use the pieces $P'_w$ and the balls $D_w$ of the identification fixed in Section~\ref{subsec:tori}.
A point $c \in C$ is the intersection of the balls $D_u$ over the prefixes $u$ of an infinite reduced word, called the \emph{word} of $c$.
These balls form a neighborhood basis of $c$, since their radii tend to $0$ by Lemma~\ref{lem:model}.
A piece $P'_w$ lies in $D_u$ if $u$ is a prefix of $w$, and misses the interior of $D_u$ otherwise.
Through words, we identify $C$ with the set of ends of $\mathrm{Cay}(\FF)$.
Since $D_u \cap C$ is the set of points whose word has the prefix $u$, the topology of $C$ is the topology of the ends.

For a subset $Z$ of $\widetilde\MM$, we write $C \cap Z$ for the set of points of $C$ in the closure of $Z$ in $\mathbb{S}^3$, and call it the \emph{Cantor part} of $Z$.
Since $Z$ itself misses $C$, this causes no ambiguity; for a subset of $\mathbb{S}^3$ that may meet $C$, such as an open subset $W$ of $\mathbb{S}^3$ or a side defined below, $C \cap W$ is the literal intersection.

\begin{lemma} \label{lem:prefix}
	Let $x_n \in P'_{w_n}$ and $c \in C$.
	Then $x_n \to c$ in $\mathbb{S}^3$ if and only if the length of the common prefix of $w_n$ and the word of $c$ tends to infinity.
\end{lemma}

\begin{proof}
	Suppose that the common prefix has length tending to infinity.
	Then $x_n$ lies in balls $D_{u_n}$ around $c$ with $\lvert u_n \rvert \to \infty$, so $x_n \to c$.
	Conversely, suppose that $x_n \to c$, and let $u$ be a prefix of the word of $c$.
	Then $x_n \in \operatorname{int} D_u$ for all large $n$, so $u$ is a prefix of $w_n$.
\end{proof}

\begin{lemma} \label{lem:bounded-limit}
	Let $x_n, z_n \in \widetilde\MM$ with $d(x_n, z_n)$ bounded.
	If $x_n \to c \in C$ in $\mathbb{S}^3$, then $z_n \to c$.
\end{lemma}

\begin{proof}
	Let $w_n$ and $w'_n$ be words with $x_n \in P'_{w_n}$ and $z_n \in P'_{w'_n}$.
	Finitely many pieces meet a bounded neighborhood of a piece, and the action is by isometries.
	So the distance of $w_n$ and $w'_n$ in $\mathrm{Cay}(\FF)$ is bounded by some $d$.
	By Lemma~\ref{lem:prefix}, the words $w_n$ share with the word of $c$ prefixes of length tending to infinity.
	By the bound $d$, so do the words $w'_n$, and $z_n \to c$ by Lemma~\ref{lem:prefix}.
\end{proof}

We now extend equivariant homeomorphisms of $\widetilde\MM$ to $\mathbb{S}^3$.
An automorphism $\varphi$ of $\FF$ is a bi-Lipschitz map of the vertex set of $\mathrm{Cay}(\FF)$, since $\lvert \varphi(w) \rvert \le L \lvert w \rvert$ and $\lvert \varphi^{-1}(w) \rvert \le L \lvert w \rvert$ for some $L$.
By the following theorem, it extends to a homeomorphism $\partial\varphi$ of the ends of $\mathrm{Cay}(\FF)$: if vertices $w_n$ converge to an end $\xi$, then $\varphi(w_n)$ converge to $\partial\varphi(\xi)$.

\begin{theorem}[{\cite[Theorem III.H.3.9]{MR1744486}}] \label{thm:qi}
	A quasi-isometry between proper geodesic Gromov-hyperbolic spaces extends to a homeomorphism of their Gromov boundaries.
\end{theorem}

\begin{lemma} \label{lem:qi-extends}
	Let $f$ be a homeomorphism of $\widetilde\MM$ and $\varphi$ an automorphism of $\FF$ with $f(gx) = \varphi(g) f(x)$ for all $g \in \FF$ and $x \in \widetilde\MM$.
	Then $f$ extends to a homeomorphism of $\mathbb{S}^3$ that acts on $C$ as $\partial\varphi$.
\end{lemma}

\begin{proof}
	Let $\bar f$ be the map equal to $f$ on $\widetilde\MM$ and to $\partial\varphi$ on $C$.
	Fix $x_0 \in P'_1$, and let $x_n \in P'_{w_n}$ converge to a point $c \in C$ with word $\xi$.
	Then $w_n^{-1} x_n \in P'_1$ by Proposition~\ref{prop:cantor}, so $f(x_n) = \varphi(w_n) f(w_n^{-1} x_n)$ lies in $\varphi(w_n) f(P'_1)$.
	Since $f(P'_1)$ is compact and $\FF$ acts by isometries, $f(x_n)$ lies within bounded distance of $\varphi(w_n) x_0 \in P'_{\varphi(w_n)}$.
	By Lemma~\ref{lem:prefix}, $w_n \to \xi$, so $\varphi(w_n) \to \partial\varphi(\xi)$, and $\varphi(w_n) x_0 \to \partial\varphi(c)$ by Lemma~\ref{lem:prefix} again.
	So $f(x_n) \to \partial\varphi(c)$ by Lemma~\ref{lem:bounded-limit}.

	Since $\partial\varphi$ is continuous on $C$ and $\widetilde\MM$ is open, $\bar f$ is continuous.
	The same holds for $f^{-1}$ and $\varphi^{-1}$, and $\partial(\varphi^{-1}) = (\partial\varphi)^{-1}$, so $\bar f$ is a homeomorphism.
\end{proof}

\begin{lemma} \label{lem:qi-extends-cantor}
	Let $f$ be a homeomorphism of $\widetilde\MM$ and $\varphi$ an automorphism of $\FF$ with $f(gx) = \varphi(g) f(x)$ for all $g \in \FF$ and $x \in \widetilde\MM$, and write $f$ also for the extension of $f$ to $\mathbb{S}^3$ given by Lemma~\ref{lem:qi-extends}, which is unique since $\widetilde\MM$ is dense in $\mathbb{S}^3$.
	Then $f(C \cap A) = C \cap f(A)$ for every $A \subset \widetilde\MM$.
\end{lemma}

\begin{proof}
	The extension is a homeomorphism of $\mathbb{S}^3$ preserving $C$, so it carries the closure of $A$ onto the closure of $f(A)$, and $C \cap \overline{A}$ onto $C \cap \overline{f(A)}$.
\end{proof}

\begin{lemma} \label{lem:extend-isotopy}
	Let $\Xi$ be an ambient isotopy of $\MM$ and $\tilde\Xi$ its lift with $\tilde\Xi_0$ the identity.
	Then each $\tilde\Xi_s$ extends by the identity on $C$ to a homeomorphism of $\mathbb{S}^3$.
\end{lemma}

\begin{proof}
	The lift $\tilde\Xi$ is $\FF$-equivariant, since $\tilde\Xi_0$ is the identity and the lifts of $\Xi_s \circ \rho$ through a given point are unique.
	So Lemma~\ref{lem:qi-extends} applies with $\varphi$ the identity, for which $\partial\varphi$ is the identity of $C$.
\end{proof}

Every deck transformation $h$ satisfies $h(gx) = (h g h^{-1}) h(x)$.
So by Lemma~\ref{lem:qi-extends}, with $\varphi$ the conjugation by $h$, it extends to a homeomorphism of $\mathbb{S}^3$.
On $C$ it sends an end $\xi$ to $h \xi$, since $h w_n h^{-1}$ is within distance $\lvert h \rvert$ of $h w_n$ in $\mathrm{Cay}(\FF)$.
Thus the action of $\FF$ on $C$ is its action on the ends of $\mathrm{Cay}(\FF)$, and the ends $h^{\pm\infty}$ of $h \in \FF \setminus \{ 1 \}$ are the ends of its axis.

\begin{lemma} \label{lem:end-in-side}
	Let $g \in \FF \setminus \{1\}$, $x \in \widetilde\MM$ and $B \subset \widetilde\MM$ with $g^n x \in B$ for all large $n$.
	Then $g^{+\infty} \in C \cap B$.
\end{lemma}

\begin{proof}
	This is Lemma~\ref{lem:closure-limit} for $K = \{ x \}$.
\end{proof}

\begin{lemma} \label{lem:common-end}
	Let $h, h' \in \FF \setminus \{1\}$ have a common end.
	Then $h$ and $h'$ have a common nonzero power.
\end{lemma}

\begin{proof}
	The ends of $h$ are the ends of its axis in $\mathrm{Cay}(\FF)$, and likewise for $h'$.
	Two axes with a common end share a ray.
	Let $v$ be a vertex on both axes and $x \in P'_v$.
	The orbits $\langle h \rangle x$ and $\langle h' \rangle x$ have compact quotients, and they lie within bounded distance of each other along that ray, since the pieces of the vertices of the ray are at bounded distance from each other.
	So Lemma~\ref{lem:axes} gives a common nonzero power.
\end{proof}

\begin{lemma} \label{lem:end-stabilizer}
	Let $y$ be an element of a free basis of $\FF$, and let $g \in \FF$ carry an end of $y$ to an end of $y$.
	Then $g \in \langle y \rangle$.
\end{lemma}

\begin{proof}
	The element $g y g^{-1}$ shares an end with $y$, so $g y g^{-1} = y^k$ for some $k$ by Lemma~\ref{lem:common-end} and Theorem~\ref{thm:free}\eqref{item:free:primitive}.
	Conjugation preserves translation length on $\mathrm{Cay}(\FF)$, so $k = \pm 1$.
	An element of a free group is not conjugate to its inverse, so $g y g^{-1} = y$ and $g \in \langle y \rangle$.
\end{proof}

\begin{lemma} \label{lem:north-south}
	Let $h \in \FF \setminus \{1\}$ and $K \subset C \setminus \{ h^{-\infty} \}$ compact.
	Then $h^n K \to h^{+\infty}$ uniformly as $n \to +\infty$.
\end{lemma}

\begin{proof}
	The action of $\FF$ on $C$ is its action on the ends of $\mathrm{Cay}(\FF)$, with the topology given by the balls $D_u$ by Lemma~\ref{lem:prefix}.
	For the action of $h$ on the ends of its axis, this is Theorem~\ref{thm:trees}\eqref{item:trees:ends}.
\end{proof}

Let $u$ be an essential sphere or torus, $\tilde u$ a lift with halfspaces $Y$ and $Y^\ast$, and $\hat u$ the closure of $\tilde u$ in $\mathbb{S}^3$.
The ends of $\tilde u$ are the two points of $\hat u \cap C$ if $u$ is a torus, and there are none if $u$ is a sphere.

\begin{lemma} \label{lem:sides}
	$\mathbb{S}^3 \setminus \hat u$ has exactly two components $W$ and $W^\ast$, with $W \setminus C = Y$ and $W^\ast \setminus C = Y^\ast$.
\end{lemma}

\begin{proof}
	Let $W$ be the component containing $Y$.
	Since $C$ is closed and totally disconnected, $W \setminus C$ is connected: a closed subset of dimension $0$ does not separate a connected manifold of dimension at least $2$ \cite[Chapter~IV, Corollary~1, p.~48]{MR6493}.
	It contains $Y$ and misses $\tilde u$, and $Y$ is a component of $\widetilde\MM \setminus \tilde u$, so $W \setminus C = Y$.
	The same holds for the component $W^\ast$ containing $Y^\ast$, and $W \ne W^\ast$ since $Y \ne Y^\ast$.
	Every component of $\mathbb{S}^3 \setminus \hat u$ contains a point outside $C$, since $C$ has empty interior, and this point lies in $Y$ or in $Y^\ast$.
	So there are no other components.
\end{proof}

We keep the notation $W$ and $W^\ast$ of Lemma~\ref{lem:sides} for the next two lemmas.

\begin{lemma} \label{lem:sides-nonempty}
	Let $W$ and $W^\ast$ be the two components of $\mathbb{S}^3 \setminus \hat u$, where $\hat u$ is the closure in $\mathbb{S}^3$ of a lift $\tilde u$ of an essential sphere or torus $u$.
	Then both $W$ and $W^\ast$ contain points of $C$.
\end{lemma}

\begin{proof}
	The set $\hat u$ splits $C$, by Lemma~\ref{lem:sphere-lift} if $u$ is a sphere and by definition if $u$ is a torus.
	By Lemma~\ref{lem:sides}, $W$ and $W^\ast$ are the components of $\mathbb{S}^3 \setminus \hat u$, so each contains a point of $C$.
\end{proof}

\begin{lemma} \label{lem:sides-partition}
	Let $W$ and $W^\ast$ be the two components of $\mathbb{S}^3 \setminus \hat u$, where $\hat u$ is the closure in $\mathbb{S}^3$ of a lift $\tilde u$ of an essential sphere or torus $u$.
	Then the sets $C \cap W$ and $C \cap W^\ast$ are open in $C$, and $C$ minus the ends of $\tilde u$ is their disjoint union.
\end{lemma}

\begin{proof}
	The sets $W$ and $W^\ast$ are open and disjoint, so $C \cap W$ and $C \cap W^\ast$ are open in $C$ and disjoint.
	By Lemma~\ref{lem:sides}, their union is $C \setminus \hat u$, and $\hat u$ meets $C$ exactly in the ends of $\tilde u$.
\end{proof}

We call $C(Y) := C \cap W$ and $C(Y^\ast) := C \cap W^\ast$ the two \emph{sides} of $\tilde u$.
By Lemma~\ref{lem:sides-partition}, every point of $C$ other than the ends of $\tilde u$ lies on exactly one side of $\tilde u$.
Since $W \subset \overline{Y} \subset W \cup \hat u$, the Cantor part of $Y$ is $C \cap Y = C(Y) \cup \{ \text{ends of } \tilde u \}$; for a sphere, $C(Y) = C \cap Y$.

\begin{lemma} \label{lem:ends-side}
	Let $\tilde u$ be a lift of a primitive torus with stabilizer $\langle y \rangle$, and let $P$ be a side of $\tilde u$.
	Let $\tilde v$ be a lift of an essential sphere or torus with a side $Z$, such that $\tilde u$ and $\tilde v$ have no common end and $y^{+\infty} \in Z$.
	Then $Z \cap P$ contains a point that is an end neither of $\tilde u$ nor of $\tilde v$.
\end{lemma}

\begin{proof}
	The side $P$ is a nonempty open subset of $C$ by Lemmas~\ref{lem:sides-nonempty} and~\ref{lem:sides-partition}, hence infinite, so it contains a point $\xi$ that is an end of neither lift.
	Since $y$ preserves the halfspace of $P$, it preserves $P$, so $y^k \xi \in P$ for all $k$.
	By Lemma~\ref{lem:north-south}, $y^k \xi \to y^{+\infty}$ as $k \to +\infty$.
	Since $Z$ is open in $C$ and contains $y^{+\infty}$, we have $y^k \xi \in Z \cap P$ for all large $k$.
	The points $y^k \xi$ are distinct, they are not ends of $\tilde u$, and at most two of these points are ends of $\tilde v$.
\end{proof}

The following two lemmas, Lemmas~\ref{lem:core-no-penetrate} and~\ref{lem:orbits-bound} below, and Lemma~\ref{lem:far-sides} of Section~\ref{subsec:crossing-sides} are used only from Section~\ref{sec:dehn-twist} on.

\begin{lemma} \label{lem:far-halfspace}
	Let $\tilde u_n$ be lifts of essential spheres or tori leaving every compact subset of $\widetilde\MM$, let $x_0 \in \widetilde\MM$, and let $Q_n$ be the halfspace of $\tilde u_n$ not containing $x_0$.
	Then the diameters of the sets $Q_n \cup \tilde u_n$ in $\mathbb{S}^3$ tend to $0$.
\end{lemma}

\begin{proof}
	Let $x_0 \in P'_v$.
	The deck transformation $v^{-1}$ extends to a homeomorphism of $\mathbb{S}^3$, whose inverse is uniformly continuous, and it carries $x_0$ into $P'_1$ and the lifts $\tilde u_n$ to lifts leaving every compact set.
	So we may assume that $x_0 \in P'_1$, and it suffices to find words $u_n$ with $\lvert u_n \rvert \to \infty$ and $Q_n \cup \tilde u_n \subset D_{u_n}$, since the radii of these balls tend to $0$ by Lemma~\ref{lem:model}.

	Let $W_n$ be the set of words $w$ with $P'_w \cap \tilde u_n \ne \varnothing$.
	Finitely many pieces have words of bounded length, so $R_n := \min \{ \lvert w \rvert : w \in W_n \}$ tends to infinity.
	Since $\tilde u_n$ is connected and the pieces are closed, $W_n$ is connected in $\mathrm{Cay}(\FF)$.
	So it lies in one component of the complement of the ball of radius $R_n - 1$ about $1$, and all its words share a prefix $u_n$ of length $R_n$.
	Hence $\tilde u_n \subset D_{u_n}$.

	Let $w$ be the word of a piece meeting $Q_n$, and let $w = w_0, w_1, \dots, w_k = 1$ be the prefixes of $w$ in decreasing length.
	The union of the pieces $P'_{w_i}$ is connected and joins $Q_n$ to $x_0$, so it meets $\tilde u_n$; thus some $w_i$ lies in $W_n$ and has the prefix $u_n$, and so does $w$.
	Hence $Q_n \subset D_{u_n}$.
\end{proof}

\begin{lemma} \label{lem:far-halfspace-limit}
	Let $\tilde u_n$ be lifts of essential spheres or tori leaving every compact subset of $\widetilde\MM$, let $x_0 \in \widetilde\MM$, and let $Q_n$ be the halfspace of $\tilde u_n$ not containing $x_0$.
	Then, after passing to a subsequence, the sets $C \cap (Q_n \cup \tilde u_n)$ converge to a single point of $C$.
\end{lemma}

\begin{proof}
	The closure of $Q_n$ meets $C$, as $Q_n$ is unbounded by Lemma~\ref{lem:halfspace-unbounded}, so $C \cap (Q_n \cup \tilde u_n)$ is nonempty.
	Since $\mathbb{S}^3$ is compact, after passing to a subsequence a point of $C \cap (Q_n \cup \tilde u_n)$ converges to a point $\zeta \in \mathbb{S}^3$, and $\zeta \in C$ since $C$ is closed.
	By Lemma~\ref{lem:far-halfspace}, the diameters of the sets $C \cap (Q_n \cup \tilde u_n)$ tend to $0$, so these sets converge to $\zeta$.
\end{proof}

\subsection{Lifts of homotopic tori} \label{subsec:homotopic-lifts}

We use the orientations and positive sides of tori fixed in Section~\ref{subsec:doubled_handlebody}.
Let $u$ and $u'$ be homotopic essential tori, and $G$ a homotopy from the inclusion of $u$ to a homeomorphism $G_1 \colon u \to u'$.
We may and do assume, after possibly reversing the orientation of $u'$, that $G_1$ preserves orientation; the conclusions of this subsection are symmetric in the two sides.

Let $\tilde G$ be a lift of $G$, from a lift $\tilde u$ of $u$ to a lift $\tilde u'$ of $u'$, and let $\langle y \rangle$ be the stabilizer of $\tilde u$.
Since the lifts of $G \circ (\rho \times \mathrm{id})$ through a given point are unique, $\tilde G$ is $\langle y \rangle$-equivariant, so $\langle y \rangle$ is contained in the stabilizer of $\tilde u'$.

Conversely, the stabilizer of a lift of a connected subset of $\MM$ is the image of its fundamental group in $\FF$.
Here the base point is transported along a path whose lift from the base point of $\widetilde\MM$ ends in the given lift.
Transport the base point of $u$ to $\tilde u$ along a path $\beta$ ending at $p \in u$.
Transport that of $u'$ to $\tilde u'$ along $\beta$ followed by the track $G(p, \cdot)$; the lift of this track from the end point of the lift of $\beta$ ends in $\tilde u'$.
Since $G$ is a homotopy from the inclusion of $u$ to $G_1$, the transported image of $\pi_1(u', G_1(p)) = (G_1)_\ast \pi_1(u, p)$ equals the transported image of $\pi_1(u, p)$.
So both lifts have the stabilizer $\langle y \rangle$; this uses only that $G$ is a homotopy ending in a homeomorphism, not that $u'$ is essential.

Since $\tilde G$ is equivariant with compact quotient, it is proper, and its image lies in $N_r(\tilde u)$ for some $r \ge 0$.
Let $Y$ and $Y'$ be the halfspaces of $\tilde u$ and $\tilde u'$ on their positive sides, and $Y^\ast$ and $Y'^\ast$ the other halfspaces.

\begin{lemma} \label{lem:same-outside}
	$Y \setminus N_r(\tilde u) = Y' \setminus N_r(\tilde u)$ and $Y^\ast \setminus N_r(\tilde u) = Y'^\ast \setminus N_r(\tilde u)$.
\end{lemma}

\begin{proof}
	Since $\tilde G$ is proper and $G_1$ preserves orientation, $\tilde G(\tilde u \times [0,1])$ is a locally finite $3$-chain with boundary $\tilde u' - \tilde u$.
	Hence a path with endpoints outside $N_r(\tilde u)$ has the same algebraic intersection number with $\tilde u$ and with $\tilde u'$.
	For a path with endpoints off $\tilde u$, this number is the indicator of the endpoint lying in $Y$ minus that of the starting point; likewise for $\tilde u'$ and $Y'$.
	By Lemma~\ref{lem:halfspace-unbounded}, there are points $x \in Y$ and $x' \in Y^\ast$ outside $N_r(\tilde u)$.
	A path from $x'$ to $x$ has intersection number $1$ with $\tilde u$, hence with $\tilde u'$; so $x' \in Y'^\ast$ and $x \in Y'$.
	For $z \notin N_r(\tilde u)$, a path from $x'$ to $z$ then shows that $z \in Y$ if and only if $z \in Y'$.
\end{proof}

\begin{lemma} \label{lem:same-cantor}
	$C \cap Y = C \cap Y'$ and $C \cap Y^\ast = C \cap Y'^\ast$.
\end{lemma}

\begin{proof}
	By Lemma~\ref{lem:closure}, the Cantor part of $N_r(\tilde u)$ is $\{ y^{\pm\infty} \}$.
	These two points lie in all four sets, since they lie in the closures of $\tilde u$ and of $\tilde u'$.

	So $C \cap Y \setminus \{ y^{\pm\infty} \}$ is the set of limit points in $C$ of $Y \setminus N_r(\tilde u)$ other than $y^{\pm\infty}$.
	Indeed, a point of $C \cap Y$ other than $y^{\pm\infty}$ lies in the closure of $Y$ but not in the closure of $N_r(\tilde u)$, so it is a limit of points of $Y \setminus N_r(\tilde u)$.
	Conversely, every limit point of $Y \setminus N_r(\tilde u)$ lies in the closure of $Y$.
	By Lemma~\ref{lem:same-outside}, the limit points of $Y \setminus N_r(\tilde u)$ are those of $Y' \setminus N_r(\tilde u)$, so $C \cap Y \setminus \{ y^{\pm\infty} \} = C \cap Y' \setminus \{ y^{\pm\infty} \}$.
	The same holds for $Y^\ast$ and $Y'^\ast$.
\end{proof}

\begin{lemma} \label{lem:invariant-sphere}
	Let $u$ be an essential torus and $u'$ a homotopic torus, not assumed essential, and let $G$ be a homotopy from the inclusion of $u$ to a homeomorphism $G_1 \colon u \to u'$.
	Let $\tilde u'$ be the lift of $u'$ at which the lift of $G$ from a lift $\tilde u$ of $u$ ends, so that $\tilde u$ and $\tilde u'$ have the same stabilizer $\langle y \rangle$.
	Then $\tilde u'$ is an open annulus, and its closure $\hat u'$ in $\mathbb{S}^3$ is $\tilde u' \cup \{ y^{\pm\infty} \}$, an embedded $2$-sphere meeting $C$ exactly in $y^{+\infty}$ and $y^{-\infty}$.
\end{lemma}

\begin{proof}
	The lift $\tilde u'$ is an open annulus since its stabilizer is infinite cyclic.
	There is a compact $K \subset \tilde u'$ with $\tilde u' = \langle y \rangle K$.
	A sequence tending to one end of the annulus $\tilde u'$ lies in sets $y^{m_n} K$ with $m_n \to +\infty$, say, and a sequence tending to the other in sets $y^{m_n} K$ with $m_n \to -\infty$.
	So the two ends of the annulus $\tilde u'$ converge to $y^{+\infty}$ and $y^{-\infty}$ respectively, by Lemma~\ref{lem:closure-limit}, and $\hat u'$ is the two-point compactification of $\tilde u'$, an embedded $2$-sphere meeting $C$ in $y^{\pm\infty}$.
\end{proof}

\begin{lemma} \label{lem:invariant}
	Whether a torus is essential depends only on its homotopy class.
\end{lemma}

\begin{proof}
	Let $u$ be an essential torus and $u'$ a homotopic torus, and let $G$ be a homotopy from the inclusion of $u$ to a homeomorphism $G_1 \colon u \to u'$.
	Let $\tilde u'$ be the lift of $u'$ at which the lift of $G$ from a lift $\tilde u$ of $u$ ends, and let $\langle y \rangle$ be the stabilizer of $\tilde u$, which is also the stabilizer of $\tilde u'$ by the beginning of this subsection; the lifts are open annuli since the stabilizers agree.
	Let $Y$ and $Y'$ be the halfspaces of $\tilde u$ and $\tilde u'$ on their positive sides, and $Y^\ast$ and $Y'^\ast$ the other halfspaces.
	Since $H_1(\widetilde{\MM}; \ZZ/2) = 0$, the properly embedded annulus $\tilde u'$ has two complementary components, which is all that the proofs of Lemmas~\ref{lem:same-outside} and~\ref{lem:same-cantor} use about $u'$.
	By Lemma~\ref{lem:invariant-sphere}, the closure $\hat u'$ of $\tilde u'$ in $\mathbb{S}^3$ is an embedded $2$-sphere meeting $C$ in $y^{\pm\infty}$.
	We show that $\hat u$ splits $C$ if and only if $\hat u'$ does.
	The proof of Lemma~\ref{lem:sides} uses only these properties of $\tilde u'$ and $\hat u'$, not that $u'$ is essential.
	So $\mathbb{S}^3 \setminus \hat u'$ has two components $W'$ and $W'^\ast$ with $W' \setminus C = Y'$ and $W'^\ast \setminus C = Y'^\ast$.
	Moreover $\hat u'$ splits $C$ if and only if $C \cap W'$ and $C \cap W'^\ast$ are both nonempty.
	Since $Y'$ is dense in $W'$ and $\overline{Y'} \subset W' \cup \hat u'$, the set $C \cap W'$ is the Cantor part of $Y'$ minus $y^{\pm\infty}$, and likewise for $W'^\ast$; the same holds for $\tilde u$ by Lemma~\ref{lem:sides}.
	By Lemma~\ref{lem:same-cantor}, the Cantor parts of the halfspaces of $\tilde u'$ are those of the halfspaces of $\tilde u$, so $\hat u$ splits $C$ if and only if $\hat u'$ does.
\end{proof}

\begin{lemma} \label{lem:invariant-crossing}
	Let $u$ and $u'$ be homotopic essential tori, let $G$ be a homotopy from the inclusion of $u$ to a homeomorphism $G_1 \colon u \to u'$, and let $\tilde u'$ be the lift of $u'$ at which the lift of $G$ from a lift $\tilde u$ of $u$ ends.
	Let $\tilde v$ be a lift of an essential sphere or torus.
	Then $\tilde u$ crosses $\tilde v$ if and only if $\tilde u'$ does.
\end{lemma}

\begin{proof}
	By Lemma~\ref{lem:same-outside}, the halfspaces of $\tilde u$ and of $\tilde u'$ agree outside $N_r(\tilde u)$, and $\tilde u \subset N_r(\tilde u')$.
	A set has points arbitrarily far from $\tilde u \cup \tilde v$ exactly when it has points arbitrarily far from $\tilde u' \cup \tilde v$.
	So the corners of $(\tilde u, \tilde v)$ are coarsely contained in $\tilde u \cup \tilde v$ exactly when those of $(\tilde u', \tilde v)$ are coarsely contained in $\tilde u' \cup \tilde v$.
\end{proof}

\begin{lemma} \label{lem:invariant-crossing-both}
	Let $u$ and $u'$ be homotopic essential tori, let $G$ be a homotopy from the inclusion of $u$ to a homeomorphism $G_1 \colon u \to u'$, and let $\tilde u'$ be the lift of $u'$ at which the lift of $G$ from a lift $\tilde u$ of $u$ ends.
	Let $v$ and $v'$ be homotopic essential spheres or homotopic essential tori, let $G'$ be a homotopy from the inclusion of $v$ to a homeomorphism $v \to v'$, and let $\tilde v'$ be the lift of $v'$ at which the lift of $G'$ from a lift $\tilde v$ of $v$ ends.
	Then $\tilde u$ crosses $\tilde v$ if and only if $\tilde u'$ crosses $\tilde v'$.
\end{lemma}

\begin{proof}
	By Lemma~\ref{lem:invariant-crossing}, $\tilde u$ crosses $\tilde v$ if and only if $\tilde u'$ crosses $\tilde v$.
	When $v$ is a torus, the argument of Lemma~\ref{lem:invariant-crossing} with the roles of the two lifts exchanged shows that $\tilde u'$ crosses $\tilde v$ if and only if $\tilde u'$ crosses $\tilde v'$.
	When $v$ is a sphere, the lift $\tilde v$ and the lifted homotopy are compact, and the Cantor part of a bounded neighborhood of $\tilde v$ is empty.
	So Lemmas~\ref{lem:same-outside} and~\ref{lem:same-cantor} hold for $\tilde v$ and $\tilde v'$ with the same proofs, and the argument of Lemma~\ref{lem:invariant-crossing} applies again.
\end{proof}

In terms of sides, Lemma~\ref{lem:same-cantor} says that the sides of a lift depend only on the homotopy class of the torus and on the lift up to bounded distance.
We make this precise.
For each lift $\tilde u$ of $u$, the lift of the homotopy $G$ from $\tilde u$ ends at a lift of $u'$; we call the resulting map from the lifts of $u$ to the lifts of $u'$ the \emph{map given by $G$}.
It is $\FF$-equivariant, since for $h \in \FF$ the lift of $G$ from $h \tilde u$ is $h \tilde G$, which ends at $h \tilde u'$.

\begin{lemma} \label{lem:representatives}
	Let $u$ and $u'$ be homotopic essential tori, and let $G$ be a homotopy from the inclusion of $u$ to a homeomorphism $G_1 \colon u \to u'$.
	Then the map given by $G$, which carries a lift $\tilde u$ of $u$ to the lift of $u'$ at which the lift of $G$ from $\tilde u$ ends, is an $\FF$-equivariant bijection from the lifts of $u$ to the lifts of $u'$, carrying each lift $\tilde u$ to a lift $\tilde u'$ with the same stabilizer, the same ends and the same sides.
\end{lemma}

\begin{proof}
	The map is $\FF$-equivariant, since for $h \in \FF$ the lift of $G$ from $h \tilde u$ is the translate by $h$ of the lift of $G$ from $\tilde u$.
	It carries $\tilde u$ to a lift $\tilde u'$ with the same stabilizer $\langle y \rangle$: the lift of $G$ is $\langle y \rangle$-equivariant, and both stabilizers are the image of $\pi_1(u)$ transported along $G$, by the beginning of this subsection.
	Since $\FF$ acts transitively on the lifts of $u$ and on the lifts of $u'$, with the stabilizer $\langle y \rangle$ at $\tilde u$ and at $\tilde u'$, the map is a bijection.
	The ends are determined by the stabilizer.
	The sides agree by Lemma~\ref{lem:same-cantor}.
\end{proof}

\begin{lemma} \label{lem:representatives-crossing}
	Let $u$ and $u'$ be homotopic essential tori, and let $G$ be a homotopy from the inclusion of $u$ to a homeomorphism $G_1 \colon u \to u'$.
	A lift of $u$ crosses a lift $\tilde v$ of an essential sphere or torus if and only if its image under the map given by $G$ does.
\end{lemma}

\begin{proof}
	The image of a lift $\tilde u$ is the lift $\tilde u'$ obtained from $\tilde u$ by lifting $G$, so this is Lemma~\ref{lem:invariant-crossing}.
\end{proof}

\subsection{Crossing in terms of sides} \label{subsec:crossing-sides}

Let $\tilde u$ and $\tilde v$ be lifts of essential spheres or tori $u$ and $v$ with no common end, with halfspaces $Y$, $Y^\ast$ and $Z_0$, $Z_0^\ast$ and sides $P := C(Y)$, $P^\ast := C(Y^\ast)$ and $Z := C(Z_0)$, $Z^\ast := C(Z_0^\ast)$.
We say that a point of $C$ is \emph{an end} if it is an end of $\tilde u$ or of $\tilde v$.

\begin{lemma} \label{lem:corner-ends-bounded}
	Let $\tilde u$ and $\tilde v$ be lifts of essential spheres or tori with no common end, and suppose that $v$ is a torus.
	Then $\tilde u \cap N_d(\tilde v)$ is bounded for every $d \ge 0$.
\end{lemma}

\begin{proof}
	If $u$ is a sphere, then $\tilde u$ is compact.
	If $u$ is a torus, let $\langle y \rangle$ and $\langle w \rangle$ be the stabilizers of $\tilde u$ and $\tilde v$.
	Then $y$ and $w$ have no common nonzero power: by Theorem~\ref{thm:trees}\eqref{item:trees:power}, such a power would have the same axis in $\mathrm{Cay}(\FF)$ as $y$ and as $w$, and the two lifts would have the same ends.
	So $\tilde u \cap N_d(\tilde v)$ is bounded by Lemma~\ref{lem:axes}.
\end{proof}

\begin{lemma} \label{lem:corner-ends}
	Suppose that $\tilde u$ and $\tilde v$ cross.
	Then the Cantor part of each of the four corners of $(\tilde u, \tilde v)$ contains a point that is not an end.
\end{lemma}

\begin{proof}
	Let $R := Y \cap Z_0$; the other corners are symmetric.
	The corner $R$ is not coarsely contained in $\tilde u \cup \tilde v$, so it is unbounded and $C \cap R \ne \varnothing$.
	Suppose, for a contradiction, that $C \cap R$ consists of ends of $\tilde u$ and $\tilde v$.
	By symmetry we may assume $w^{-\infty} \in C \cap R$, where $\langle w \rangle$ is the stabilizer of $\tilde v$; so $v$ is a torus.
	By Lemma~\ref{lem:corner-ends-bounded}, $\tilde u \cap N_d(\tilde v)$ is bounded for every $d \ge 0$.

	We first find an unbounded end $E$ of $\tilde v$ lying in $Y$.
	The covering $\tilde v \to v$ with deck group $\langle w \rangle$ is an annular cover of the torus $v$, so there is a homeomorphism $\tilde v \cong \mathbb{S}^1 \times \mathbb{R}$ under which $w$ acts as $(\theta, s) \mapsto (\theta, s + 1)$.
	Let $K := \mathbb{S}^1 \times [a, b] \subset \tilde v$ be a compact sub-annulus containing the bounded set $\tilde u \cap \tilde v$.
	Then $\tilde v \setminus K$ has the two components $\mathbb{S}^1 \times (-\infty, a)$ and $\mathbb{S}^1 \times (b, \infty)$, both unbounded, and the first, which we call $E$, contains $w^{-n} p$ for all large $n$ whenever it contains $p$, with $w^{-n} p \to w^{-\infty}$ by Lemma~\ref{lem:closure-limit}.
	The connected set $E$ misses $\tilde u$, so it lies in $Y$ or in $Y^\ast$, and $w^{-\infty}$ lies in its Cantor part.
	In fact $E \subset Y$: otherwise $w^{-\infty}$ would lie in $C \cap Y^\ast$ and in $C \cap Y$, hence be an end of $\tilde u$.

	Next we find a sphere lift $\sigma$ in $Z_0$ with a halfspace $U \subset Z_0$.
	Fix an essential sphere $s_0$ in $\MM$; its lifts are compact, of bounded diameter, and every point of $\widetilde\MM$ lies within bounded distance of one of these lifts.
	By Lemma~\ref{lem:halfspace-unbounded}, $Z_0$ contains points arbitrarily far from $\tilde v$, so $Z_0$ contains a lift $\sigma$ of $s_0$ disjoint from $\tilde v$.
	The connected set $\tilde v$ lies in one halfspace of $\sigma$; let $U$ be the closed halfspace of $\sigma$ not containing $\tilde v$.
	It is connected, meets $Z_0$ and misses the boundary $\tilde v$ of $Z_0$, so $U \subset Z_0$.
	Then $w^{-n} U \subset Z_0$ for all $n$, and $\tilde v$ misses $w^{-n} U$.

	We now show that $\tilde u$ misses $w^{-n} U$ for large $n$.
	Choose $p \in E$ and a path $\beta$ from $p$ to $\sigma$.
	The compact set $w^{-n}(\sigma \cup \beta)$ lies in $N_d(\tilde v)$ for a fixed $d$ and converges to $w^{-\infty}$, so it leaves every bounded set.
	Since $\tilde u \cap N_d(\tilde v)$ is bounded by Lemma~\ref{lem:corner-ends-bounded}, $w^{-n}(\sigma \cup \beta)$ misses $\tilde u$ for large $n$.
	Crossing lifts intersect, since disjoint lifts do not cross by Section~\ref{subsec:lifts}.
	So $\tilde u$ is connected, meets $\tilde v$, and misses the sphere $w^{-n} \sigma$; as $\tilde v$ lies outside $w^{-n} U$, so does $\tilde u$.

	Hence $w^{-n} U \cup w^{-n} \beta$ is a connected set missing $\tilde u$ and containing $w^{-n} p \in Y$, so $w^{-n} U \subset Y \cap Z_0 = R$.
	By Lemma~\ref{lem:halfspace-unbounded}, $w^{-n} U$ is unbounded, so $C \cap w^{-n} U \ne \varnothing$.
	It misses the ends of $\tilde u$ and $\tilde v$, since the other halfspace of $w^{-n} \sigma$ contains $\tilde u \cup \tilde v$ and the Cantor part of a sphere lift is empty.
	This contradicts the assumption that $C \cap R$ consists of ends.
\end{proof}

\begin{lemma} \label{lem:crossing-criterion}
	Let $\tilde u$ and $\tilde v$ be lifts of essential spheres or tori with no common end, let $P$ and $P^\ast$ be the sides of $\tilde u$, and let $Z$ and $Z^\ast$ be the sides of $\tilde v$.
	Then $\tilde u$ and $\tilde v$ cross if and only if each of the four sets $P \cap Z$, $P \cap Z^\ast$, $P^\ast \cap Z$, $P^\ast \cap Z^\ast$ contains a point that is an end neither of $\tilde u$ nor of $\tilde v$.
\end{lemma}

\begin{proof}
	Let $Y$, $Y^\ast$ be the halfspaces of $\tilde u$ and $Z_0$, $Z_0^\ast$ the halfspaces of $\tilde v$, named so that $P = C(Y)$, $P^\ast = C(Y^\ast)$, $Z = C(Z_0)$ and $Z^\ast = C(Z_0^\ast)$, and say that a point of $C$ is an end if it is an end of $\tilde u$ or of $\tilde v$.
	Suppose that $\tilde u$ and $\tilde v$ cross.
	By Lemma~\ref{lem:corner-ends}, the closure of the corner $Y \cap Z_0$ contains a point $\xi \in C$ that is not an end.
	Then $\xi \in C \cap Y$ and $\xi \in C \cap Z_0$, so $\xi \in P \cap Z$.
	The other three sets are treated in the same way.

	Suppose that they do not cross, say $Y \cap Z_0 \subset N_r(\tilde u \cup \tilde v)$.
	By Lemma~\ref{lem:closure}, and since the $r$-neighborhood of a lift of a sphere is compact, the closure of $N_r(\tilde u \cup \tilde v)$ in $\mathbb{S}^3$ meets $C$ only in ends.
	Let $\xi \in P \cap Z$.
	By Lemma~\ref{lem:sides}, there are a component $W$ of $\mathbb{S}^3 \setminus \hat u$ with $W \setminus C = Y$ and a component $W_Z$ of $\mathbb{S}^3 \setminus \hat v$ with $W_Z \setminus C = Z_0$, and $\xi \in W \cap W_Z$ since $P = C \cap W$ and $Z = C \cap W_Z$.
	Since $C$ has empty interior, the set $(W \cap W_Z) \setminus C = Y \cap Z_0$ is dense in the open set $W \cap W_Z$.
	So $\xi$ lies in the closure of $Y \cap Z_0$, and $\xi$ is an end.
\end{proof}

\begin{lemma} \label{lem:far-sides}
	Let $u$ be a primitive torus and $\langle y \rangle$ the stabilizer of $\tilde u$.
	Suppose that $\tilde u$ and $\tilde v$ do not cross, that $y^{+\infty} \in Z$, and that $P^\ast \cap Z^\ast$ contains a point that is not an end.
	Then every point of $P$ lies in $Z$ or is an end of $\tilde v$.
\end{lemma}

\begin{proof}
	By Lemma~\ref{lem:ends-side}, $Z \cap P$ and $Z \cap P^\ast$ contain points that are not ends.
	By hypothesis, so does $Z^\ast \cap P^\ast$.
	By Lemma~\ref{lem:crossing-criterion}, every point of $Z^\ast \cap P$ is an end.
	A point of $P$ is not an end of $\tilde u$, and if it is not an end of $\tilde v$, then it lies in $Z$ or in $Z^\ast$, hence in $Z$.
\end{proof}

\subsection{Lifts of curves} \label{subsec:curve-lifts}

This subsection defines when a lift of a curve penetrates a lift of a sphere or torus, and relates the penetrated lifts to the dual tree and to the intersection number.
A lift of a curve is recorded by its two ends in $C$.

Let $\gamma$ be the homotopy class of an essential simple closed curve, and $g \in \FF$ an element representing $\gamma$.
A \emph{lift} of $\gamma$ is a pair $\tilde\gamma = \{ g^{+\infty}, g^{-\infty} \}$ for such a $g$; its elements are the \emph{ends} of $\tilde\gamma$.
By Lemma~\ref{lem:common-end} and Theorem~\ref{thm:free}\eqref{item:free:commute}, the elements of $\FF$ with the same ends as $g$ form the maximal cyclic subgroup containing $g$.
So $\tilde\gamma$ determines $g$ up to inversion, since $g$ is required to represent $\gamma$; we say that $g$ \emph{represents} $\tilde\gamma$.

The lifts of $\gamma$ are the translates $h \tilde\gamma = \{ h g^{\pm\infty} \}$, $h \in \FF$.
By Lemma~\ref{lem:closure}, the ends of $\tilde\gamma$ are the two limit points in $C$ of any lift of a representative curve with stabilizer $\langle g \rangle$.
Thus a lift of the class $\gamma$ is a pair of points of $C$, while a lift of a representative curve $c$ is a line in $\widetilde\MM$.
When we speak of the lift of $c$ with stabilizer $\langle g \rangle$, we mean the latter; Lemma~\ref{lem:closure} relates the two.

Throughout this subsection, let $u$ be an essential sphere or a primitive torus, and $\tilde u$ a lift with halfspaces $Y$, $Y^\ast$.
The statements below that are cited from later sections repeat this hypothesis, since primitivity is used through Lemma~\ref{lem:end-stabilizer}.
If $u$ is a torus, let $\langle y \rangle$ be the stabilizer of $\tilde u$; by Lemma~\ref{lem:primitive-basis}, $y$ belongs to a free basis.

We say that $\tilde\gamma$ is a \emph{core} of $\tilde u$ if $u$ is a torus and the ends of $\tilde\gamma$ are the ends of $\tilde u$, that is, if $g \in \langle y \rangle$.
Among the cores of $\tilde u$, the lift of the core of $u$ is called \emph{the core} of $\tilde u$; its ends are $y^{\pm\infty}$.
By Lemma~\ref{lem:common-end}, if $\tilde\gamma$ is not a core of $\tilde u$, then no end of $\tilde\gamma$ is an end of $\tilde u$.
We say that $\tilde\gamma$ \emph{penetrates} $\tilde u$ if one end of $\tilde\gamma$ lies in $C(Y)$ and the other in $C(Y^\ast)$.
A core of $\tilde u$ does not penetrate $\tilde u$.

The \emph{dual tree} $X_u$ of $u$ is the tree whose vertices are the components of $\widetilde\MM \setminus \rho^{-1}(u)$, which we call the \emph{complementary regions} of $u$, or simply its \emph{regions}, and whose edges are the lifts of $u$, the edge $\tilde u$ joining the two regions adjacent to $\tilde u$.
It is a tree since $\widetilde\MM$ is simply connected and each lift separates it.
The action of $\FF$ on $\widetilde\MM$ induces a simplicial action on $X_u$ without edge inversions, transitive on edges, and the stabilizer of the edge $\tilde u$ is the stabilizer of $\tilde u$.
There are no edge inversions: a lift of a sphere has trivial stabilizer, and the stabilizer of a lift of a torus preserves its two halfspaces.
The latter holds by Section~\ref{subsec:lifts}, since deck transformations preserve the orientations of $\widetilde\MM$ and of the lift.
So $X_u$ is a one-edge $\ZZ$-splitting if $u$ is a torus, and a tree with one orbit of edges and trivial edge stabilizers if $u$ is a sphere.

We write $\tilde u_e$ for the lift at the edge $e$.
For a component $\kappa$ of $X_u \setminus e$, the regions in $\kappa$ lie in one halfspace of $\tilde u_e$, which we denote by $Y_\kappa$.
The two components of $X_u \setminus e$ give the two halfspaces.
The stabilizer of $e$ is the stabilizer of $\tilde u_e$, so $\tilde\gamma$ is a core of $\tilde u_e$ if and only if $g$ fixes $e$.

\begin{lemma} \label{lem:nested-sides}
	Let $u$ be an essential sphere or a primitive torus, and let $e \ne e'$ be edges of $X_u$, let $\kappa$ be the component of $X_u \setminus e$ containing $e'$, and let $\kappa'$ be the component of $X_u \setminus e'$ not containing $e$.
	Then the ends of $\tilde u_{e'}$ lie in $C(Y_\kappa)$.
\end{lemma}

\begin{proof}
	The lift $\tilde u_{e'}$ lies in $Y_\kappa$, so its ends lie in $C \cap Y_\kappa$, and they are not ends of $\tilde u_e$ by Lemma~\ref{lem:end-stabilizer}.
\end{proof}

\begin{lemma} \label{lem:nested-sides-inclusion}
	Let $u$ be an essential sphere or a primitive torus, let $e \ne e'$ be edges of $X_u$, let $\kappa$ be the component of $X_u \setminus e$ containing $e'$, and let $\kappa'$ be the component of $X_u \setminus e'$ not containing $e$.
	Then $C(Y_{\kappa'}) \subset C(Y_\kappa)$.
\end{lemma}

\begin{proof}
	Since $\kappa' \subset \kappa$, we have $Y_{\kappa'} \subset Y_\kappa$, so $C \cap Y_{\kappa'} \subset C \cap Y_\kappa$.
	By Lemma~\ref{lem:nested-sides} with the roles of $e$ and $e'$ exchanged, the ends of $\tilde u_e$ lie in the side of $\tilde u_{e'}$ other than $C(Y_{\kappa'})$.
	So a point of $C(Y_{\kappa'})$ is an end neither of $\tilde u_{e'}$ nor of $\tilde u_e$, and it lies in $C(Y_\kappa)$.
\end{proof}

Suppose that $\tilde\gamma$ is a core of no lift of $u$, that is, $g$ fixes no edge of $X_u$.
By Theorem~\ref{thm:trees}\eqref{item:trees:dichotomy}, either $g$ fixes exactly one vertex of $X_u$, or $g$ is hyperbolic with an axis $A$ on which it translates by $\lVert g \rVert > 0$.
In both cases $X_u$ has a unique minimal nonempty $\langle g \rangle$-invariant subtree: the fixed vertex, respectively the axis.
We call this subtree the \emph{minimal subtree} of $\tilde\gamma$ on $X_u$.
We orient the axis in the direction of translation by $g$.

\begin{lemma} \label{lem:subtree-sides}
	Let $u$ be an essential sphere or a primitive torus, let $\tilde\gamma$ be a core of no lift of $u$, and let $e$ be an edge of $X_u$ and $A$ the minimal subtree of $\tilde\gamma$.
	If $e \notin A$, then both ends of $\tilde\gamma$ lie in $C(Y_\kappa)$, where $\kappa$ is the component of $X_u \setminus e$ containing $A$.
	If $e \in A$, then $g^{+\infty} \in C(Y_{\kappa^+})$ and $g^{-\infty} \in C(Y_{\kappa^-})$, where $\kappa^\pm$ is the component of $X_u \setminus e$ containing the positive, respectively negative, ray of $A$.
\end{lemma}

\begin{proof}
	Let $v$ be a vertex of $A$ and $x$ a point of the region $v$.
	Then $g^n x$ lies in the region $g^n v$.
	If $e \notin A$, then $g^n v \in A \subset \kappa$ for all $n$, so $g^n x \in Y_\kappa$ and $g^{\pm\infty} \in C \cap Y_\kappa$ by Lemma~\ref{lem:end-in-side}.
	If $e \in A$, then $g^n v \in \kappa^+$ for large $n$ and $g^{n} v \in \kappa^-$ for large $-n$, so $g^{+\infty} \in C \cap Y_{\kappa^+}$ and $g^{-\infty} \in C \cap Y_{\kappa^-}$ by the same lemma.
	In both cases the ends of $\tilde\gamma$ are not ends of $\tilde u_e$, since $\tilde\gamma$ is not a core of $\tilde u_e$, so they lie in the sides.
\end{proof}

\begin{lemma} \label{lem:subtree-sides-penetrate}
	Let $u$ be an essential sphere or a primitive torus, let $\tilde\gamma$ be a core of no lift of $u$, let $e$ be an edge of $X_u$, and let $A$ be the minimal subtree of $\tilde\gamma$ on $X_u$.
	Then $\tilde\gamma$ penetrates $\tilde u_e$ if and only if $e \in A$.
\end{lemma}

\begin{proof}
	If $e \notin A$, then both ends of $\tilde\gamma$ lie in the side $C(Y_\kappa)$ of $\tilde u_e$ by Lemma~\ref{lem:subtree-sides}, so $\tilde\gamma$ does not penetrate $\tilde u_e$.
	If $e \in A$, then $\kappa^+ \ne \kappa^-$, so $C(Y_{\kappa^+})$ and $C(Y_{\kappa^-})$ are the two sides of $\tilde u_e$, and $\tilde\gamma$ penetrates $\tilde u_e$ by Lemma~\ref{lem:subtree-sides}.
\end{proof}

The proof of Proposition~\ref{prop:orbit-count} realizes the minimal intersection with a fixed representative of $U$ by a curve that forms no bigon with this representative.
It then shows that this minimum is $i(\gamma, U)$, using Laudenbach's theorem for spheres and Lemma~\ref{lem:representatives} for tori.
Let $c$ be an essential simple closed curve and $s$ an essential sphere or torus, transverse to $c$.
We say that $c$ and $s$ \emph{form a bigon} if there are distinct points $p_1, p_2 \in c \cap s$, an oriented arc $a \subset c$ from $p_1$ to $p_2$ and a path $a'$ in $s$ from $p_2$ to $p_1$ such that the loop $a \cdot a'$ is null-homotopic in $\MM$.
The arc $a$ may contain further points of $c \cap s$, and $a'$ need not be embedded.

\begin{lemma} \label{lem:curve-bigon-fewer}
	Let $c$ be an essential simple closed curve and $s$ an essential sphere or torus transverse to $c$.
	If $c$ and $s$ form a bigon, then there is a simple closed curve $c'$ homotopic to $c$ and transverse to $s$ with $\lvert c' \cap s \rvert < \lvert c \cap s \rvert$.
\end{lemma}

\begin{proof}
	Since $c$ and $s$ form a bigon, there are distinct points $p_1, p_2 \in c \cap s$, an oriented arc $a \subset c$ from $p_1$ to $p_2$ and a path $a'$ in $s$ from $p_2$ to $p_1$ such that the loop $a \cdot a'$ is null-homotopic in $\MM$.
	Since $a \cdot a'$ is null-homotopic, $a$ is homotopic relative to its endpoints to the reverse of $a'$.
	Replace $a$ in $c$ by a copy of the reverse of $a'$ pushed slightly off $s$: near $p_1$ to the side from which $c \setminus a$ arrives, and near $p_2$ to the side into which $c \setminus a$ leaves.
	If these sides differ, the pushed copy crosses $s$ once; otherwise it misses $s$.
	By general position in dimension $3$, the result is a simple closed curve $c'$ homotopic to $c$ and transverse to $s$.
	The replacement removes the points of $a \cap s$, of which there are at least two, and creates at most one, so $\lvert c' \cap s \rvert < \lvert c \cap s \rvert$.
\end{proof}

\begin{proposition}[Bigon criterion] \label{prop:curve-bigon}
	Let $c$ be an essential simple closed curve and $s$ an essential sphere or torus transverse to $c$.
	Then $c$ and $s$ form no bigon if and only if $\lvert c \cap s \rvert \le \lvert c' \cap s \rvert$ for every simple closed curve $c'$ homotopic to $c$ and transverse to $s$.
\end{proposition}

The proof is adapted from Farb and Margalit \cite{MR2850125}.
Only the implication from minimality to the absence of bigons is used, in Proposition~\ref{prop:orbit-count}.
Since a bigon is defined through a null-homotopy in $\MM$ rather than an embedded disk, the argument uses only the homotopy between the curves and general position, and no input specific to surfaces such as the Jordan curve theorem.

\begin{proof}
	If $c$ and $s$ form a bigon, then $\lvert c \cap s \rvert$ is not minimal, by Lemma~\ref{lem:curve-bigon-fewer}.

	Conversely, let $c'$ be a simple closed curve homotopic to $c$, transverse to $s$, with $\lvert c' \cap s \rvert < \lvert c \cap s \rvert$.
	Let $A := \mathbb{S}^1 \times [0,1]$ with boundary circles $\partial_0 A$ and $\partial_1 A$, and let $h \colon A \to \MM$ be a homotopy from $c$ to $c'$, a homeomorphism on each boundary circle.
	After a homotopy relative to $\partial A$, $h$ is transverse to $s$, and $h^{-1}(s)$ is a compact $1$-submanifold of $A$ consisting of circles and of properly embedded arcs.
	For $j \le k$ in $\{ 0, 1 \}$, let $n_{jk}$ be the number of arcs with endpoints on $\partial_j A$ and $\partial_k A$.
	Then $\lvert c \cap s \rvert = 2 n_{00} + n_{01}$ and $\lvert c' \cap s \rvert = 2 n_{11} + n_{01}$, so $n_{00} > n_{11} \ge 0$.
	Hence some arc $\alpha$ of $h^{-1}(s)$ has both endpoints on $\partial_0 A$.

	Such an arc gives a bigon.
	The arc $\alpha$ cuts off a disk $D \subset A$ with $\partial D = \alpha \cup \alpha_0$ for an arc $\alpha_0 \subset \partial_0 A$.
	Put $a := h(\alpha_0)$, an arc of $c$ whose endpoints are distinct points of $c \cap s$ since $h$ is injective on $\partial_0 A$, and let $a'$ be the path $h|_\alpha$ in $s$, oriented so that $a \cdot a'$ is a loop.
	Then $h|_D$ is a null-homotopy of $a \cdot a'$, so $c$ and $s$ form a bigon.
\end{proof}

\begin{lemma} \label{lem:orbit-count-reduced}
	Let $u$ be an essential sphere or a primitive torus, and suppose that $\tilde\gamma$ is a core of no lift of $u$.
	Let $c$ be a simple closed curve representing $\gamma$, transverse to $u$ and forming no bigon with $u$, and let $\tilde c$ be the lift of $c$ with stabilizer $\langle g \rangle$.
	Following $\tilde c$ through the regions gives a bi-infinite edge path in $X_u$, in which each point of $\tilde c \cap \rho^{-1}(u)$ contributes one edge, and this path is reduced.
\end{lemma}

\begin{proof}
	If the path is not reduced, then $\tilde c$ meets some lift $\tilde u_e$ in two consecutive points $p_1$, $p_2$ of $\tilde c \cap \rho^{-1}(u)$.
	The arc of $\tilde c$ from $p_1$ to $p_2$ and a path in $\tilde u_e$ from $p_2$ to $p_1$ form a loop, which is null-homotopic since $\widetilde\MM$ is simply connected.
	Its projection is a bigon formed by $c$ and $u$, provided $\rho(p_1) \ne \rho(p_2)$.
	If $\rho(p_1) = \rho(p_2)$, then $p_2 = g^k p_1$ for some $k \ne 0$, so $g^k$ stabilizes $\tilde u_e$ and $\tilde\gamma$ is a core of $\tilde u_e$, contrary to the hypothesis.
	So $c$ and $u$ form a bigon, contrary to the hypothesis on $c$.
\end{proof}

\begin{proposition} \label{prop:orbit-count}
	Let $u$ be an essential sphere or a primitive torus, let $U$ be the homotopy class of $u$, and suppose that $\tilde\gamma$ is a core of no lift of $u$.
	Then the lifts of $u$ penetrated by $\tilde\gamma$ form exactly $i(\gamma, U)$ orbits under $\langle g \rangle$.
\end{proposition}

\begin{proof}
	Let $c$ be a simple closed curve representing $\gamma$, transverse to $u$, with $\lvert c \cap u \rvert \le \lvert c' \cap u \rvert$ for every simple closed curve $c'$ representing $\gamma$ and transverse to $u$.
	By Proposition~\ref{prop:curve-bigon}, $c$ and $u$ form no bigon.
	Let $\tilde c$ be the lift of $c$ with stabilizer $\langle g \rangle$.
	By Lemma~\ref{lem:orbit-count-reduced}, following $\tilde c$ through the regions gives a reduced bi-infinite edge path in $X_u$, in which each point of $\tilde c \cap \rho^{-1}(u)$ contributes one edge.

	We first show that the number of $\langle g \rangle$-orbits of penetrated lifts is $\lvert c \cap u \rvert$.
	If $c \cap u = \varnothing$, the path is a vertex fixed by $g$, and no lift of $u$ is penetrated by Lemma~\ref{lem:subtree-sides-penetrate}.
	Otherwise the path is a $g$-invariant line, hence the axis $A$ of $g$, and $\tilde c$ meets each lift on $A$ exactly once.
	By Lemma~\ref{lem:subtree-sides-penetrate}, the penetrated lifts are the edges of $A$, so they correspond to the points of $\tilde c \cap \rho^{-1}(u)$.
	Since $\tilde c \to c$ is the quotient by $\langle g \rangle$, their $\langle g \rangle$-orbits correspond to the points of $c \cap u$.
	So in both cases, the number of $\langle g \rangle$-orbits of penetrated lifts is $\lvert c \cap u \rvert$, the minimum of $\lvert c' \cap u \rvert$ over simple closed curves $c'$ representing $\gamma$ and transverse to $u$.

	It remains to show that this minimum is $i(\gamma, U)$, which is also minimized over the representatives of $U$.
	Let $u'$ be another representative of $U$.
	If $u$ is a torus, then by Lemma~\ref{lem:representatives}, the lifts of homotopic representatives of $U$ correspond equivariantly with the same sides, so the number of penetrated orbits does not depend on the representative $u$.
	If $u$ is a sphere, then $u'$ is isotopic to $u$ by Theorem~\ref{thm:laudenbach}\eqref{item:laudenbach:spheres}, and the lift of the isotopy carries the lifts of $u$ equivariantly to the lifts of $u'$ with the same sides, by Lemma~\ref{lem:extend-isotopy}.
	In both cases the hypothesis on $\tilde\gamma$ also holds for $u'$, since corresponding lifts have the same stabilizers.
	Hence the number of penetrated orbits, and the minimum of $\lvert c' \cap u \rvert$ over the simple closed curves $c'$ representing $\gamma$ and transverse to $u$, are the same for $u$ and $u'$, and both equal $i(\gamma, U)$.
\end{proof}

\begin{proposition} \label{prop:orbit-count-translation}
	Let $u$ be an essential sphere or a primitive torus, let $U$ be the homotopy class of $u$, and suppose that $\tilde\gamma$ is a core of no lift of $u$.
	Then $i(\gamma, U)$ is the translation length of $g$ on $X_u$.
\end{proposition}

\begin{proof}
	By Lemma~\ref{lem:subtree-sides-penetrate}, the lifts of $u$ penetrated by $\tilde\gamma$ are the edges of the minimal subtree of $\tilde\gamma$ on $X_u$.
	If $g$ fixes a vertex, the minimal subtree is that vertex, no lift is penetrated, and the translation length is $0$.
	If $g$ is hyperbolic, the minimal subtree is its axis, and the number of $\langle g \rangle$-orbits of edges of the axis is the translation length.
	In both cases the number of $\langle g \rangle$-orbits of penetrated lifts is the translation length, and it is $i(\gamma, U)$ by Proposition~\ref{prop:orbit-count}.
\end{proof}

\begin{lemma} \label{lem:core-no-penetrate}
	Let $u$ be a primitive torus.
	A core of a lift of $u$ penetrates no lift of $u$.
\end{lemma}

\begin{proof}
	Let $\tilde\gamma$ be a core of $\tilde u$.
	Its ends are the ends of $\tilde u$, which lie on no side of $\tilde u$.
	Let $\tilde u' \ne \tilde u$ be another lift.
	Then $\tilde u$ lies in a halfspace $Y'$ of $\tilde u'$, so both ends of $\tilde\gamma$ lie in $C \cap Y'$, and they are not ends of $\tilde u'$ by Lemma~\ref{lem:end-stabilizer}.
	So both lie in $C(Y')$.
\end{proof}

The group $\FF$ permutes the lifts of $u$, and for $g \in \FF$ we speak of the orbits of $\langle g \rangle$ on the set of lifts of $u$.
By Proposition~\ref{prop:orbit-count}, the lifts penetrated by a lift $\tilde\gamma$ of a curve with stabilizer $\langle g \rangle$ form finitely many orbits of $\langle g \rangle$, each infinite, if $\tilde\gamma$ is a core of no lift of $u$.
The following lemma shows that two such families, for non-commuting elements, have a finite intersection.

\begin{lemma} \label{lem:orbits-bound}
	Let $u$ be an essential sphere or a primitive torus, let $a, b \in \FF$ be elements that do not commute, and let $\mathcal{E}$ be a set of lifts of $u$.
	Suppose that the lifts in $\mathcal{E}$ lie in at most $m$ orbits of $\langle a \rangle$ and in at most $m'$ orbits of $\langle b \rangle$ on the set of lifts of $u$.
	Then $\lvert \mathcal{E} \rvert \le 2 m m'$.
\end{lemma}

\begin{proof}
	Suppose that $\lvert \mathcal{E} \rvert \ge 2 m m' + 1$.
	Then three distinct members $E_0$, $E_1$, $E_2$ of $\mathcal{E}$ lie in one orbit of $\langle a \rangle$ and in one orbit of $\langle b \rangle$.
	Write $E_s = a^{\alpha_s} E_0 = b^{\beta_s} E_0$ for $s = 1, 2$.
	The integers $\alpha_1, \alpha_2$ are nonzero and distinct, since $E_0, E_1, E_2$ are distinct, and so are $\beta_1, \beta_2$.
	The elements $v_s := b^{-\beta_s} a^{\alpha_s}$ lie in the stabilizer of $E_0$, which is trivial or infinite cyclic, so they commute; and they are nontrivial, since $\alpha_s, \beta_s \ne 0$ and $a$, $b$ freely generate a free group of rank $2$ by Theorem~\ref{thm:free}\eqref{item:free:basis}.

	In this free group, $v_1$ and $v_2$ are reduced words with two syllables, so neither is a proper power, and commuting elements of a free group are powers of a common element by Theorem~\ref{thm:free}\eqref{item:free:commute}.
	Hence $v_1 = v_2^{\pm 1}$.
	But $v_1 \ne v_2$ since $\alpha_1 \ne \alpha_2$, and $v_1 \ne v_2^{-1} = a^{-\alpha_2} b^{\beta_2}$ since the reduced words begin with different letters.
	This is a contradiction.
\end{proof}

\section{Classes of tori and crossing numbers} \label{sec:crossing}

The goal of this section is the crossing number $N(S,T)$ of two classes of essential spheres or primitive tori: the number of orbits of crossing pairs of their lifts.
It is the quantity through which intersections of surfaces enter the rest of the paper, and Lemma~\ref{lem:count-crossing} gives the form in which it is used.
Before defining the crossing number, we fix what a class of primitive tori is: two primitive tori are in the same class if their dual trees are equivariantly isomorphic.
By Lemma~\ref{lem:tree-sides}, the class determines the ends and sides of the lifts, and by Lemma~\ref{lem:sides-tree}, it is determined by these ends and sides.
Homotopic tori are in the same class, and in the same oriented class if the orientations are carried along the homotopy, by Lemma~\ref{lem:oriented-homotopic}.

\subsection{Classes of primitive tori} \label{subsec:classes}

% This subsection defines classes of primitive tori and shows that the ends, sides and crossings of their lifts depend only on the class.

By Theorem~\ref{thm:laudenbach}\eqref{item:laudenbach:spheres}, the homotopy classes of essential spheres are their isotopy classes, and we call them the \emph{classes} of essential spheres.
By Lemma~\ref{lem:invariant}, being an essential torus depends only on the homotopy class, and by Section~\ref{subsec:tori} so does the core; hence being a primitive torus depends only on the homotopy class.

Let $t$ and $t'$ be primitive tori, with the dual trees $X_t$ and $X_{t'}$ defined in Section~\ref{subsec:curve-lifts}.
A \emph{correspondence} between $t$ and $t'$ is an $\FF$-equivariant bijection $\tilde u \mapsto \tilde u'$ from the lifts of $t$ to the lifts of $t'$ such that $\tilde u$ and $\tilde u'$ have the same stabilizer and the same sides, for every lift $\tilde u$.
Since the ends of a lift are determined by its stabilizer, corresponding lifts also have the same ends.
By Lemma~\ref{lem:crossing-criterion}, a lift of $t$ crosses a lift $\tilde v$ of an essential sphere or torus with no common end if and only if the corresponding lift of $t'$ does.

\begin{lemma} \label{lem:no-leaf}
	$X_t$ has no vertex of valence one.
\end{lemma}

\begin{proof}
	Suppose that the vertex $v$ of $X_t$ has a single edge $e$, and let $\langle y \rangle$ be the stabilizer of $e$.
	The stabilizer of $v$ fixes its only edge, and the stabilizer of $e$ fixes its endpoints, so the stabilizer of $v$ is $\langle y \rangle$.
	Hence $\langle y \rangle$ acts cocompactly on the closure of the region $v$, which covers, with deck group $\langle y \rangle$, a component of the compact manifold obtained by cutting $\MM$ open along $t$.
	Since $v$ has no other edge, this closure is the closure $Y \cup \tilde u_e$ of the halfspace $Y$ of $\tilde u_e$ containing the region $v$.
	So the Cantor part of $Y$ is $\{ y^{\pm\infty} \}$ by Lemma~\ref{lem:closure}, and the side $C(Y)$ of $\tilde u_e$ is empty, contrary to Lemma~\ref{lem:sides-nonempty}.
\end{proof}

\begin{lemma} \label{lem:no-cyclic-vertex}
	No vertex stabilizer of $X_t$ is cyclic.
\end{lemma}

\begin{proof}
	Suppose that the vertex $v$ has cyclic stabilizer $\langle z \rangle$; it contains the infinite cyclic stabilizer of an edge at $v$, so $z \ne 1$.
	The stabilizer of an edge $e$ at $v$ is generated by a conjugate of $y$, which is not a proper power by Lemma~\ref{lem:primitive-basis} and Theorem~\ref{thm:free}\eqref{item:free:primitive}; being contained in $\langle z \rangle$, it is $\langle z \rangle$.
	So $\langle z \rangle$ fixes every edge at $v$, and the valence of $v$ is the number of edge ends at the image of $v$ in the quotient graph $X_t / \FF$, which is one or two.
	Valence one contradicts Lemma~\ref{lem:no-leaf}.
	If the valence is two, the quotient graph is a loop, every vertex is a translate of $v$, and the stabilizers of adjacent vertices are conjugates of $\langle z \rangle$ containing the stabilizer of the edge joining the two vertices, which is a conjugate of $\langle z \rangle$ as well; since $z$ is not a proper power, these three groups coincide.
	So all vertex stabilizers equal $\langle z \rangle$, which is then normal in $\FF$; but a free group of rank at least two has no nontrivial normal cyclic subgroup, since the centralizer of $z$, which is cyclic by Theorem~\ref{thm:free}\eqref{item:free:commute}, would have index at most two.
\end{proof}

\begin{lemma} \label{lem:sides-tree-halfspace}
	Let $\tilde u_1$ and $\tilde u_2$ be distinct lifts of $t$, and let $Z$ be a halfspace of $\tilde u_2$.
	Then $\tilde u_1 \subset Z$ if and only if some side of $\tilde u_1$ is contained in the side $C(Z)$ of $\tilde u_2$.
\end{lemma}

\begin{proof}
	Since $\tilde u_1$ is connected and disjoint from $\tilde u_2$, it lies in a halfspace $Y_2$ of $\tilde u_2$; we write $Y_2^\ast$ for the other halfspace of $\tilde u_2$, and $Y_1$ for the halfspace of $\tilde u_1$ not containing $\tilde u_2$.
	By Lemma~\ref{lem:nested-sides-inclusion}, the side $C(Y_1)$ of $\tilde u_1$ away from $\tilde u_2$ is contained in $C(Y_2)$, and by Lemma~\ref{lem:nested-sides} the other side of $\tilde u_1$ contains the ends of $\tilde u_2$.
	The ends of $\tilde u_2$ lie in no side of $\tilde u_2$, so this other side is contained in no side of $\tilde u_2$.
	Since $C(Y_1)$ is nonempty by Lemma~\ref{lem:sides-nonempty} and the two sides of $\tilde u_2$ are disjoint, $C(Y_1)$ is not contained in $C(Y_2^\ast)$.
	Hence $\tilde u_1 \subset Z$ if and only if some side of $\tilde u_1$ is contained in $C(Z)$, for $Z = Y_2$ and for $Z = Y_2^\ast$.
\end{proof}

\begin{lemma} \label{lem:sides-tree}
	A correspondence $\tilde u \mapsto \tilde u'$ between $t$ and $t'$ induces an $\FF$-equivariant isomorphism $X_t \to X_{t'}$ carrying the edge $\tilde u$ to the edge $\tilde u'$ for every lift $\tilde u$.
\end{lemma}

\begin{proof}
	Let $\tilde u_1$ and $\tilde u_2$ be distinct lifts of $t$.
	By Lemma~\ref{lem:sides-tree-halfspace}, whether a third lift separates $\tilde u_1$ from $\tilde u_2$, that is, whether $\tilde u_1$ and $\tilde u_2$ lie in different halfspaces of the third lift, is determined by the sides of the three lifts.
	So a correspondence preserves this relation.

	In the tree $X_t$, two distinct edges share a vertex if and only if no edge separates them.
	The sets of edges at a vertex are exactly the maximal sets of pairwise adjacent edges, since three pairwise adjacent edges of a tree share a vertex.
	Every vertex of $X_t$ and of $X_{t'}$ has at least two edges by Lemma~\ref{lem:no-leaf}, so the set of edges at a vertex is not contained in the set of edges at another vertex.
	So the correspondence carries the edges at each vertex of $X_t$ onto the edges at a vertex of $X_{t'}$, which defines the isomorphism; it is $\FF$-equivariant since the correspondence is.
\end{proof}

\begin{lemma} \label{lem:tree-sides-criterion}
	Let $\tilde u$ be a lift of $t$ at the edge $e$, with a halfspace $Y$, and let $\kappa$ be the component of $X_t \setminus e$ with $Y_\kappa = Y$.
	Let $\xi \in C$ be a point that is not an end of $\tilde u$, and let $w_n$ be the prefixes of its word.
	Let $x_0$ be a point of $\widetilde\MM$ off the lifts of $t$, and let $v_0$ be the region containing $x_0$.
	Then $\xi \in C(Y)$ if and only if $w_n v_0 \in \kappa$ for all large $n$.
\end{lemma}

\begin{proof}
	For $w \in \FF$, the point $w x_0$ lies in the region $w v_0$.
	Fix a point $x_1$ of the piece $P'_1$ of Section~\ref{sec:univ}.
	Since $w_n x_1 \in P'_{w_n}$, we have $w_n x_1 \to \xi$ by Lemma~\ref{lem:prefix}, and since $d(w_n x_0, w_n x_1) = d(x_0, x_1)$ is bounded, $w_n x_0 \to \xi$ by Lemma~\ref{lem:bounded-limit}.
	Let $W$ be the component of $\mathbb{S}^3 \setminus \hat u$ containing $Y$, so that $C(Y) = C \cap W$ and $W \setminus C = Y$ by Lemma~\ref{lem:sides}.
	If $\xi \in C(Y)$, then $\xi$ lies in the open set $W$, so $w_n x_0 \in W \setminus C = Y$ for all large $n$.
	Then the region $w_n v_0$, which is connected and misses $\tilde u$, lies in $Y$, that is, $w_n v_0 \in \kappa$.
	If $\xi \notin C(Y)$, then $\xi$ lies in the other side of $\tilde u$, and the same argument shows that $w_n v_0$ lies in the other component of $X_t \setminus e$ for all large $n$.
\end{proof}

\begin{lemma} \label{lem:tree-sides}
	Let $\psi \colon X_t \to X_{t'}$ be an $\FF$-equivariant isomorphism.
	Then the lift of $t'$ at $\psi(\tilde u)$ has the same stabilizer and the same sides as $\tilde u$, for every lift $\tilde u$ of $t$.
	So $\tilde u \mapsto \psi(\tilde u)$ is a correspondence.
\end{lemma}

\begin{proof}
	The stabilizers agree since $\psi$ is equivariant, and hence so do the ends.

	For the sides, let $\tilde u$ be a lift of $t$ at the edge $e$, with a halfspace $Y$, and let $\kappa$ be the component of $X_t \setminus e$ with $Y_\kappa = Y$.
	Let $\xi \in C$ be a point that is not an end of $\tilde u$, and let $w_n$ be the prefixes of its word.
	Choose a point $x_0$ off the lifts of $t$, and let $v_0$ be the region containing $x_0$.
	By Lemma~\ref{lem:tree-sides-criterion}, $\xi \in C(Y)$ if and only if $w_n v_0 \in \kappa$ for all large $n$.

	Apply Lemma~\ref{lem:tree-sides-criterion} to the torus $t'$, with its own lifts and regions.
	Take the lift $\tilde u'$ at $\psi(e)$, its halfspace $Y' := Y'_{\psi(\kappa)}$ for the component $\psi(\kappa)$ of $X_{t'} \setminus \psi(e)$, and a base point $x_0'$ in the region $\psi(v_0)$ of $t'$.
	This is possible since $\xi$ is not an end of $\tilde u'$, the ends of $\tilde u$ and $\tilde u'$ being the same.
	For $w \in \FF$, the point $w x_0'$ lies in the region $w \psi(v_0) = \psi(w v_0)$.
	So $\xi \in C(Y')$ if and only if $\psi(w_n v_0) \in \psi(\kappa)$ for all large $n$, that is, if and only if $w_n v_0 \in \kappa$ for all large $n$, if and only if $\xi \in C(Y)$.
	Since neither side contains an end, $C(Y') = C(Y)$, and likewise for the other halfspaces.
\end{proof}

Two primitive tori $t$ and $t'$ are in the same \emph{class} if their dual trees are $\FF$-equivariantly isomorphic, equivalently, by Lemmas~\ref{lem:sides-tree} and~\ref{lem:tree-sides}, if there is a correspondence between $t$ and $t'$.
By Lemmas~\ref{lem:representatives} and~\ref{lem:sides-tree}, homotopic primitive tori are in the same class.
We write $T$ for the class of $t$ and $t \in T$.

If the core of $t$ is oriented, let $o_t$ be the following orientation of $X_t$, in the sense of Section~\ref{subsec:dehn_prelim}.
It assigns to each edge, oriented toward the positive side of its lift, the generator of its stabilizer representing the oriented core.
We write $\vec X_t := (X_t, o_t)$ for the dual tree with this orientation, the \emph{oriented dual tree} of $t$.
The \emph{oriented class} $\vec T$ of $t$ is its class together with $o_t$, and two tori with oriented cores are in the same oriented class if some equivariant isomorphism of their dual trees carries $o_t$ to $o_{t'}$.

By Section~\ref{subsec:dehn_prelim}, the orientations of $X_t$ are $o_t$ and $o_t^{-1}$, and reversing the orientation of the core of $t$ replaces $o_t$ by $o_t^{-1}$.
An equivariant isomorphism $X_t \to X_{t_1}$ pulls $o_{t_1}$ back to an orientation of $X_t$.
So for every oriented class $\vec T$ over $T$ and every $t \in T$, the core of $t$ admits an orientation with $\vec t \in \vec T$.
The next lemma says that homotopic primitive tori are in the same oriented class once the orientations are carried along the homotopy.

\begin{lemma} \label{lem:oriented-homotopic}
	Let $t$ and $t'$ be homotopic primitive tori, and let $G$ be a homotopy from the inclusion of $t$ to a homeomorphism $G_1 \colon t \to t'$.
	Orient $t'$ so that $G_1$ preserves orientation, and orient the core of $t'$ as the image under $G_1$ of an oriented core curve of $t$ lying on $t$.
	Then $\vec t$ and $\vec t'$ are in the same oriented class.
\end{lemma}

If $t'$ carries the orientation for which $G_1$ reverses orientation, apply the lemma after reversing the orientations of $t'$ and of its core, which leaves $o_{t'}$ unchanged, since reversing either one replaces $o_{t'}$ by its inverse.
So in that case $\vec t'$ lies in the oriented class of $\vec t$ when the core of $t'$ is oriented as the image under $G_1$ of the oriented core of $t$ with its orientation reversed.

\begin{proof}
	We first show that the map given by $G$, which is a correspondence by Lemma~\ref{lem:representatives}, preserves positive sides.
	Lift $G$ to a homotopy $\tilde G$ from a lift $\tilde t$ of $t$ to a lift $\tilde t'$ of $t'$; by Section~\ref{subsec:homotopic-lifts}, both lifts have the same stabilizer $\langle y \rangle$, and $\tilde G$ is $\langle y \rangle$-equivariant.
	Let $Y$ and $Y'$ be the halfspaces of $\tilde t$ and $\tilde t'$ on their positive sides.
	Since $G_1$ preserves orientation, Lemma~\ref{lem:same-cantor} gives $C \cap Y = C \cap Y'$, and the ends of the two lifts agree, so $C(Y) = C(Y')$, and likewise for the other halfspaces.
	Extending $\tilde t \mapsto \tilde t'$ equivariantly gives the map given by $G$ of Lemma~\ref{lem:representatives}.
	Since deck transformations carry positive sides to positive sides, it carries the positive side of every lift of $t$ to the positive side of the corresponding lift of $t'$.

	Let $\psi \colon X_t \to X_{t'}$ be the equivariant isomorphism induced by this correspondence by Lemma~\ref{lem:sides-tree}; we show that $\psi$ carries edges oriented toward positive sides to edges oriented toward positive sides.
	Let $\vec e$ be the edge $e$ of $\tilde t$ oriented toward the positive side of $\tilde t$, that is, toward the component $\kappa$ of $X_t \setminus e$ with $Y_\kappa = Y$.
	By Lemma~\ref{lem:sides-tree-halfspace}, an edge of $X_t$ other than $e$ lies in $\kappa$ if and only if one of its sides is contained in $C(Y)$, and $\kappa$ contains edges by Lemma~\ref{lem:no-leaf}.
	So $\psi(\kappa)$ is the component of $X_{t'} \setminus \psi(e)$ whose edges have a side in $C(Y) = C(Y')$, that is, $Y'_{\psi(\kappa)} = Y'$, and $\psi(\vec e)$ is the edge of $\tilde t'$ oriented toward the positive side of $\tilde t'$.

	It remains to show that $o_{t'} \circ \psi = o_t$.
	Let $c$ be the oriented core of $t$, and let $\tilde c$ be its lift in $\tilde t$ from a point $x$, which ends at $o_t(\vec e) \, x$.
	Then $\tilde G_1 \circ \tilde c$ is the lift from $\tilde G_1(x)$ of the oriented core $G_1 \circ c$ of $t'$, and it ends at $\tilde G_1(o_t(\vec e) \, x) = o_t(\vec e) \, \tilde G_1(x)$ by equivariance; so $o_{t'}(\psi(\vec e)) = o_t(\vec e)$.
	Since $\psi$ is equivariant, $o_{t'} \circ \psi$ is an orientation of $X_t$, and an orientation is determined by its value on one oriented edge by Section~\ref{subsec:dehn_prelim}.
	Hence $o_{t'} \circ \psi = o_t$, and $\vec t$ and $\vec t'$ are in the same oriented class.
\end{proof}

By Lemma~\ref{lem:tree-sides}, the ends and sides of the lifts of a torus, and whether a lift crosses a lift of another essential sphere or torus with no common end, depend only on the class.
For a class of curves $\gamma$, put $i(\gamma, T) := i(\gamma, [t])$ for $t \in T$, where $[t]$ is the homotopy class of $t$.
This does not depend on $t$.
If a lift of $\gamma$ is a core of a lift of $t$, then $i(\gamma, [t]) = 0$, since $\gamma$ is represented by a curve in a regular neighborhood of $t$ disjoint from $t$.
Otherwise, by Proposition~\ref{prop:orbit-count}, $i(\gamma, [t])$ is the number of $\langle g \rangle$-orbits of lifts of $t$ penetrated by a lift of $\gamma$ represented by $g$, which depends only on the stabilizers and sides of the lifts.
The core of $t$ is the conjugacy class of the stabilizer of a lift, so it also depends only on $T$.

Two classes of essential spheres or of primitive tori are \emph{core-independent} if, when both are classes of tori, their cores are not freely homotopic; a family of such classes is \emph{pairwise core-independent} if any two of its members are.

\subsection{Crossing numbers} \label{subsec:crossing-numbers}

This subsection defines the crossing number of two classes and relates it to the circles of intersection of representatives.

\begin{lemma} \label{lem:joint-stabilizer}
	Let $u$ and $v$ be essential spheres or primitive tori whose classes are core-independent.
	Then no nontrivial element of $\FF$ stabilizes both a lift of $u$ and a lift of $v$.
\end{lemma}

\begin{proof}
	If $u$ or $v$ is a sphere, this holds since lifts of spheres have trivial stabilizers.
	Otherwise the cores are represented by elements $x$ and $y$ of free bases of $\FF$, by Lemma~\ref{lem:primitive-basis}.
	A common nontrivial element gives conjugates of $x$ and $y$ with a common nonzero power.
	By Theorem~\ref{thm:free}\eqref{item:free:commute}, nontrivial elements of $\FF$ with a common nonzero power are powers of a common element, and by Theorem~\ref{thm:free}\eqref{item:free:primitive}, elements of free bases are not proper powers.
	So $x$ is conjugate to $y^{\pm 1}$, and the cores of $u$ and $v$ are freely homotopic, contrary to core-independence.
\end{proof}

Let $S$ and $T$ be core-independent classes of essential spheres or of primitive tori, with representatives $s$ and $t$.
Let $N(S,T)$ be the number of $\FF$-orbits of crossing pairs $(\tilde s, \tilde t)$ of lifts of $s$ and $t$.
Here $\FF$ acts on the crossing pairs: a deck transformation is an isometry of $\widetilde\MM$ carrying lifts to lifts and halfspaces to halfspaces, so it preserves crossing.

The number $N(S,T)$ does not depend on the choice of $s$ and $t$.
For spheres, the lift of an isotopy given by Theorem~\ref{thm:laudenbach}\eqref{item:laudenbach:spheres} carries lifts to lifts with the same sides, by Lemma~\ref{lem:extend-isotopy}.
For tori, Lemma~\ref{lem:tree-sides} gives a correspondence.
Lifts of $s$ and $t$ have no common end, by Lemmas~\ref{lem:common-end} and~\ref{lem:joint-stabilizer}, so crossing pairs correspond by Lemma~\ref{lem:crossing-criterion}.

We call $N(S,T)$ the \emph{crossing number} of $S$ and $T$.
Since crossing is symmetric in the two lifts, $N(S,T) = N(T,S)$.
Two classes of spheres are core-independent, so $N(S,S)$ is defined for a class $S$ of essential spheres, and it is $0$ by Lemma~\ref{lem:count-disjoint} below, since a sphere is disjoint from a parallel copy.
It is modeled on the intersection number of two splittings of Scott \cite{MR1608688}, in which the crossing of lifts is replaced by the crossing of almost invariant subsets of $\FF$; we do not use that notion.

\begin{lemma} \label{lem:count-crossing}
	Let $\tilde t$ be a lift of $t$ with stabilizer $G$.
	Then the lifts of $s$ crossing $\tilde t$ form exactly $N(S, T)$ orbits under $G$.
\end{lemma}

\begin{proof}
	Since $\FF$ acts transitively on the lifts of $t$, the $\FF$-orbits of crossing pairs of lifts of $s$ and $t$ correspond to the $G$-orbits of lifts of $s$ crossing $\tilde t$.
\end{proof}

Now let $s$ and $t$ be transverse.

\begin{lemma} \label{lem:circles-trivial}
	Every circle $c$ of $s \cap t$ is null-homotopic in $\MM$, so that each lift of $c$ is a circle.
\end{lemma}

\begin{proof}
	The deck transformation given by the class of $c$, based at a point of a lift $\tilde c$, stabilizes the lifts of $s$ and $t$ through that point, so it is trivial by Lemma~\ref{lem:joint-stabilizer}.
	Hence $c$ is null-homotopic, and each lift of $c$ is a circle.
\end{proof}

\begin{lemma} \label{lem:circles-trivial-lifts}
	Each lift $\tilde c$ of a circle $c$ of $s \cap t$ lies on exactly one lift of $s$ and on exactly one lift of $t$.
\end{lemma}

\begin{proof}
	The connected set $\tilde c$ lies in $\rho^{-1}(s)$, whose components are the lifts of $s$, so it lies on exactly one lift of $s$, the lift through any given point of $\tilde c$.
	The same holds for $t$.
\end{proof}

Assign to each circle $c$ of $s \cap t$ the $\FF$-orbit of the pair of lifts through a lift of $c$.
By Lemmas~\ref{lem:circles-trivial} and~\ref{lem:circles-trivial-lifts}, this is well-defined, since the lifts of $c$ are the $\FF$-translates of any one lift.

\begin{lemma} \label{lem:count}
	$\lvert \pi_0(s \cap t) \rvert \ge \#\{\FF\text{-orbits of intersecting pairs of lifts}\} \ge N(S,T)$.
\end{lemma}

\begin{proof}
	The assignment is onto the orbits of intersecting pairs, which gives the first inequality.
	Crossing pairs intersect, since disjoint lifts do not cross; this gives the second.
\end{proof}

\begin{lemma} \label{lem:count-disjoint}
	If $S$ and $T$ have disjoint representatives, then $N(S,T) = 0$.
\end{lemma}

\begin{proof}
	Let $s' \in S$ and $t' \in T$ be disjoint; they are transverse, so Lemma~\ref{lem:count} applied to $s'$ and $t'$ gives $N(S,T) \le \lvert \pi_0(s' \cap t') \rvert = 0$.
\end{proof}

\section{Dehn twists} \label{sec:dehn-twist}

The goal of this section is Theorems~\ref{thm:commuting} and~\ref{thm:commuting-crossing}: two Dehn twists along primitive tori with core-independent classes commute if the classes have disjoint representatives, and only if the crossing number of the two classes is zero.
On the way we show in Proposition~\ref{prop:coincide} that a Dehn twist in $\Out(\FF)$ determines the oriented class of its torus.
Proposition~\ref{prop:twist-periodic} shows that no nonzero power of a Dehn twist fixes a class with positive crossing number with the torus.

\subsection{Dehn twist homeomorphisms} \label{subsec:twist-homeo}

This subsection defines the Dehn twist homeomorphism along a primitive torus and identifies its image in $\Out(\FF)$.

We use the one-edge $\ZZ$-splittings, their orientations $\vec T$ and the Dehn twists $\delta_{\vec T} \in \Out(\FF)$ of Section~\ref{subsec:dehn_prelim}, and the model $\widetilde\MM = \mathbb{S}^3 \setminus C$ with its ends and Cantor parts of Section~\ref{subsec:ends}.

Let $t$ be a primitive torus; its dual tree $X_t$ is a one-edge $\ZZ$-splitting by Section~\ref{subsec:curve-lifts}.

Recall from Section~\ref{subsec:doubled_handlebody} that $\MM$ and every torus in $\MM$ are oriented, so that $t$ has a positive side.
Orient the core of $t$, and write $\vec t$ for $t$ with this orientation.
Let $N(t) \cong [0,1] \times \mathbb{S}^1 \times \mathbb{S}^1$ be a closed regular neighborhood of $t$, parametrized by an orientation-preserving homeomorphism so that $t = \{1/2\} \times \mathbb{S}^1 \times \mathbb{S}^1$, the circle $\{1/2\} \times \{0\} \times \mathbb{S}^1$ is an oriented core of $\vec t$ and $\{1\} \times \mathbb{S}^1 \times \mathbb{S}^1$ lies on the positive side of $t$.
Let $\delta$ be the homeomorphism of $[0,1] \times \mathbb{S}^1 \times \mathbb{S}^1$, where $\mathbb{S}^1 = \mathbb{R}/\mathbb{Z}$, given by $\delta(r, \theta, \phi) = (r, \theta, \phi + r)$.
The \emph{Dehn twist homeomorphism} $\delta_{\vec t}$ of $\MM$ is $\delta$ on $N(t)$ and the identity elsewhere.

We compute the lifts of $\delta_{\vec t}$ to the lifts of $N(t)$.
The circle $\{1/2\} \times \mathbb{S}^1 \times \{0\}$ lies on a torus parallel to $t$, so by Section~\ref{subsec:tori} it represents, as a loop based at $(1/2, 0, 0)$, the $b$-th power of the class of the oriented core $\{1/2\} \times \{0\} \times \mathbb{S}^1$ in $\pi_1(\MM, (1/2, 0, 0))$ for some $b \in \ZZ$.
Composing the parametrization with the shear $(r, \theta, \phi) \mapsto (r, \theta, \phi - b \theta)$, which preserves orientation, fixes $\{1/2\} \times \{0\} \times \mathbb{S}^1$ pointwise and commutes with $\delta$, we may assume that $b = 0$ without changing $\delta_{\vec t}$.

Let $\tilde t$ be a lift of $t$ and $w$ the generator of its stabilizer representing the oriented core.
Then the lift of $N(t)$ containing $\tilde t$ is the covering $[0,1] \times \mathbb{S}^1 \times \mathbb{R} \to N(t)$ given by the quotient map $\mathbb{R} \to \mathbb{R} / \ZZ = \mathbb{S}^1$ in the last coordinate.
On this covering, $w$ acts by $(r, \theta, \tilde\phi) \mapsto (r, \theta, \tilde\phi + 1)$.
The lift of $\delta$ that is the identity on $\{0\} \times \mathbb{S}^1 \times \mathbb{R}$ is $(r, \theta, \tilde\phi) \mapsto (r, \theta, \tilde\phi + r)$, and on $\{1\} \times \mathbb{S}^1 \times \mathbb{R}$, which faces the positive side of $\tilde t$, it is $w$.

\begin{lemma}[{\cite[Lemma 2.6]{MR3677931}}] \label{lem:gul}
	Let $o_t$ be the orientation of $X_t$, which assigns to the edge of a lift $\tilde t$ oriented toward the positive side of $\tilde t$ the generator of the stabilizer of $\tilde t$ representing the oriented core.
	Then the image of $\delta_{\vec t}$ under the homomorphism $\MCG(\MM) \to \Out(\FF)$ is $\delta_{\vec X_t}$.
\end{lemma}

\begin{proof}
	We include the proof, which fixes the sign convention.
	Each component of $\widetilde\MM$ minus the lifts of $N(t)$ lies in a unique region of $t$, and we identify its closure with that region, a vertex of $X_t$.
	Let $\Omega_0$ be a region and $\tilde\delta$ the lift of $\delta_{\vec t}$ that is the identity on $\Omega_0$.
	Since $\delta_{\vec t}$ is the identity off $N(t)$, the restriction of $\tilde\delta$ to a region $\Omega$ is a deck transformation $g_\Omega$, and $\tilde\delta$ permutes the lifts of $N(t)$, so it induces an automorphism of $X_t$ carrying $\Omega$ to $g_\Omega \Omega$.

	Let $\Omega$ and $\Omega'$ be the regions adjacent to a lift $\tilde t'$ of $t$, and $\vec e$ the edge $\tilde t'$ oriented from $\Omega$ to $\Omega'$.
	Consider the lift of $\delta_{\vec t}$ that is the identity on the region on the negative side of $\tilde t'$.
	By the computation of the lift of $\delta$ on the lift $[0,1] \times \mathbb{S}^1 \times \mathbb{R}$ of $N(t)$ containing $\tilde t'$, made before this lemma, this lift of $\delta_{\vec t}$ is, on the region on the positive side, the generator of the stabilizer of $\tilde t'$ representing the oriented core.
	If $\Omega'$ is on the positive side, this generator is $o_t(\vec e)$, so this lift is the identity on $\Omega$ and $o_t(\vec e)$ on $\Omega'$.
	If $\Omega$ is on the positive side, this generator is $o_t(\ast \vec e)$, and composing that lift with $o_t(\ast \vec e)^{-1} = o_t(\vec e)$ gives the lift that is the identity on $\Omega$ and $o_t(\vec e)$ on $\Omega'$.
	In both cases $\tilde\delta$ is $g_\Omega$ composed with the lift that is the identity on $\Omega$, so $g_{\Omega'} = g_\Omega \, o_t(\vec e)$.

	By induction along the geodesic $\vec e_1, \dots, \vec e_k$ from $\Omega_0$ to $\Omega$, we get $g_\Omega = o_t(\vec e_1) \dots o_t(\vec e_k)$, so $\tilde\delta$ induces $\delta_{\vec X_t, \Omega_0}$ on $X_t$.

	Let $\varphi \in \Aut(\FF)$ be defined by $\tilde\delta(g x) = \varphi(g) \tilde\delta(x)$, so that the image of $\delta_{\vec t}$ in $\Out(\FF)$ is the outer class of $\varphi$.
	On $X_t$, this reads $\varphi(g) = \delta_{\vec X_t, \Omega_0} \, g \, \delta_{\vec X_t, \Omega_0}^{-1}$, and since $\FF$ acts faithfully on $X_t$ by Section~\ref{subsec:dehn_prelim}, $\varphi$ is the automorphism $\varphi_{\Omega_0}$ of that section, whose outer class is $\delta_{\vec X_t}$.
\end{proof}

By Lemma~\ref{lem:gul}, the image of $\delta_{\vec t}$ in $\Out(\FF)$ is $\delta_{\vec X_t}$, which depends only on the oriented class $\vec T$ of $\vec t$ of Section~\ref{subsec:classes}.
Indeed, an equivariant isomorphism $\psi \colon \vec X \to \vec X'$ of oriented trees satisfies $\psi \, \delta_{\vec X, v_0} \, \psi^{-1} = \delta_{\vec X', \psi(v_0)}$, so the two twists induce the same automorphism of $\FF$.
We write $\delta_{\vec T} \in \Out(\FF)$ for this image; this agrees with the notation of Section~\ref{subsec:dehn_prelim} for the oriented splitting $\vec X_t$.
Reversing the orientation of the core, or that of the torus, replaces $o_t$ by its inverse and $\delta_{\vec T}$ by its inverse.

Let $h$ be an orientation-preserving homeomorphism of $\MM$ and $\tilde h$ a lift, with $\tilde h(g x) = \varphi(g) \tilde h(x)$ for some $\varphi \in \Aut(\FF)$ in the outer class of $h$.
For a primitive torus $t$, the lift $\tilde h$ carries the regions and lifts of $t$ to the regions and lifts of $h(t)$.
So it induces an isomorphism $X_t \to X_{h(t)}$ that is equivariant with respect to $\varphi$ and carries $o_t$ to $o_{h(t)}$ for the transported orientations.
Hence the oriented class of $h(\vec t)$ is that of $\vec t$ with the action precomposed with $\varphi^{-1}$.
It depends only on $\vec T$ and on the outer class of $\varphi$, since precomposing with conjugation by $a \in \FF$ gives an isomorphic oriented tree via $x \mapsto a x$.

For an essential sphere $s$, the class of $h(s)$ is the homotopy class of $h(s)$, and the kernel of $\MCG(\MM) \to \Out(\FF)$ is generated by sphere twists, which preserve the class of every essential sphere, by Section~\ref{subsec:doubled_handlebody}.
So an element of the kernel fixes every class of essential spheres, and, since its $\varphi$ is inner, every oriented class of primitive tori.
Hence $\MCG(\MM)$ acts on classes of essential spheres and on oriented classes of primitive tori through $\Out(\FF)$.
For a class $S$ of essential spheres or of primitive tori, we write $\delta_{\vec T}(S)$ for the class of $\delta_{\vec t}(s)$ for any $s \in S$, a sphere or a torus.
For an oriented class $\vec S$ of primitive tori, we write $\delta_{\vec T}(\vec S)$ for the oriented class of $\delta_{\vec t}(\vec s)$, with the transported orientations, for any $\vec s \in \vec S$.

\begin{proposition} \label{prop:coincide}
	Let $\vec T$ and $\vec T'$ be oriented classes of primitive tori.
	Then $\delta_{\vec T} = \delta_{\vec T'}$ in $\Out(\FF)$ if and only if $\vec T = \vec T'$.
\end{proposition}

The proof is given at the end of Section~\ref{subsec:twist-periodic}.

\subsection{The lift of a Dehn twist} \label{subsec:twist-lift}

This subsection describes the lifts of powers of a Dehn twist to $\widetilde\MM$, and how they act on the Cantor set $C$ and on the sides of lifts of tori.

Let $t$ be a primitive torus with oriented core, and $\delta := \delta_{\vec t}$, supported in the regular neighborhood $N(t)$.
In this section and in Sections~\ref{sec:growth} and~\ref{sec:fellow}, a \emph{region} of $t$ is the closure of a component of $\widetilde\MM$ minus the lifts of $N(t)$.
It is a region of $t$ in the sense of Section~\ref{subsec:curve-lifts} with the open regular neighborhoods of the adjacent lifts removed, and we identify the two.

Fix a lift $\tilde t$ and the generator $w$ of its stabilizer representing the oriented core.
Let $W_0$ and $W_1$ be the regions adjacent to the lifts of $\{0\} \times \mathbb{S}^1 \times \mathbb{S}^1$ and of $\{1\} \times \mathbb{S}^1 \times \mathbb{S}^1$ in the lift of $N(t)$ containing $\tilde t$, and $Y^\ast \supset W_0$ and $Y \supset W_1$ the halfspaces of $\tilde t$.
The stabilizer $\Stab(W_j)$ of $W_j$ acts cocompactly on $W_j$, so $C \cap W_j$ is the limit set $\partial \Stab(W_j)$ of $\Stab(W_j)$ by Lemma~\ref{lem:bounded-limit}, and $\langle w \rangle \le \Stab(W_0) \cap \Stab(W_1)$.

For a lift $\tilde t'$ of $t$, let $N(\tilde t')$ be the component of $\rho^{-1}(N(t))$ containing $\tilde t'$.
Since $\tilde t'$ separates $\widetilde\MM$, the complement of $N(\tilde t')$ has two components.
For a lift $\tilde t' \ne \tilde t$ adjacent to $W_0$ or to $W_1$, let $U_{\tilde t'}$ be the closure of the component not containing $\tilde t$; it is the part of $\widetilde\MM$ beyond $\tilde t'$ as seen from $\tilde t$.

\begin{lemma} \label{lem:halfspace-decomposition}
	The halfspace $Y$ decomposes as
	\[
		Y = \bigl( N(\tilde t) \cap Y \bigr) \,\cup\, W_1 \,\cup\, \bigcup_{\tilde t'} \bigl( N(\tilde t') \cup U_{\tilde t'} \bigr) ,
	\]
	the union over the lifts $\tilde t' \ne \tilde t$ of $t$ adjacent to $W_1$.
\end{lemma}

\begin{proof}
	The boundary of $W_1$ consists of one boundary component of $N(\tilde t'')$ for each lift $\tilde t''$ of $t$ adjacent to $W_1$, including $\tilde t$; each is a lift of a torus parallel to $t$.
	By Lemma~\ref{lem:boundary} applied to $V = W_1$, the space $\widetilde\MM$ is the union of $W_1$ and, for each such $\tilde t''$, the closed set beyond that boundary component.
	For $\tilde t'' = \tilde t$ this set is $N(\tilde t) \cup \overline{Y^\ast}$, which meets $Y$ in $N(\tilde t) \cap Y$.
	For $\tilde t'' = \tilde t' \ne \tilde t$ it is $N(\tilde t') \cup U_{\tilde t'}$, which lies in $Y$, since it is connected, misses $\tilde t$ and meets $W_1 \subset Y$.
	Intersecting with $Y$ gives the formula.
\end{proof}

\begin{lemma} \label{lem:halfspace-decomposition-cantor}
	$C \cap Y = \partial \Stab(W_1) \cup \bigcup_{\tilde t'} (C \cap U_{\tilde t'})$, the union over the lifts $\tilde t' \ne \tilde t$ of $t$ adjacent to $W_1$.
\end{lemma}

\begin{proof}
	The right side is contained in $C \cap Y$, since $W_1$ and the sets $U_{\tilde t'}$ lie in $Y$ by Lemma~\ref{lem:halfspace-decomposition}.
	Conversely, let $\xi \in C \cap Y$, the limit of a sequence $z_n \in Y$.
	If infinitely many $z_n$ lie in $W_1$ or in the adjacent lifts of $N(t)$, then $\xi \in \partial \Stab(W_1)$ by Lemma~\ref{lem:bounded-limit}, and if infinitely many lie in one $U_{\tilde t'}$, then $\xi \in C \cap U_{\tilde t'}$.

	Otherwise, after passing to a subsequence, $z_n \in U_{\tilde t'_n}$ with the $\tilde t'_n$ pairwise distinct, so the $\tilde t'_n$ leave every compact set by local finiteness of the lifts; fix $x_1 \in W_1$.
	The set $U_{\tilde t'_n}$ lies in the halfspace of $\tilde t'_n$ not containing $x_1$.
	By Lemma~\ref{lem:far-halfspace}, the distance in $\mathbb{S}^3$ from $z_n$ to a point $p_n \in \tilde t'_n$ tends to $0$, so $p_n \to \xi$.
	Since $\tilde t'_n$ is adjacent to $W_1$, there is $q_n \in W_1$ at bounded distance from $p_n$, and $q_n \to \xi$ by Lemma~\ref{lem:bounded-limit}.
	Thus $\xi \in C \cap W_1 = \partial \Stab(W_1)$.
\end{proof}

Let $\tilde\delta$ be the lift of $\delta$ that is the identity on $W_0$.
On a region $\Omega$ it is a deck transformation $g_\Omega$.
By the proof of Lemma~\ref{lem:gul}, across the lift of $N(t)$ at any lift $\tilde t'$ of $t$, possibly $\tilde t$ itself, it changes by a generator of the stabilizer of $\tilde t'$.
Namely $g_{\Omega'} = g_\Omega \, o_t(\vec e)$ for the regions $\Omega$, $\Omega'$ adjacent to $\tilde t'$ and the edge $\vec e$ of $\tilde t'$ oriented from $\Omega$ to $\Omega'$.
In particular $g_{W_1} = o_t(\vec e) = w$ for the edge $\vec e$ of $\tilde t$ oriented from $W_0$ to $W_1$, that is, toward the positive side of $\tilde t$.
Since $\tilde\delta$ permutes the lifts of $N(t)$, it carries $U_{\tilde t'}$ onto $g_\Omega U_{\tilde t'}$, where $\Omega$ is the region adjacent to $\tilde t'$ in $U_{\tilde t'}$.

Let $\varphi \in \Aut(\FF)$ be defined by $\tilde\delta(gx) = \varphi(g) \tilde\delta(x)$.
Then $\varphi$ is the identity on $\Stab(W_0)$ and conjugation by $w$ on $\Stab(W_1)$; in particular $\varphi(w) = w$.
By Lemmas~\ref{lem:qi-extends} and~\ref{lem:qi-extends-cantor}, $\tilde\delta$ extends to a homeomorphism of $\mathbb{S}^3$ with $\tilde\delta(C \cap A) = C \cap \tilde\delta(A)$.

\begin{lemma} \label{lem:twist-fixes}
	$\tilde\delta$ preserves $\tilde t$, $Y$ and $Y^\ast$.
\end{lemma}

\begin{proof}
	$\tilde\delta(W_0) = W_0$ and $\tilde\delta(W_1) = w W_1 = W_1$, so $\tilde\delta$ preserves the lift of $N(t)$ between $W_0$ and $W_1$ and its two sides.
	Since $\delta$ preserves each torus $\{r\} \times \mathbb{S}^1 \times \mathbb{S}^1$ of $N(t)$, in particular $t$, the lift $\tilde\delta$ preserves $\tilde t$, hence its halfspaces $Y \supset W_1$ and $Y^\ast \supset W_0$.
\end{proof}

\begin{lemma} \label{lem:twist-fixes-negative}
	$\tilde\delta$ fixes $\partial \Stab(W_0)$ pointwise.
\end{lemma}

\begin{proof}
	Let $\xi = \lim a_n x_0 \in \partial \Stab(W_0)$ with $a_n \in \Stab(W_0)$ and $x_0 \in W_0$.
	Then $\tilde\delta(\xi) = \lim \varphi(a_n) \tilde\delta(x_0) = \lim a_n x_0 = \xi$.
\end{proof}

\begin{lemma} \label{lem:twist-fixes-positive}
	$\tilde\delta$ equals $w$ on $\partial \Stab(W_1)$.
\end{lemma}

\begin{proof}
	Let $\xi = \lim a_n x_1 \in \partial \Stab(W_1)$ with $a_n \in \Stab(W_1)$ and $x_1 \in W_1$.
	Then $\tilde\delta(\xi) = \lim w a_n w^{-1} \tilde\delta(x_1) = \lim w a_n x_1 = w \xi$, since $\tilde\delta(x_1) = w x_1$.
\end{proof}

The construction applies to any region $W$ in place of $W_0$ and to any power of $\delta$.
For $n \in \ZZ$, let $\tilde\delta_{W,n}$ be the lift of $\delta^{\,n}$ that is the identity on $W$; it is the $n$-th power of the lift of $\delta$ that is the identity on $W$.
By Lemmas~\ref{lem:qi-extends} and~\ref{lem:qi-extends-cantor}, $\tilde\delta_{W,n}$ extends to a homeomorphism of $\mathbb{S}^3$ with $\tilde\delta_{W,n}(C \cap A) = C \cap \tilde\delta_{W,n}(A)$, and by Lemma~\ref{lem:twist-fixes-negative} it fixes $C \cap W$ pointwise.
The lemmas from here to Lemma~\ref{lem:twist-curve} are used in Sections~\ref{sec:growth} and~\ref{sec:fellow}.

For the next three lemmas, let $\tilde t'$ be any lift of $t$ adjacent to $W$, possibly $\tilde t$ itself, and $\Omega$ the region on the other side of $\tilde t'$.

\begin{lemma} \label{lem:twist-adjacent}
	$\tilde\delta_{W,n}$ equals $u^n$ on $\Omega$ for a generator $u$ of the stabilizer of $\tilde t'$.
\end{lemma}

\begin{proof}
	The lift $\tilde\delta_{W,1}$ is a deck transformation $g_\Omega$ on $\Omega$, and $g_\Omega$ is a generator $u$ of the stabilizer of $\tilde t'$, since $g_W = 1$.
	So $\tilde\delta_{W,1}(\Omega) = u \Omega = \Omega$, and $\tilde\delta_{W,n} = \tilde\delta_{W,1}^n$ equals $u^n$ on $\Omega$.
\end{proof}

\begin{lemma} \label{lem:twist-adjacent-sides}
	$\tilde\delta_{W,n}$ preserves $\tilde t'$, its halfspaces and its sides.
\end{lemma}

\begin{proof}
	By Lemma~\ref{lem:twist-adjacent}, $\tilde\delta_{W,n}$ is the identity on $W$ and equals $u^n$ on $\Omega$, where $u$ generates the stabilizer of $\tilde t'$.
	So $\tilde\delta_{W,n}(W) = W$ and $\tilde\delta_{W,n}(\Omega) = u^n \Omega = \Omega$, since the stabilizer of $\tilde t'$ preserves the two halfspaces of $\tilde t'$ by Section~\ref{subsec:curve-lifts}, hence the two regions adjacent to $\tilde t'$.
	Since $\tilde\delta_{W,n}$ permutes the lifts of $N(t)$, it preserves the lift of the regular neighborhood $N(t)$ between $W$ and $\Omega$, which is the only lift of $N(t)$ adjacent to both $W$ and $\Omega$; since $\delta$ preserves the torus $t = \{1/2\} \times \mathbb{S}^1 \times \mathbb{S}^1$ of $N(t)$, it preserves $\tilde t'$.
	The halfspace of $\tilde t'$ containing $W$ is preserved since $W$ is fixed, hence so is the other, and the sides are preserved by Lemma~\ref{lem:qi-extends-cantor}.
\end{proof}

\begin{lemma} \label{lem:twist-adjacent-image}
	Let $u$ be the generator of Lemma~\ref{lem:twist-adjacent}.
	For a lift $\tilde t''$ of $t$ adjacent to $\Omega$, we have $\tilde\delta_{W,n}(\tilde t'') = u^n \tilde t''$.
\end{lemma}

\begin{proof}
	The lift $\tilde\delta_{W,1}$ carries each lift adjacent to $\Omega$ to its image under $g_\Omega = u$, so $\tilde\delta_{W,n} = \tilde\delta_{W,1}^n$ gives $\tilde\delta_{W,n}(\tilde t'') = u^n \tilde t''$.
\end{proof}

Let $\tilde t_1$ and $\tilde t_2$ be disjoint lifts of primitive tori $t_1$ and $t_2$ whose stabilizers have no common nontrivial element.
Then $\tilde t_1$ and $\tilde t_2$ have no common end by Lemma~\ref{lem:common-end}.
The lift $\tilde t_2$ lies in one halfspace $Y_0$ of $\tilde t_1$, so its ends lie in $C \cap Y_0$, which is the side $C(Y_0)$ together with the ends of $\tilde t_1$ by Section~\ref{subsec:ends}; hence they lie in $C(Y_0)$.
So the \emph{side of $\tilde t_1$ not containing the ends of $\tilde t_2$} is well defined, and it is the side of the halfspace of $\tilde t_1$ not containing $\tilde t_2$.

For the next two lemmas and the corollary, let $\tilde t_1$, $\tilde t_2$, $\tilde t_3$ be lifts of primitive tori $t_1$, $t_2$, $t_3$, not necessarily distinct, such that $\tilde t_1$ and $\tilde t_3$ are disjoint from $\tilde t_2$ and their stabilizers have no common nontrivial element with the stabilizer of $\tilde t_2$.
Let $P$ be the side of $\tilde t_1$, and $R$ the side of $\tilde t_3$, not containing the ends of $\tilde t_2$; both are well defined, since each of the pairs $(\tilde t_1, \tilde t_2)$ and $(\tilde t_3, \tilde t_2)$ consists of disjoint lifts whose stabilizers have no common nontrivial element.
Let $f$ be a lift of a homeomorphism of $\MM$ that preserves $\tilde t_1$ with its sides and $\tilde t_2$, and satisfies $f(\tilde t_3) = w^n \tilde t_3$, where $w$ stabilizes $\tilde t_2$ and $n \in \ZZ$.
By Lemmas~\ref{lem:qi-extends} and~\ref{lem:qi-extends-cantor}, $f$ extends to a homeomorphism of $\mathbb{S}^3$ with $f(C \cap A) = C \cap f(A)$ for every $A \subset \widetilde\MM$.

\begin{lemma} \label{lem:twist-sides}
	$f(P) = P$.
\end{lemma}

\begin{proof}
	Since $f$ preserves the sides of $\tilde t_1$, we have $f(P) = P$.
\end{proof}

\begin{lemma} \label{lem:twist-sides-moved}
	$f(R) = w^n R$, and this is the side of $w^n \tilde t_3$ not containing the ends of $\tilde t_2$.
\end{lemma}

\begin{proof}
	Let $Y''$ be the halfspace of $\tilde t_3$ not containing $\tilde t_2$, so that $R = C(Y'')$ is $C \cap Y''$ minus the ends of $\tilde t_3$.
	Since $f(\tilde t_2) = \tilde t_2$ and $f(\tilde t_3) = w^n \tilde t_3$, the homeomorphism $f$ carries $Y''$ onto the halfspace of $w^n \tilde t_3$ not containing $\tilde t_2$, which is $w^n Y''$ since $w^n \tilde t_2 = \tilde t_2$.
	It carries the ends of $\tilde t_3$, which form $C \cap \tilde t_3$, to the ends of $w^n \tilde t_3$, so $f(R) = C(w^n Y'') = w^n R$.
	The deck transformation $w^n$ fixes $\tilde t_2$ and its ends, so it carries $R$ to the side of $w^n \tilde t_3$ not containing the ends of $\tilde t_2$, and this side is $w^n R$.
\end{proof}

\begin{corollary} \label{lem:twist-sides-curve}
	$f$ carries the ends of a lift of a curve with one end in $P$ and one end in $R$ to a point of $P$ and a point of $w^n R$.
\end{corollary}

\begin{proof}
	This follows from $f(P) = P$ and $f(R) = w^n R$, by Lemmas~\ref{lem:twist-sides} and~\ref{lem:twist-sides-moved}.
\end{proof}

\begin{lemma} \label{lem:twist-curve}
	Let $\tilde\gamma$ be a lift of a class $\gamma$ with ends $g^{\pm\infty}$.
	Then $\tilde\delta_{W,n}(\tilde\gamma) := \{ \tilde\delta_{W,n}(g^{+\infty}), \tilde\delta_{W,n}(g^{-\infty}) \}$ is a lift of $\delta_{\vec T}^{\,n}(\gamma)$, where $T$ is the class of $t$.
\end{lemma}

\begin{proof}
	Let $\varphi_{W,n} \in \Aut(\FF)$ satisfy $\tilde\delta_{W,n}(h x) = \varphi_{W,n}(h) \tilde\delta_{W,n}(x)$.
	Then $\tilde\delta_{W,n}(g^m x) = \varphi_{W,n}(g)^m \tilde\delta_{W,n}(x)$, so $\tilde\delta_{W,n}(g^{\pm\infty}) = \varphi_{W,n}(g)^{\pm\infty}$ by Lemma~\ref{lem:closure-limit}.
	The outer class of $\varphi_{W,n}$ is the image of $\delta^{\,n}$ in $\Out(\FF)$, which is $\delta_{\vec T}^{\,n}$ by Lemma~\ref{lem:gul}; so $\varphi_{W,n}(g)$ represents $\delta_{\vec T}^{\,n}(\gamma)$.
\end{proof}

\subsection{Periodic classes} \label{subsec:twist-periodic}

The goal of this subsection is Proposition~\ref{prop:twist-periodic}: no nonzero power of $\delta_{\vec T}$ fixes $S$ if $N(S,T) > 0$.
We then show in Lemma~\ref{lem:dual-sphere-crosses} and Corollary~\ref{lem:dual-sphere-crosses-positive} that a core penetrating a surface yields a crossing pair, deduce that $\delta_{\vec T}$ has infinite order, and prove Proposition~\ref{prop:coincide}.

Let $\vec T$ be an oriented class of primitive tori, with underlying class $T$, and $S$ a class of essential spheres or of primitive tori such that $S$ and $T$ are core-independent.
Until the proof of Proposition~\ref{prop:twist-periodic}, suppose that $N(S,T) > 0$, and let $\delta = \delta_{\vec t}$ for a torus $\vec t$ with oriented core in the oriented class $\vec T$.
Choose $s \in S$ transverse to $t$ and crossing lifts $\tilde s$, $\tilde t$.

We use the notation of Section~\ref{subsec:twist-lift} for the lift $\tilde t$.
The element $w$ is the generator of the stabilizer of $\tilde t$ representing the oriented core.
The regions $W_0$ and $W_1$ are the two regions adjacent to $N(\tilde t)$, on the negative and on the positive side of $\tilde t$, and $Y^\ast \supset W_0$ and $Y \supset W_1$ are the halfspaces of $\tilde t$.
For a lift $\tilde t' \ne \tilde t$ of $t$ adjacent to $W_0$ or to $W_1$, the set $U_{\tilde t'}$ is the part of $\widetilde\MM$ beyond $N(\tilde t')$ as seen from $\tilde t$.
The lift $\tilde\delta$ of $\delta$ is the identity on $W_0$, and it equals the deck transformation $g_\Omega$ on each region $\Omega$; in particular $g_{W_0} = 1$ and $g_{W_1} = w$.

We first record how $\tilde\delta$ moves the points of $C$ on the two sides of $\tilde t$.

\begin{lemma} \label{lem:twist-attracts}
	Let $\xi \in (C \cap Y) \setminus \{ w^{-\infty} \}$.
	Then $\tilde\delta^k(\xi) \to w^{+\infty}$ as $k \to +\infty$.
\end{lemma}

\begin{proof}
	If $\xi \in \partial \Stab(W_1)$, then $\tilde\delta^k(\xi) = w^k \xi \to w^{+\infty}$ by Lemmas~\ref{lem:twist-fixes-positive} and~\ref{lem:north-south}.
	Otherwise $\xi \in C \cap U_{\tilde t'}$ for a lift $\tilde t'$ adjacent to $W_1$, by Lemma~\ref{lem:halfspace-decomposition-cantor}.
	The region $\Omega$ adjacent to $\tilde t'$ in $U_{\tilde t'}$ has $g_\Omega = w u^{\pm 1}$, with $u$ generating the stabilizer of $\tilde t'$, so $\tilde\delta(U_{\tilde t'}) = g_\Omega U_{\tilde t'} = w U_{\tilde t'}$.
	By induction, $\tilde\delta^k(U_{\tilde t'}) = w^k U_{\tilde t'}$, since $\varphi(w) = w$.

	The compact set $C \cap U_{\tilde t'}$ misses $w^{-\infty}$.
	Indeed, $w^{-\infty}$ lies in $\partial \Stab(W_1) = C \cap W_1$, and $W_1$ and $U_{\tilde t'}$ lie in the closures of opposite halfspaces of $\tilde t'$, so by Lemma~\ref{lem:sides-partition} their Cantor parts meet only in the ends of $\tilde t'$.
	But $w^{-\infty}$ is not an end of $\tilde t'$, since $\tilde t \ne \tilde t'$ and distinct lifts have no common end by Lemma~\ref{lem:end-stabilizer}.
	So $w^k (C \cap U_{\tilde t'}) \to w^{+\infty}$ by Lemma~\ref{lem:north-south}.
\end{proof}

\begin{lemma} \label{lem:twist-repels}
	Let $\xi \in (C \cap Y^\ast) \setminus \{ w^{+\infty} \}$.
	Then $\xi$ lies in a closed $\tilde\delta$-invariant subset of $C$ not containing $w^{+\infty}$.
\end{lemma}

\begin{proof}
	If $\xi \in \partial \Stab(W_0)$, take $\{ \xi \}$, which is fixed by Lemma~\ref{lem:twist-fixes-negative}.
	Otherwise $\xi \in C \cap U_{\tilde t'}$ for a lift $\tilde t'$ adjacent to $W_0$, by Lemma~\ref{lem:halfspace-decomposition-cantor} applied to $Y^\ast$.
	Here $g_\Omega = u^{\pm 1}$, so $\tilde\delta(U_{\tilde t'}) = U_{\tilde t'}$, and $C \cap U_{\tilde t'}$ is $\tilde\delta$-invariant.
	It misses $w^{+\infty}$.
	Indeed, $w^{+\infty}$ lies in $\partial \Stab(W_0) = C \cap W_0$, and $W_0$ and $U_{\tilde t'}$ lie in the closures of opposite halfspaces of $\tilde t'$, so by Lemma~\ref{lem:sides-partition} their Cantor parts meet only in the ends of $\tilde t'$.
	But $w^{+\infty}$ is not an end of $\tilde t'$, since $\tilde t \ne \tilde t'$ and distinct lifts have no common end by Lemma~\ref{lem:end-stabilizer}.
\end{proof}

Let $H$ and $H^\ast$ be the Cantor parts of the halfspaces of $\tilde s$.
Then $H \cap H^\ast$ is the Cantor part of $\tilde s$, which misses $w^{\pm\infty}$ by Lemmas~\ref{lem:common-end} and~\ref{lem:joint-stabilizer}.
By Lemma~\ref{lem:corner-ends}, there are
\[
	\xi \in H \cap (C \cap Y) , \quad
	\eta \in H^\ast \cap (C \cap Y) , \quad
	\xi' \in H \cap (C \cap Y^\ast) , \quad
	\eta' \in H^\ast \cap (C \cap Y^\ast) ,
\]
none of which is an end of $\tilde s$ or of $\tilde t$.

\begin{lemma} \label{lem:periodic-distinct}
	Suppose that $\delta^m(S) = S$ with $m \ge 1$.
	For $n \ge 1$, there is a lift $\tilde s_n$ of $s$ whose halfspaces have Cantor parts $\tilde\delta^{mn}(H)$ and $\tilde\delta^{mn}(H^\ast)$.
\end{lemma}

\begin{proof}
	For $n \ge 1$, $\tilde\delta^{mn}(\tilde s)$ is a lift of $\delta^{mn}(s)$, which is in the class $S$ of $s$ since $\delta^m(S) = S$.
	By Lemma~\ref{lem:qi-extends-cantor}, its halfspaces have Cantor parts $\tilde\delta^{mn}(H)$ and $\tilde\delta^{mn}(H^\ast)$.
	If $s$ is a sphere, then $\delta^{mn}(s)$ is isotopic to $s$ by Theorem~\ref{thm:laudenbach}\eqref{item:laudenbach:spheres}, and the lifted isotopy carries $\tilde\delta^{mn}(\tilde s)$ to a lift $\tilde s_n$ of $s$ with the same Cantor parts, by Lemma~\ref{lem:extend-isotopy}.
	If $s$ is a torus, the correspondence of Lemma~\ref{lem:tree-sides} carries $\tilde\delta^{mn}(\tilde s)$ to a lift $\tilde s_n$ of $s$ with the same ends and sides.
	Since the Cantor part of a halfspace is its side together with the ends, $\tilde s_n$ has the same Cantor parts of halfspaces.
\end{proof}

\begin{lemma} \label{lem:periodic-distinct-distinct}
	Suppose that $\delta^m(S) = S$ with $m \ge 1$, and for $n \ge 1$ let $\tilde s_n$ be a lift of $s$ whose halfspaces have Cantor parts $\tilde\delta^{mn}(H)$ and $\tilde\delta^{mn}(H^\ast)$, which exists by Lemma~\ref{lem:periodic-distinct}.
	Then the lifts $\tilde s_n$ are pairwise distinct.
\end{lemma}

\begin{proof}
	Suppose that $\tilde s_n = \tilde s_{n'}$ for some $n < n'$.
	Then $\tilde\delta^{m'}$ permutes $\{ H, H^\ast \}$ for $m' := m(n' - n)$, and $H \ne H^\ast$ since $H \cap H^\ast$ is finite.
	So $\tilde\delta^{2m'}$ fixes $H$ and $H^\ast$.
	By Lemma~\ref{lem:twist-attracts}, $\tilde\delta^{2m'j}(\xi) \to w^{+\infty}$ and $\tilde\delta^{2m'j}(\eta) \to w^{+\infty}$ as $j \to \infty$.
	Since $H$ and $H^\ast$ are closed, $w^{+\infty} \in H \cap H^\ast$, contradicting that $H \cap H^\ast$ misses $w^{\pm\infty}$.
\end{proof}

\begin{proposition} \label{prop:twist-periodic}
	Suppose that $N(S,T) > 0$.
	Then $\delta_{\vec T}^{\,m}(S) \ne S$ for every $m \ne 0$.
\end{proposition}

\begin{proof}
	Suppose that $\delta^m(S) = S$ with $m \ne 0$; replacing $m$ by $-m$, we may assume $m \ge 1$.
	Let $\tilde s_n$ be the lifts of Lemma~\ref{lem:periodic-distinct}, and $Q_n$ the halfspace of $\tilde s_n$ not containing a base point $x_0$.
	The lifts $\tilde s_n$ are pairwise distinct by Lemma~\ref{lem:periodic-distinct-distinct} and the lifts of $s$ form a locally finite family, so they leave every compact set.
	By Lemma~\ref{lem:far-halfspace-limit}, after passing to a subsequence, the Cantor parts $C \cap Q_n$ converge to a point $\zeta \in C$.
	Passing to a further subsequence, either $C \cap Q_n = \tilde\delta^{mn}(H)$ for all $n$ or $C \cap Q_n = \tilde\delta^{mn}(H^\ast)$ for all $n$; write $H^\epsilon$ for the set $H$ or $H^\ast$ that occurs.

	Let $\xi_\epsilon$ be the one of the points $\xi \in H \cap (C \cap Y)$ and $\eta \in H^\ast \cap (C \cap Y)$ that lies in $H^\epsilon$, and let $\xi'_\epsilon$ be the one of the points $\xi' \in H \cap (C \cap Y^\ast)$ and $\eta' \in H^\ast \cap (C \cap Y^\ast)$ that lies in $H^\epsilon$; none of these four points is an end of $\tilde s$ or of $\tilde t$.
	By Lemma~\ref{lem:twist-attracts}, $\tilde\delta^{mn}(\xi_\epsilon) \to w^{+\infty}$, and it lies in $C \cap Q_n$; so $\zeta = w^{+\infty}$.
	By Lemma~\ref{lem:twist-repels}, $\tilde\delta^{mn}(\xi'_\epsilon)$ stays in a closed set not containing $w^{+\infty}$, while lying in $C \cap Q_n \to \{ w^{+\infty} \}$.
	This is a contradiction.
\end{proof}

The next lemma produces crossing pairs from intersections of a core with a surface.
It is used for the dual sphere of a core below, and in Section~\ref{sec:fellow} for two tori.

\begin{lemma} \label{lem:dual-sphere-crosses}
	Let $t$ be a primitive torus with class $T$ and core $\gamma_T$, and $u$ an essential sphere or a primitive torus whose class $U$ is core-independent from $T$.
	Let $\tilde t$ be a lift of $t$, and $\tilde u$ a lift of $u$ penetrated by the core of $\tilde t$.
	Then $\tilde u$ and $\tilde t$ cross.
\end{lemma}

\begin{proof}
	Let $\langle y \rangle$ be the stabilizer of $\tilde t$, with sides $P$ and $P^\ast$, so that the core of $\tilde t$ has the ends $y^{\pm\infty}$.
	The lifts $\tilde t$ and $\tilde u$ have no common end.
	Indeed, a sphere lift has no ends, and if $u$ is a torus, a common end would give a common nonzero power of $y$ and of a generator of the stabilizer of $\tilde u$ by Lemma~\ref{lem:common-end}, contradicting Lemma~\ref{lem:joint-stabilizer}.

	Since the core penetrates $\tilde u$, the sides $Z$ and $Z^\ast$ of $\tilde u$ can be named so that $y^{+\infty} \in Z$ and $y^{-\infty} \in Z^\ast$.
	By Lemma~\ref{lem:ends-side}, $Z \cap P$ and $Z \cap P^\ast$ contain points that are not ends.
	The same lemma with $y^{-1}$ in place of $y$, which generates the same stabilizer and has $(y^{-1})^{+\infty} = y^{-\infty}$, shows that $Z^\ast \cap P$ and $Z^\ast \cap P^\ast$ contain points that are not ends.
	So $\tilde u$ and $\tilde t$ cross by Lemma~\ref{lem:crossing-criterion}.
\end{proof}

\begin{corollary} \label{lem:dual-sphere-crosses-positive}
	Let $t$ be a primitive torus with class $T$ and core $\gamma_T$, and $u$ an essential sphere or a primitive torus whose class $U$ is core-independent from $T$.
	If $i(\gamma_T, U) > 0$, then $N(U, T) > 0$.
\end{corollary}

\begin{proof}
	Let $\tilde t$ be a lift of $t$ with stabilizer $\langle y \rangle$; its core $\{ y^{\pm\infty} \}$ is a lift of $\gamma_T$.
	It is a core of no lift of $u$: lifts of spheres have trivial stabilizers, and if $u$ is a torus this follows from Lemma~\ref{lem:joint-stabilizer}.
	If $i(\gamma_T, U) > 0$, it penetrates some lift $\tilde u$ of $u$ by Proposition~\ref{prop:orbit-count}, and $(\tilde u, \tilde t)$ is a crossing pair by Lemma~\ref{lem:dual-sphere-crosses}.
	Since crossing is symmetric in the two lifts, $N(U, T) = N(T, U)$, and both are positive.
\end{proof}

\begin{corollary} \label{cor:twist-infinite-order}
	Let $\vec T$ be an oriented class of primitive tori.
	Then $\delta_{\vec T}$ has infinite order in $\Out(\FF)$.
\end{corollary}

\begin{proof}
	The core $\gamma_T$ of $T$ is the class of a primitive curve, so $i(\gamma_T, S_0) = 1$ for some class $S_0$ of essential spheres, by the definition in Section~\ref{subsec:primitive}.
	For instance, let $\{ a, b_2, \dots, b_g \}$ be a free basis of $\FF$ with $a$ representing $\gamma_T$.
	One may take for $S_0$ the class of the nonseparating sphere dual to this basis element.
	Its dual tree is the Bass--Serre tree of the HNN splitting of $\FF$ with vertex group $\langle b_2, \dots, b_g \rangle$ and stable letter $a$, on which $a$ translates by $1$.
	A class of spheres and a class of tori are core-independent, so $N(S_0, T) > 0$ by Corollary~\ref{lem:dual-sphere-crosses-positive}.
	By Proposition~\ref{prop:twist-periodic}, $\delta_{\vec T}^{\,m}(S_0) \ne S_0$ for every $m \ne 0$, so $\delta_{\vec T}^{\,m} \ne 1$.
\end{proof}

\begin{corollary} \label{cor:dehn_disj}
	Let $\vec T$ be an oriented class of primitive tori, with underlying class $T$, and $S$ a class of essential spheres or of primitive tori that is core-independent from $T$.
	If $\delta_{\vec T}^{\,m}(S) = S$ for some $m \ne 0$, then $N(S,T) = 0$.
\end{corollary}

\begin{proof}
	This is the contrapositive of Proposition~\ref{prop:twist-periodic}.
\end{proof}

The proof of Proposition~\ref{prop:coincide} rests on a theorem of Cohen and Lustig on efficient Dehn twists on graphs of groups \cite[Theorem~1.1]{MR1691946}.
We first identify our Dehn twists with the Dehn twists of Cohen and Lustig, and then check efficiency.
Let $t$ be a primitive torus with oriented core, $o := o_t$, and $\mathcal{G} := X_t / \FF$ the quotient graph of groups, with $\pi_1 \mathcal{G} = \FF$ given by the action.
The \emph{twistors} $z_e$ are the images of $o(\vec e)$ in the edge groups; they satisfy $z_{\bar e} = z_e^{-1}$ since $o(\ast \vec e) = o(\vec e)^{-1}$.

\begin{lemma} \label{lem:coincide-twist}
	% Let $t$, $o$, $\mathcal{G}$ and $(z_e)$ be as above.
	$\delta_{\vec X_t}$ is the outer class of the Dehn twist $D(\mathcal{G}, (z_e))$ of \cite[Definition~5.1]{MR1691946}.
\end{lemma}

\begin{proof}
	The Dehn twist $D(\mathcal{G}, (z_e))$ of \cite[Definition~5.1]{MR1691946} is the automorphism of $\pi_1 \mathcal{G}$ that fixes the vertex groups and sends each edge letter $t_e$ to $t_e z_e$.
	Describe the Bass--Serre tree by $\mathcal{G}$-paths from a base vertex $v_0$, and let $[p]$ be the vertex at the end of a path $p = g_0 t_{e_1} g_1 \dots t_{e_k} g_k$ with edges $\vec e_1, \dots, \vec e_k$.
	Then $D(\mathcal{G}, (z_e))$ carries $[p]$ to $[g_0 t_{e_1} z_{e_1} g_1 \dots t_{e_k} z_{e_k} g_k] = o(\vec e_1) \dots o(\vec e_k) \, [p]$, where $o(\vec e_i) = p_i z_{e_i} p_i^{-1}$ for the initial segment $p_i$ of $p$ ending with $t_{e_i}$.
	So by Lemma~\ref{lem:twist-path} it induces $\delta_{\vec X_t, v_0}$ on $X_t$.
	This is the Dehn twist of the graph of groups $\mathcal{G}$ in the sense of Cohen and Lustig \cite{MR1341810}, computed on the Bass--Serre tree of \cite{MR1341810}.
	So $\delta_{\vec X_t}$ is the outer class of $D(\mathcal{G}, (z_e))$.
\end{proof}

\begin{lemma} \label{lem:coincide-efficient}
	% Let $t$, $o$, $\mathcal{G}$ and $(z_e)$ be as above.
	The Dehn twist $D(\mathcal{G}, (z_e))$ is efficient in the sense of \cite[Definition~6.2]{MR1691946}.
\end{lemma}

\begin{proof}
	We name the five conditions of \cite[Definition~6.2]{MR1691946} as we verify them.
	Let $\langle y \rangle$, $y := o(\vec e)$ for an oriented edge $\vec e$, be the edge group, contained in a vertex group $G_v$.

	The twistor of the edge is nontrivial, so there is no unused edge.

	There are no proper powers: if $r \in G_v$ and $r^p \in \langle y \rangle$ for some $p \ne 0$, then $r \in \langle y \rangle$.
	Indeed, if $r^p = y^q$ with $q = 0$, then $r = 1$.
	If $q \ne 0$, then $r$ and $y$ are powers of a common element by Theorem~\ref{thm:free}\eqref{item:free:commute}.
	Since $y$ is not a proper power by Lemma~\ref{lem:primitive-basis} and Theorem~\ref{thm:free}\eqref{item:free:primitive}, this element is $y^{\pm 1}$ and $r \in \langle y \rangle$.

	The graph of groups $\mathcal{G} = X_t / \FF$ is minimal: a vertex of valence one of $X_t / \FF$ whose edge group maps onto its vertex group would lift to a vertex of valence one of $X_t$, and $X_t$ has none by Lemma~\ref{lem:no-leaf}.

	If the quotient graph is a loop, that is, $\FF = A \ast_{\langle y \rangle}$ with stable letter $r$ and $r y r^{-1} \in A$, then no nonzero powers of $y$ and of $r y r^{-1}$ are conjugate in $A$.
	Indeed, let $c y^{m} c^{-1} = r y^{n} r^{-1}$ with $c \in A$ and $m, n \ne 0$.
	Then $m = n$, since conjugation preserves translation length on $\mathrm{Cay}(\FF)$ and $y^m$ is not conjugate to its inverse.
	So $r^{-1} c$ commutes with $y^m$, and $r^{-1} c \in \langle y \rangle \le A$ by Theorem~\ref{thm:free}\eqref{item:free:commute} and~\eqref{item:free:primitive}.
	Hence $r \in A$, which is false.
	So the two ends of the edge at the vertex are not positively bonded.

	Finally, no vertex is invisible.
	A vertex of valence two with both edge groups mapping onto its vertex group occurs only if the quotient graph is a loop with cyclic vertex group $A$.
	Then $\FF$ would be $\ZZ^2$ or the fundamental group of the Klein bottle.
	If the quotient graph is a segment, its vertices have valence one, and the last two conditions are vacuous.
\end{proof}

\begin{proof}[Proof of Proposition~\ref{prop:coincide}]
	If $\vec T = \vec T'$, then $\delta_{\vec T} = \delta_{\vec T'}$, since $\delta_{\vec T}$ depends only on $\vec T$ by Section~\ref{subsec:twist-homeo}.

	Conversely, suppose that $\delta_{\vec T} = \delta_{\vec T'}$ in $\Out(\FF)$, and let $\vec t$ and $\vec t'$ be tori with oriented cores in $\vec T$ and $\vec T'$.
	By Lemma~\ref{lem:gul}, these are the Dehn twists along the oriented one-edge $\ZZ$-splittings $\vec X_t$ and $\vec X_{t'}$; we write $o := o_t$ and $o' := o_{t'}$.
	We use the theorem of Cohen and Lustig \cite[Theorem~1.1]{MR1691946}.
	Let $D_1$, $D_2$ be efficient Dehn twists on graphs of groups $\mathcal{G}_1$, $\mathcal{G}_2$, in the sense of \cite[Definition~6.2]{MR1691946}, and $h \colon \pi_1 \mathcal{G}_1 \to \pi_1 \mathcal{G}_2$ an isomorphism with $h D_1 h^{-1} = D_2$ in $\Out(\pi_1 \mathcal{G}_2)$.
	Then there is an isomorphism of graphs of groups $H \colon \mathcal{G}_1 \to \mathcal{G}_2$ that induces $h$ up to an inner automorphism and carries each twistor $z_e$ to $z_{H(e)}$.

	We apply this to $\mathcal{G}_1 = X_t / \FF$ and $\mathcal{G}_2 = X_{t'} / \FF$, with $\pi_1 \mathcal{G}_i = \FF$ given by the actions, with the twistors $z_e$ defined before Lemma~\ref{lem:coincide-twist} from $o$ and $o'$, and with $h = \mathrm{id}$.
	The identification $\pi_1 \mathcal{G}_i = \FF$ is canonical only up to inner automorphisms, which is why the theorem is used with $h$ determined up to inner automorphisms.
	By Lemma~\ref{lem:coincide-twist}, $\delta_{\vec X_t}$ is the outer class of $D(\mathcal{G}_1, (z_e))$, and $\delta_{\vec X_{t'}}$ that of $D(\mathcal{G}_2, (z_e))$.
	Both twists are efficient by Lemma~\ref{lem:coincide-efficient}.

	So $H$ induces an isomorphism $X_t \to X_{t'}$ that is equivariant with respect to conjugation by some $c \in \FF$.
	Composing it with $c^{-1}$ gives an $\FF$-equivariant isomorphism, which carries $o$ to $o'$ since $H$ carries $z_e$ to $z_{H(e)}$.
	By the definition of the oriented class in Section~\ref{subsec:classes}, $\vec T = \vec T'$.
\end{proof}

\subsection{Commuting Dehn twists} \label{subsec:commuting}

% The goal of this subsection is Theorems~\ref{thm:commuting} and~\ref{thm:commuting-crossing}: two Dehn twists along core-independent classes of primitive tori commute if the classes have disjoint representatives, and only if their crossing number vanishes.
If two homeomorphisms have disjoint supports, they commute, and so do their mapping classes.
On surfaces, Dehn twists commute if and only if they have representatives with disjoint supports \cite[Fact~3.9]{MR2850125}.
For primitive tori in $\MM$, we prove that disjoint representatives give commuting twists, and that commuting twists have crossing number $N(T,T') = 0$ in the sense of Section~\ref{sec:crossing}.

\begin{theorem} \label{thm:commuting}
	Let $T$ and $T'$ be core-independent classes of primitive tori, and $\vec T$, $\vec T'$ oriented classes over $T$ and $T'$.
	If $T$ and $T'$ have disjoint representatives, then $\delta_{\vec T}$ and $\delta_{\vec T'}$ commute.
\end{theorem}

\begin{proof}
	Let $t \in T$ and $t' \in T'$ be disjoint, with regular neighborhoods $N(t)$ and $N(t')$ chosen disjoint, and with oriented cores and orientations such that $\vec t$ and $\vec t'$ lie in the oriented classes $\vec T$ and $\vec T'$.
	Then $\delta_{\vec t}$ and $\delta_{\vec t'}$ have disjoint supports, so they commute, and so do their images $\delta_{\vec T}$ and $\delta_{\vec T'}$.
\end{proof}

\begin{theorem} \label{thm:commuting-crossing}
	Let $T$ and $T'$ be core-independent classes of primitive tori, and $\vec T$, $\vec T'$ oriented classes over $T$ and $T'$.
	If $\delta_{\vec T}$ and $\delta_{\vec T'}$ commute, then $N(T, T') = 0$.
\end{theorem}

\begin{proof}
	Suppose that $\delta_{\vec T}$ and $\delta_{\vec T'}$ commute, and let $\vec t$ and $\vec t'$ be tori with oriented cores in $\vec T$ and $\vec T'$.
	The homeomorphism $h := \delta_{\vec t'}$ conjugates $\delta_{\vec t}$ to the Dehn twist homeomorphism along $h(\vec t)$, with the regular neighborhood $h(N(t))$ and the transported parametrization.
	This parametrization satisfies the conventions of Section~\ref{subsec:twist-homeo} for the transported orientations of $h(t)$ and of its core, since $h$ preserves the orientation of $\MM$.
	By Lemma~\ref{lem:gul}, the image of $h \delta_{\vec t} h^{-1}$ in $\Out(\FF)$ is the Dehn twist along the oriented class $\delta_{\vec T'}(\vec T)$, and it equals $\delta_{\vec T'} \delta_{\vec T} \delta_{\vec T'}^{-1} = \delta_{\vec T}$.
	By Proposition~\ref{prop:coincide}, $\delta_{\vec T'}(\vec T) = \vec T$, so in particular $\delta_{\vec T'}(T) = T$.
	By Corollary~\ref{cor:dehn_disj}, $N(T, T') = 0$.
\end{proof}

\section{Growth under multitwists} \label{sec:growth}

The goal of this section is Proposition~\ref{prop:twist-growth}: under a multitwist along disjoint primitive tori, the intersection number of a curve $\alpha$ with a sphere $S$ grows linearly in the twisting numbers.
The slope in the $v$-th twisting number is $i(\alpha, T_v) \, i(\gamma_v, S)$, where $\gamma_v$ is the core of $T_v$.
As a consequence, Corollary~\ref{cor:disjoint_twists_free_abelian} shows that the Dehn twists along disjoint core-independent primitive tori generate a free abelian group, and Proposition~\ref{prop:pairwise-free-abelian} extends this to classes with pairwise disjoint representatives.
The estimate rests on a count of lifts, Lemma~\ref{lem:block-count}, which is also used in Section~\ref{sec:fellow}.
The lift of a power of a Dehn twist moves one of three consecutive lifts of the torus far along the middle one, and a curve joining the outer two must then penetrate many lifts of a third surface.

\subsection{Counting penetrated lifts} \label{subsec:spirals}

Let $u$ be an essential sphere or a primitive torus, with class $U$ and dual tree $X_u$.
Let $\tilde t_2$ be a lift of a primitive torus $t_2$ whose class is core-independent from $U$ if $u$ is a torus, with stabilizer $\langle w \rangle$ and core $\tilde\gamma'$ of class $\gamma'$.
Then $\tilde\gamma'$ is a core of no lift of $u$: for a torus $u$ by Lemma~\ref{lem:joint-stabilizer}, and for a sphere $u$ because sphere lifts have no cores.
Put $\tau := i(\gamma', U)$ and suppose that $\tau \ge 1$.
Then $w$ is hyperbolic on $X_u$ with translation length $\tau$ by Proposition~\ref{prop:orbit-count-translation}, and we let $A$ be its axis.

For a lift $\tilde\alpha$ of a curve class that is a core of no lift of $u$, with minimal subtree $B$ on $X_u$, let $\pi_A(B)$ be the nearest-point projection of $B$ to $A$.
It is $A \cap B$ if this is nonempty, and a single vertex otherwise.
For an edge $e \in A \setminus B$, the minimal subtree $B$ lies in the component of $X_u \setminus e$ containing $\pi_A(B)$, since $B$ is connected and is joined to $\pi_A(B)$ by a path missing $e$.
We say that a lift $\tilde u_e$ of $u$ \emph{separates} two lifts of curves if the ends of one lie on one side of $\tilde u_e$ and the ends of the other on the other side.
By Lemma~\ref{lem:subtree-sides}, $\tilde u_e$ separates $\tilde\alpha$ from $\tilde\beta$ if and only if $e$ lies on neither minimal subtree and the two minimal subtrees lie in different components of $X_u \setminus e$.

For the next two lemmas, let $\tilde\alpha$ and $\tilde\beta$ be lifts of curve classes, cores of no lift of $u$, whose representing elements do not commute with $w$, with minimal subtrees $B_\alpha$ and $B_\beta$ on $X_u$.

\begin{lemma} \label{lem:separating-count}
	$A \cap B_\alpha$ and $A \cap B_\beta$ are finite.
\end{lemma}

\begin{proof}
	Let $g_\alpha$ represent $\tilde\alpha$.
	If $A \cap B_\alpha$ were infinite, then $B_\alpha$ would be the axis of $g_\alpha$ and would share a ray with $A$.
	Let $p$ be the translation length of $g_\alpha$.
	Then for a suitable sign the element $w^{p} g_\alpha^{\mp \tau}$ would fix every edge far out on this ray.

	Edge stabilizers of $X_u$ are trivial if $u$ is a sphere, and distinct edges of $X_u$ have stabilizers with trivial intersection by Lemma~\ref{lem:end-stabilizer} if $u$ is a torus.
	So $w^{p} = g_\alpha^{\pm \tau}$, and $w$ and $g_\alpha$ would commute by Theorem~\ref{thm:free}\eqref{item:free:commute}.
	Hence $A \cap B_\alpha$ is finite, and likewise $A \cap B_\beta$.
\end{proof}

\begin{lemma} \label{lem:separating-count-bound}
	Let $\lambda$ be the number of edges of $A$ in $B_\alpha$ plus the number of edges of $A$ in $B_\beta$, and let $b$ be the number of edges of $A$ separating $\tilde\alpha$ from $\tilde\beta$.
	Then for $n \ge 1$, at least $n \tau - \lambda - b$ lifts of $u$ are penetrated by $\tilde\gamma'$ and separate $\tilde\alpha$ from $w^n \tilde\beta$.
\end{lemma}

\begin{proof}
	Index the vertices of $A$ by $\ZZ$ so that $w$ acts as translation by $\tau$.
	Let $\pi_A(B_\alpha)$ run from $p$ to $p'$ and $\pi_A(B_\beta)$ from $q$ to $q'$, with $p \le p'$ and $q \le q'$.
	The minimal subtree of $w^n \tilde\beta$ is $w^n B_\beta$, with projection running from $q + n\tau$ to $q' + n\tau$.

	Every edge of $A$ between $p'$ and $q + n\tau$ lies on neither minimal subtree and has $\pi_A(B_\alpha)$ on one side and $w^n \pi_A(B_\beta)$ on the other.
	So it separates $\tilde\alpha$ from $w^n \tilde\beta$, and it is penetrated by $\tilde\gamma'$ by Lemma~\ref{lem:subtree-sides-penetrate}.
	There are $q + n\tau - p'$ such edges, and
	\[
		p' - q = (p' - p) + (q' - q) + (p - q') \le \lambda + b ,
	\]
	since $p' - p$ and $q' - q$ count edges of $A$ in $B_\alpha$ and in $B_\beta$, and the edges of $A$ strictly between $q'$ and $p$, if $p > q'$, separate $\tilde\beta$ from $\tilde\alpha$.
\end{proof}

Now let $\tilde t_1$ be a lift of a primitive torus $t_1$ whose class is core-independent from $U$ if $u$ is a torus, such that $\tilde t_1$ and $\tilde t_2$ are disjoint and their stabilizers have no common nontrivial element.
Let $P$ be the side of $\tilde t_1$ not containing the ends of $\tilde t_2$, which is well defined by Section~\ref{subsec:twist-lift}.

\begin{lemma} \label{lem:far-penetrate}
	Let $\tilde u$ be a lift of $u$ penetrated by $\tilde\gamma'$ that does not cross $\tilde t_1$.
	Then the ends of $\tilde t_1$ lie in one side of $\tilde u$.
\end{lemma}

\begin{proof}
	The lifts $\tilde t_1$ and $\tilde u$ have no common end: for a torus $u$ this follows from core-independence by Lemmas~\ref{lem:common-end} and~\ref{lem:joint-stabilizer}, and a sphere lift has no ends.
	Suppose that the ends of $\tilde t_1$ lay in different sides of $\tilde u$.
	Then by Lemma~\ref{lem:ends-side} each of the four intersections of a side of $\tilde u$ with a side of $\tilde t_1$ would contain a point that is not an end, so $\tilde u$ and $\tilde t_1$ would cross by Lemma~\ref{lem:crossing-criterion}.
	So the ends of $\tilde t_1$ lie in one side of $\tilde u$.
\end{proof}

\begin{lemma} \label{lem:far-penetrate-side}
	Let $\tilde u$ be a lift of $u$ penetrated by $\tilde\gamma'$ that does not cross $\tilde t_1$, and $Z$ the side of $\tilde u$ containing the ends of $\tilde t_1$, which exists by Lemma~\ref{lem:far-penetrate}.
	Then every point of $P$ lies in $Z$ or is an end of $\tilde u$.
\end{lemma}

\begin{proof}
	The lifts $\tilde t_1$ and $\tilde u$ have no common end: for a torus $u$ this follows from core-independence by Lemmas~\ref{lem:common-end} and~\ref{lem:joint-stabilizer}, and a sphere lift has no ends.
	Let $Z^\ast$ be the other side of $\tilde u$.
	Since $\tilde\gamma'$ penetrates $\tilde u$, one end of $\tilde t_2$ lies in $Z^\ast$.
	It lies on the side of $\tilde t_1$ other than $P$, and it is an end neither of $\tilde t_1$ nor of $\tilde u$.
	So Lemma~\ref{lem:far-sides}, applied to $\tilde t_1$ and $\tilde u$, gives the claim.
\end{proof}

Let $\tilde t_3$ be a lift of a primitive torus $t_3$ whose class is core-independent from $U$ if $u$ is a torus.
Suppose that $\tilde t_1$, $\tilde t_2$, $\tilde t_3$ are pairwise disjoint and that their stabilizers pairwise have no common nontrivial element.
Then no two of $\tilde t_1$, $\tilde t_2$, $\tilde t_3$ have a common end, by Lemma~\ref{lem:common-end}.
Let $R$ be the side of $\tilde t_3$ not containing the ends of $\tilde t_2$.
For $n \ge 1$, the lift $w^n \tilde t_3$ stands in the same relation to $\tilde t_2$ as $\tilde t_3$.
It is disjoint from $w^n \tilde t_2 = \tilde t_2$, its stabilizer is the conjugate of the stabilizer of $\tilde t_3$ by $w^n$, and $w^n R$ is its side not containing the ends of $\tilde t_2$.

\begin{lemma} \label{lem:block-separate}
	Let $X$ be $\tilde t_3$ or $w^n \tilde t_3$ for some $n \ge 1$, and $R_X$ its side not containing the ends of $\tilde t_2$.
	Let $\tilde u$ be a lift of $u$ penetrated by $\tilde\gamma'$ that crosses neither $\tilde t_1$ nor $X$, such that the ends of $\tilde t_1$ and the ends of $X$ lie in different sides of $\tilde u$.
	Then every lift of a curve class with one end in $P$ and one end in $R_X$ penetrates $\tilde u$ or is a core of $\tilde u$.
\end{lemma}

\begin{proof}
	By Lemmas~\ref{lem:far-penetrate} and~\ref{lem:far-penetrate-side}, the ends of $\tilde t_1$ lie in a side $Z$ of $\tilde u$, and every point of $P$ lies in $Z$ or is an end of $\tilde u$.
	The lift $X$ satisfies the hypotheses imposed on $\tilde t_1$, with $R_X$ in place of $P$: it is a lift of the primitive torus $t_3$, whose class is core-independent from $U$ if $u$ is a torus, it is disjoint from $\tilde t_2$, and its stabilizer, a conjugate of the stabilizer of $\tilde t_3$ by a power of $w$, has no common nontrivial element with $\langle w \rangle$.
	So by the same lemmas, the ends of $X$ lie in one side of $\tilde u$, which is the other side $Z^\ast$ by hypothesis, and every point of $R_X$ lies in $Z^\ast$ or is an end of $\tilde u$.
	Let $\tilde a$ be a lift of a curve class with one end in $P$ and one end in $R_X$, and suppose that it is not a core of $\tilde u$.
	Then no end of $\tilde a$ is an end of $\tilde u$ by Section~\ref{subsec:curve-lifts}, so one end lies in $Z$ and the other in $Z^\ast$, and $\tilde a$ penetrates $\tilde u$.
\end{proof}

Lemma~\ref{lem:block-count} is the count of lifts announced at the start of this section, and its statement collects the hypotheses of this subsection.

\begin{lemma} \label{lem:block-count}
	Let $u$ be an essential sphere or a primitive torus with class $U$.
	Let $\tilde t_1$, $\tilde t_2$, $\tilde t_3$ be pairwise disjoint lifts of primitive tori $t_1$, $t_2$, $t_3$, whose classes are core-independent from $U$ if $u$ is a torus, and whose stabilizers pairwise have no common nontrivial element.
	Let $\langle w \rangle$ be the stabilizer of $\tilde t_2$, let $\tilde\gamma'$ be the core of $\tilde t_2$, of class $\gamma'$, and suppose that $\tau := i(\gamma', U) \ge 1$.
	Let $P$ be the side of $\tilde t_1$, and $R$ the side of $\tilde t_3$, not containing the ends of $\tilde t_2$.
	Let $T$ and $T''$ be the classes of $t_1$ and $t_3$, with cores $\gamma_T$ and $\gamma_{T''}$, and put
	\[
		\theta := \max \bigl( i(\gamma_T, U), \, i(\gamma_{T''}, U) \bigr) , \qquad
		\iota := \max \bigl( N(U, T), \, N(U, T'') \bigr) .
	\]
	Let $\tilde a$ be a lift of a curve class $a$ with one end in $P$ and one end in $R$, let $n \ge 1$, and let $\tilde c$ be a lift of a curve class with one end in $P$ and one end in $w^n R$.
	Then at least
	\[
		n \tau - 4 \tau \theta - 2 \tau \, i(a, U) - 8 \tau \iota - 2
	\]
	lifts of $u$ are penetrated by both $\tilde c$ and $\tilde\gamma'$.
\end{lemma}

\begin{proof}
	We shall apply Lemma~\ref{lem:separating-count-bound} to the cores of $\tilde t_1$ and $\tilde t_3$, and first check its hypotheses.
	Let $\tilde\gamma_T$ and $\tilde\gamma_{T''}$ be the cores of $\tilde t_1$ and $\tilde t_3$.
	They are cores of no lift of $u$, and so is $\tilde\gamma'$, by the core-independence hypotheses.
	The elements representing $\tilde\gamma'$, $\tilde\gamma_T$, $\tilde\gamma_{T''}$ pairwise do not commute, since the stabilizers have no common nontrivial element.
	The element $g_a$ representing $\tilde a$ does not commute with $w$: otherwise $\tilde a = \tilde\gamma'$, whose ends are not in $P$.

	Let $A$ be the axis of $w$ on the dual tree $X_u$ of $u$, on which $w$ is hyperbolic with translation length $\tau$ by Proposition~\ref{prop:orbit-count-translation}.
	We apply Lemma~\ref{lem:separating-count-bound} with $\tilde\alpha = \tilde\gamma_T$ and $\tilde\beta = \tilde\gamma_{T''}$, whose representing elements do not commute with $w$, and bound the quantities $\lambda$ and $b$ of Lemma~\ref{lem:separating-count-bound} in turn, beginning with $\lambda \le 4 \tau \theta$.
	The edges of $A$ in the minimal subtree of $\tilde\gamma_T$ are lifts of $u$ penetrated by $\tilde\gamma'$ and by $\tilde\gamma_T$.
	By Proposition~\ref{prop:orbit-count}, these meet $\tau$ orbits of $\langle w \rangle$ and $i(\gamma_T, U)$ orbits of the element representing $\tilde\gamma_T$, so there are at most $2 \tau \theta$ such lifts by Lemma~\ref{lem:orbits-bound}.
	Likewise for $\tilde\gamma_{T''}$, so $\lambda \le 4 \tau \theta$.

	Next we show that $b \le 2 \tau \, i(a, U) + 1 + 4 \tau \iota$.
	Let $e \in A$ separate $\tilde\gamma_T$ from $\tilde\gamma_{T''}$; the lift $\tilde u_e$ is penetrated by $\tilde\gamma'$ by Lemma~\ref{lem:subtree-sides-penetrate}.
	Suppose that $\tilde u_e$ crosses neither $\tilde t_1$ nor $\tilde t_3$.
	Since $\tilde u_e$ separates the cores of $\tilde t_1$ and $\tilde t_3$, whose ends are those of $\tilde t_1$ and $\tilde t_3$, the lift $\tilde a$ penetrates $\tilde u_e$ or is a core of $\tilde u_e$ by Lemma~\ref{lem:block-separate} with $X = \tilde t_3$.
	Hence each edge of $A$ separating $\tilde\gamma_T$ from $\tilde\gamma_{T''}$ is penetrated by both $\tilde\gamma'$ and $\tilde a$, or has $\tilde a$ as a core, or crosses $\tilde t_1$, or crosses $\tilde t_3$.

	If $\tilde a$ is a core of no lift of $u$, there are at most $2 \tau \, i(a, U)$ edges of the first kind by Proposition~\ref{prop:orbit-count} and Lemma~\ref{lem:orbits-bound}, and none of the second.
	Otherwise there are none of the first kind by Lemma~\ref{lem:core-no-penetrate}, and at most one of the second.
	The lifts of $u$ crossing $\tilde t_1$ meet $N(U, T)$ orbits of the stabilizer of $\tilde t_1$ by Lemma~\ref{lem:count-crossing}, since $U$ and $T$ are core-independent.
	So at most $2 \tau \iota$ edges of $A$ are of the third kind by Lemma~\ref{lem:orbits-bound}, and likewise at most $2 \tau \iota$ of the fourth.
	So $b \le 2 \tau \, i(a, U) + 1 + 4 \tau \iota$.
	Lemma~\ref{lem:separating-count-bound} now gives at least $n \tau - 4 \tau \theta - 2 \tau \, i(a, U) - 1 - 4 \tau \iota$ lifts of $u$ penetrated by $\tilde\gamma'$ and separating $\tilde\gamma_T$ from $w^n \tilde\gamma_{T''}$.

	It remains to see that all but at most $4 \tau \iota + 1$ of these lifts are penetrated by $\tilde c$.
	Let $\tilde u$ be one of these lifts.
	Suppose that $\tilde u$ crosses neither $\tilde t_1$ nor $w^n \tilde t_3$.
	By the choice of $\tilde u$, it separates the cores $\tilde\gamma_T$ of $\tilde t_1$ and $w^n \tilde\gamma_{T''}$ of $w^n \tilde t_3$, whose ends are those of $\tilde t_1$ and $w^n \tilde t_3$.
	So the lift $\tilde c$ penetrates $\tilde u$ or is a core of $\tilde u$ by Lemma~\ref{lem:block-separate} with $X = w^n \tilde t_3$, and $\tilde c$ is a core of at most one lift of $u$.
	At most $2 \tau \iota$ of the lifts penetrated by $\tilde\gamma'$ cross $\tilde t_1$, and at most $2 \tau \iota$ cross $w^n \tilde t_3$, a lift of $t_3$, again by Lemmas~\ref{lem:count-crossing} and~\ref{lem:orbits-bound}.
	Subtracting these and the possible core gives the claimed bound.
\end{proof}

\subsection{The growth estimate} \label{subsec:growth}

Let $S$ be a class of essential spheres with representative $s$ and dual tree $X_S$, and let $\tilde\alpha$ be a lift of a curve class $\alpha$ represented by $g$.
By Lemma~\ref{lem:subtree-sides-penetrate} and Proposition~\ref{prop:orbit-count}, $\tilde\alpha$ has a minimal subtree on $X_S$, the lifts of $s$ penetrated by $\tilde\alpha$ are the edges of this subtree, and they form $i(\alpha, S)$ orbits under $\langle g \rangle$.
Thus $i(\alpha, S)$ is the translation length of $g$ on $X_S$, by Proposition~\ref{prop:orbit-count-translation} applied to $u = s$, whose lifts have no cores.
Lemmas~\ref{lem:count-crossing} and~\ref{lem:block-count} are stated for spheres as well as tori, and we apply them with $u = s$.

We need one more fact about segments of an axis.

\begin{lemma} \label{lem:segments}
	Let $\phi$ be a hyperbolic isometry of a tree with axis $B$ and translation length $\lVert \phi \rVert$.
	Let $k \ge 1$ and $c \ge 0$ be integers, and let $\Sigma_j$, $j \in \ZZ$, be segments of $B$ with $\Sigma_{j+k} = \phi \Sigma_j$, each with more than $c$ edges, such that any two of these segments have at most $c$ common edges.
	Then $\lVert \phi \rVert \ge \sum_{j=1}^{k} \bigl( \lvert \Sigma_j \rvert - c \bigr)$, where $\lvert \Sigma_j \rvert$ is the number of edges of $\Sigma_j$.
\end{lemma}

\begin{proof}
	Orient $B$ in the direction of translation by $\phi$, and let $a_j$ and $b_j$ be the initial and terminal vertices of $\Sigma_j$.
	We first show that the initial vertices are distinct and that $a_{j'}$ comes at least $\lvert \Sigma_j \rvert - c$ edges after $a_j$ whenever $a_j$ precedes $a_{j'}$.
	Let $j \ne j'$.
	If $a_j = a_{j'}$, then the two segments share $\min(\lvert \Sigma_j \rvert, \lvert \Sigma_{j'} \rvert) > c$ edges; so the initial vertices are distinct.
	Let $a_j$ precede $a_{j'}$.
	If $b_{j'}$ did not come after $b_j$, then $\Sigma_{j'} \subset \Sigma_j$ and they would share $\lvert \Sigma_{j'} \rvert > c$ edges.
	So $\Sigma_j \cap \Sigma_{j'}$ is the segment from $a_{j'}$ to $b_j$ if $a_{j'}$ precedes $b_j$, and it has at most $c$ edges.
	Hence $a_{j'}$ comes at least $\lvert \Sigma_j \rvert - c$ edges after $a_j$ in all cases.

	Now we sum these spacings along one period of $\phi$.
	Order the segments by their initial vertices.
	Since $\phi$ translates $B$, it preserves this order, and it has exactly $k$ orbits on the segments; so it carries each segment to the $k$-th following segment.
	Let $\Sigma = \Sigma^{(0)}, \Sigma^{(1)}, \dots, \Sigma^{(k)} = \phi \Sigma$ be $k+1$ consecutive segments in this order.
	Then $\Sigma^{(0)}, \dots, \Sigma^{(k-1)}$ represent the $k$ orbits, so they are translates of $\Sigma_1, \dots, \Sigma_k$ in some order, and
	\[
		\lVert \phi \rVert
		= d \bigl( a(\Sigma), a(\phi \Sigma) \bigr)
		= \sum_{i=0}^{k-1} d \bigl( a(\Sigma^{(i)}), a(\Sigma^{(i+1)}) \bigr)
		\ge \sum_{i=0}^{k-1} \bigl( \lvert \Sigma^{(i)} \rvert - c \bigr)
		= \sum_{j=1}^{k} \bigl( \lvert \Sigma_j \rvert - c \bigr) .
		\qedhere
	\]
\end{proof}

Let $t_1, \dots, t_r$ be pairwise disjoint primitive tori with oriented cores, whose classes $T_1, \dots, T_r$ are pairwise core-independent, and let $\gamma_v$ be the core of $T_v$.
For $m_1, \dots, m_r \in \ZZ$, put
\[
	\Delta := \delta_{\vec T_1}^{\,m_1} \dots \delta_{\vec T_r}^{\,m_r} \in \Out(\FF) .
\]
The twists commute, since the homeomorphisms $\delta_{\vec t_v}$ have disjoint supports.
The multitwist $\Delta$ is represented by the homeomorphism $h := \delta_{\vec t_1}^{\,m_1} \circ \dots \circ \delta_{\vec t_r}^{\,m_r}$, whose support is the union of the disjoint regular neighborhoods $N(t_v)$.

Let $\mathbf{t} := t_1 \sqcup \dots \sqcup t_r$.
A \emph{region} of $\mathbf{t}$ is the closure of a component of $\widetilde\MM$ minus the lifts of all the regular neighborhoods $N(t_v)$; for $r = 1$ this is the definition of Section~\ref{subsec:twist-lift}.
The dual tree $X$ of $\mathbf{t}$ has the regions as vertices and the lifts of the $t_v$ as edges, the edge of a lift joining the two regions adjacent to it; for $r = 1$ this is the dual tree $X_{t_1}$ defined in Section~\ref{subsec:curve-lifts}.
It is a tree since $\widetilde\MM$ is simply connected and each lift of a $t_v$ separates it.

Lemmas~\ref{lem:subtree-sides}, \ref{lem:subtree-sides-penetrate} and~\ref{lem:nested-sides} hold for $X$ with the same proofs, the last in the form needed below.
For distinct edges $e \ne e'$ of $X$, let $\kappa$ be the component of $X \setminus e$ containing $e'$; then the ends of the lift at $e'$ lie in $C(Y_\kappa)$.
The first two proofs use only regions and halfspaces.
In the last, the ends of the lift at $e'$ are not ends of the lift at $e$.
This holds by Lemma~\ref{lem:end-stabilizer} if the two lifts are lifts of one torus, and by Lemmas~\ref{lem:common-end} and~\ref{lem:joint-stabilizer} if they are lifts of different tori, since the classes are core-independent.

For the next five lemmas, let $W$ be a region of $\mathbf{t}$, and $\tilde h$ the lift of $h$ that is the identity on $W$.

\begin{lemma} \label{lem:multitwist-lift}
	$\tilde h$ is a deck transformation $h_\Omega$ on each region $\Omega$ of $\mathbf{t}$, with $h_W = 1$.
\end{lemma}

\begin{proof}
	The lift $\tilde h$ is the composition of the lifts of the factors $\delta_{\vec t_v}^{\,m_v}$ that are the identity on the region of $t_v$ containing $W$.
	Each factor is a deck transformation on each region of $t_v$ by Section~\ref{subsec:twist-lift}, so $\tilde h$ is a deck transformation $h_\Omega$ on each region $\Omega$ of $\mathbf{t}$, and $h_W = 1$.
\end{proof}

We write $h_\Omega$ for the deck transformations of Lemma~\ref{lem:multitwist-lift}.
For a region $\Omega$, the lift $h_\Omega^{-1} \tilde h$ of $h$ is the identity on $\Omega$, so it is the composition of the lifts of the factors $\delta_{\vec t_v}^{\,m_v}$ that are the identity on the regions of $t_v$ containing $\Omega$.

\begin{lemma} \label{lem:multitwist-lift-adjacent}
	If $\Omega$ and $\Omega'$ are adjacent across a lift $\tilde t_v$ of $t_v$, then $h_{\Omega'} = h_\Omega u^{\lvert m_v \rvert}$ for a suitable generator $u$ of the stabilizer of $\tilde t_v$.
\end{lemma}

\begin{proof}
	Consider the lift $h_\Omega^{-1} \tilde h$, which is the composition of the lifts of the factors that are the identity on the regions containing $\Omega$.
	For $v' \ne v$, the factor for $t_{v'}$ is the identity on the region of $t_{v'}$ containing $\Omega$, which also contains $\Omega'$.
	By Lemma~\ref{lem:twist-adjacent}, the factor for $t_v$ equals $u^{\lvert m_v \rvert}$ on $\Omega'$ for a suitable generator $u$ of the stabilizer of $\tilde t_v$, and $u$ preserves $\Omega'$; so $h_{\Omega'} = h_\Omega u^{\lvert m_v \rvert}$.
\end{proof}

\begin{lemma} \label{lem:multitwist-lift-preserves}
	For each region $\Omega$, the lift $h_\Omega^{-1} \tilde h$ of $h$ preserves every lift of $\mathbf{t}$ adjacent to $\Omega$, with its sides.
\end{lemma}

\begin{proof}
	Lemma~\ref{lem:twist-adjacent-sides}, applied to each factor of $h_\Omega^{-1} \tilde h$, shows that $h_\Omega^{-1} \tilde h$ preserves the lifts adjacent to $\Omega$ with their sides.
\end{proof}

Let $\Psi$ be the automorphism of $\FF$ with $\tilde h(g x) = \Psi(g) \tilde h(x)$.

\begin{lemma} \label{lem:multitwist-lift-class}
	The outer class of $\Psi$ is $\Delta$.
\end{lemma}

\begin{proof}
	This follows from Lemma~\ref{lem:gul} applied to each factor of $\tilde h$.
\end{proof}

\begin{lemma} \label{lem:multitwist-lift-formula}
	We have
	\begin{equation} \label{eq:twist-formula}
		\Psi(g) = h_{gW} \, g \qquad \text{for all } g \in \FF .
	\end{equation}
\end{lemma}

\begin{proof}
	For $x \in W$ we have $\Psi(g) x = \tilde h(g x) = h_{gW} \, g x$ by Lemma~\ref{lem:multitwist-lift}, which gives \eqref{eq:twist-formula}.
\end{proof}

\begin{proposition} \label{prop:twist-growth}
	Let $\alpha$ be the homotopy class of an essential simple closed curve and $S$ the class of an essential sphere.
	Then there is a constant $K$, depending on $\alpha$, $S$ and $t_1, \dots, t_r$ but not on $m_1, \dots, m_r$, such that
	\[
		\left\lvert \, i(\Delta(\alpha), S) - \sum_{v=1}^{r} \lvert m_v \rvert \, i(\alpha, T_v) \, i(\gamma_v, S) \, \right\rvert \le K .
	\]
\end{proposition}

The proof occupies the rest of this subsection.
Let $s$ be a representative of the class $S$ of essential spheres, with dual tree $X_S$, and let $\tilde\alpha$ be a lift of the curve class $\alpha$, represented by $g$.

If $\tilde\alpha$ is a core of a lift of some $t_v$, then $g$ fixes the corresponding edge of $X$.
Then $\alpha$ is a power of $\gamma_v$, represented by a curve in $N(t_v)$ disjoint from all $t_u$, so $i(\alpha, T_u) = 0$ for all $u$.
Otherwise, penetration is defined through the sides in $C$, so by Lemma~\ref{lem:subtree-sides-penetrate} the lifts of $t_v$ at the edges of the minimal subtree of $\tilde\alpha$ on $X$ are the lifts of $t_v$ penetrated by $\tilde\alpha$.
By Proposition~\ref{prop:orbit-count} applied in the dual tree $X_v$ of $t_v$, they form $i(\alpha, T_v)$ orbits under $\langle g \rangle$.

If $g$ fixes a vertex $W$ of $X$, then $i(\alpha, T_v) = 0$ for all $v$: either $\tilde\alpha$ is a core of a lift of some $t_v$, or the minimal subtree of $\tilde\alpha$ on $X$ is the vertex $W$ and $\tilde\alpha$ penetrates no lift of any $t_v$ by Lemma~\ref{lem:subtree-sides-penetrate}.
Taking this $W$ as base region, \eqref{eq:twist-formula} gives $\Psi(g) = g$, so $\Delta(\alpha) = \alpha$ and the proposition holds with $K \ge i(\alpha, S)$.
So assume that $g$ is hyperbolic on $X$ with axis $A$; then $\tilde\alpha$ is a core of no lift of any $t_v$.

Let $q$ be the translation length of $g$ on $X$, let $W = x_0$ be a vertex of $A$, and let $x_0, x_1, \dots, x_q = g x_0$ be the vertices of $A$ from $x_0$ to $g x_0$, with edges $e_1, \dots, e_q$, where $e_j$ joins $x_{j-1}$ to $x_j$.
Extend the indexing to all of $A$ by $x_{j+q} = g x_j$ and $e_{j+q} = g e_j$.
Let $\tilde t_{(j)}$ be the lift at $e_j$, a lift of $t_{v(j)}$, and put $n_j := \lvert m_{v(j)} \rvert$ and $\ell(j) := i(\gamma_{v(j)}, S)$.
Each orbit of $\langle g \rangle$ on the edges of $A$ meets $\{ e_1, \dots, e_q \}$ exactly once, so exactly $i(\alpha, T_v)$ of the indices $j \in \{ 1, \dots, q \}$ have $v(j) = v$, and
\begin{equation} \label{eq:main-term}
	\sum_{j=1}^{q} n_j \, \ell(j) = \sum_{v=1}^{r} \lvert m_v \rvert \, i(\alpha, T_v) \, i(\gamma_v, S) .
\end{equation}

Let $\tilde h$ be the lift of $h$ that is the identity on $W$, and write $h_j := h_{x_j}$.
By Lemmas~\ref{lem:multitwist-lift}, \ref{lem:multitwist-lift-adjacent} and~\ref{lem:multitwist-lift-formula}, $h_0 = 1$ and $h_j = h_{j-1} u_j^{\,n_j}$ for a generator $u_j$ of the stabilizer of $\tilde t_{(j)}$, and $\Psi(g) = h_q \, g = u_1^{\,n_1} \dots u_q^{\,n_q} \, g$.
The proof of Lemma~\ref{lem:twist-curve} applies to $\tilde h$ with $\Psi$ in place of $\varphi_{W,n}$, whose outer class is $\Delta$ by Lemma~\ref{lem:multitwist-lift-class}, so $\tilde\alpha' := \tilde h(\tilde\alpha)$ is a lift of $\Delta(\alpha)$, represented by $g' := \Psi(g)$, and $i(\Delta(\alpha), S)$ is the translation length of $g'$ on $X_S$ by Proposition~\ref{prop:orbit-count-translation}.

\begin{lemma} \label{lem:growth-upper}
	There is $K_1$, independent of $m_1, \dots, m_r$, with $i(\Delta(\alpha), S) \le \sum_{j=1}^{q} n_j \ell(j) + K_1$.
\end{lemma}

\begin{proof}
	Let $x$ be a vertex of $X_S$.
	The translation length of $g'$ is at most $d(x, g' x)$, and by the triangle inequality and the invariance of $d$,
	\[
		d(x, g' x)
		= d(x, u_1^{\,n_1} \dots u_q^{\,n_q} g x)
		\le \sum_{j=1}^{q} d(x, u_j^{\,n_j} x) + d(x, g x) .
	\]
	The element $u_j$ represents the core of $\tilde t_{(j)}$, a lift of $\gamma_{v(j)}$, so its translation length on $X_S$ is $\ell(j)$ by Proposition~\ref{prop:orbit-count-translation}.
	By Theorem~\ref{thm:trees}\eqref{item:trees:dichotomy}, $d(x, u_j^{\,n_j} x) \le n_j \ell(j) + 2 d(x, L_j)$, where $L_j$ is the axis of $u_j$ if $\ell(j) > 0$ and the fixed point set of $u_j$ otherwise.
	In the second case this is the triangle inequality $d(x, u_j^{\,n_j} x) \le d(x, x') + d(u_j^{\,n_j} x', u_j^{\,n_j} x) = 2 d(x, x')$ for a point $x'$ of $L_j$ nearest to $x$.
	So the claim holds with $K_1 := 2 \sum_{j=1}^{q} d(x, L_j) + d(x, g x)$.
\end{proof}

For the lower bound, put $L := \max_v i(\gamma_v, S)$, $I_S := \max_v N(S, T_v)$ and $c_0 := 2 L^2$.
For $j \in \ZZ$, let $\tilde\gamma'_j$ be the core of $\tilde h(\tilde t_{(j)})$, let $A'_j$ be its minimal subtree on $X_S$, and let $\Sigma_j$ be the set of lifts of $s$ penetrated by both $\tilde\alpha'$ and $\tilde\gamma'_j$.
By Lemma~\ref{lem:subtree-sides-penetrate}, $\Sigma_j$ is the set of common edges of $A'_j$ and the minimal subtree $B'$ of $\tilde\alpha'$, a segment if $B'$ is an axis.

The set $\Sigma_j$ is finite.
Indeed, if $g'$ commuted with the element representing $\tilde\gamma'_j$, then $\tilde\alpha'$ would be a core of $\tilde h(\tilde t_{(j)})$, so $\Delta(\alpha) = \gamma_{v(j)}^p$ for some $p$.
But $\Delta$ fixes $\gamma_{v(j)}$, since each twist along $t_u$ is supported in a regular neighborhood disjoint from a core of $t_{v(j)}$ or carries that core to a parallel copy.
So $\alpha = \gamma_{v(j)}^p$ and $g$ fixes an edge of $X$, which contradicts the hyperbolicity of $g$.
Hence $\Sigma_j$ has at most $2 \, i(\Delta(\alpha), S) \, \ell(j)$ edges by Lemma~\ref{lem:orbits-bound}.

\begin{lemma} \label{lem:growth-overlap}
	$\Sigma_{j+q} = g' \Sigma_j$ for all $j$.
\end{lemma}

\begin{proof}
	Since $\tilde h(\tilde t_{(j+q)}) = \tilde h(g \tilde t_{(j)}) = g' \tilde h(\tilde t_{(j)})$ and $\tilde\alpha'$ is $g'$-invariant, $\Sigma_{j+q} = g' \Sigma_j$.
\end{proof}

\begin{lemma} \label{lem:growth-overlap-bound}
	$\Sigma_j \cap \Sigma_{j'}$ has at most $c_0$ edges for $j \ne j'$.
\end{lemma}

\begin{proof}
	For $j \ne j'$, the lifts $\tilde h(\tilde t_{(j)})$ and $\tilde h(\tilde t_{(j')})$ are distinct lifts of core-independent tori or of one torus, so the elements representing their cores do not commute by Lemmas~\ref{lem:joint-stabilizer} and~\ref{lem:end-stabilizer}.
	The edges of $A'_j \cap A'_{j'}$ are lifts of $s$ penetrated by both cores, which meet $\ell(j)$ and $\ell(j')$ orbits of these elements by Proposition~\ref{prop:orbit-count}; so there are at most $2 \ell(j) \ell(j') \le c_0$ such edges by Lemma~\ref{lem:orbits-bound}.
\end{proof}

\begin{lemma} \label{lem:growth-block}
	Put $K_2 := 4 L^2 + 2 L \, i(\alpha, S) + 8 L \, I_S + 2$.
	Then $\lvert \Sigma_j \rvert \ge n_j \ell(j) - K_2$ for every $j$.
\end{lemma}

\begin{proof}
	If $\ell(j) = 0$ or $n_j = 0$ there is nothing to prove, so let $\ell(j) > 0$ and $n_j \ge 1$.
	Write $w := u_j$, $n := n_j$, and $f := h_{j-1}^{-1} \tilde h$.
	We shall apply Corollary~\ref{lem:twist-sides-curve} and Lemma~\ref{lem:block-count} with $\tilde t_1 = \tilde t_{(j-1)}$, $\tilde t_2 = \tilde t_{(j)}$ and $\tilde t_3 = \tilde t_{(j+1)}$.
	The three lifts are pairwise disjoint, being distinct lifts of the disjoint tori $t_1, \dots, t_r$.
	Their stabilizers pairwise have no common nontrivial element by Lemmas~\ref{lem:joint-stabilizer} and~\ref{lem:end-stabilizer}, so they pairwise have no common end by Lemma~\ref{lem:common-end}.
	Let $P$ be the side of $\tilde t_{(j-1)}$, and $R$ the side of $\tilde t_{(j+1)}$, not containing the ends of $\tilde t_{(j)}$.
	The lift $\tilde\alpha$ penetrates $\tilde t_{(j-1)}$, $\tilde t_{(j)}$, $\tilde t_{(j+1)}$ in this order, so by Lemmas~\ref{lem:subtree-sides} and~\ref{lem:nested-sides} it has one end in $P$ and one end in $R$.

	By Lemma~\ref{lem:multitwist-lift-preserves}, $f$ preserves $\tilde t_{(j-1)}$ and $\tilde t_{(j)}$ with their sides, since they are adjacent to $x_{j-1}$, and $h_j^{-1} \tilde h$ preserves $\tilde t_{(j+1)}$, which is adjacent to $x_j$.
	Since $h_j = h_{j-1} w^n$, we have $f(\tilde t_{(j+1)}) = w^n \tilde t_{(j+1)}$.
	By Corollary~\ref{lem:twist-sides-curve}, applied to $f$ with $\tilde t_1 = \tilde t_{(j-1)}$, $\tilde t_2 = \tilde t_{(j)}$ and $\tilde t_3 = \tilde t_{(j+1)}$, the lift $\tilde c := f(\tilde\alpha) = h_{j-1}^{-1} \tilde\alpha'$ has one end in $P$ and one end in $w^n R$.

	We apply Lemma~\ref{lem:block-count} with $u = s$, with $\tilde t_1 = \tilde t_{(j-1)}$, $\tilde t_2 = \tilde t_{(j)}$ and $\tilde t_3 = \tilde t_{(j+1)}$, with the generator $w = u_j$ of the stabilizer of $\tilde t_{(j)}$, with the sides $P$ and $R$, and with the lifts $\tilde a = \tilde\alpha$ and $\tilde c = f(\tilde\alpha)$; the core-independence hypotheses are empty since $u$ is a sphere.
	Here $\tau = \ell(j) \ge 1$, $\theta \le L$ and $\iota \le I_S$.
	So at least $n \ell(j) - 4 \ell(j) L - 2 \ell(j) \, i(\alpha, S) - 8 \ell(j) I_S - 2 \ge n \ell(j) - K_2$ lifts of $s$ are penetrated by both $\tilde c$ and the core of $\tilde t_{(j)}$.
	Since $f$ preserves $\tilde t_{(j)}$, we have $\tilde h(\tilde t_{(j)}) = h_{j-1} \tilde t_{(j)}$, so $h_{j-1}$ carries the core of $\tilde t_{(j)}$ to the core $\tilde\gamma'_j$ of $\tilde h(\tilde t_{(j)})$ and $\tilde c$ to $\tilde\alpha'$.
	Applying $h_{j-1}$, these are therefore lifts penetrated by both $\tilde\alpha'$ and $\tilde\gamma'_j$.
\end{proof}

\begin{proof}[Proof of Proposition~\ref{prop:twist-growth}]
	By Lemma~\ref{lem:growth-upper} and \eqref{eq:main-term}, it remains to prove the lower bound $i(\Delta(\alpha), S) \ge \sum_{j=1}^{q} n_j \ell(j) - q (K_2 + c_0)$.
	Let $\mathcal{J}$ be the set of $j \in \ZZ$ with $n_j \ell(j) - K_2 > c_0$; it is invariant under $j \mapsto j + q$.
	If $\mathcal{J} = \varnothing$, then $\sum_{j=1}^{q} n_j \ell(j) \le q (K_2 + c_0)$ and the bound holds.
	Otherwise $\Sigma_j \ne \varnothing$ for $j \in \mathcal{J}$ by Lemma~\ref{lem:growth-block}, so $g'$ is hyperbolic on $X_S$ and $B'$ is its axis.

	By Lemmas~\ref{lem:growth-overlap}, \ref{lem:growth-overlap-bound} and~\ref{lem:growth-block}, the segments $\Sigma_j$, $j \in \mathcal{J}$, have more than $c_0$ edges and satisfy the hypotheses of Lemma~\ref{lem:segments} for $g'$ with $c = c_0$.
	Here $k$ is the number of orbits of $j \mapsto j + q$ on $\mathcal{J}$, and the segments are indexed accordingly.
	That lemma and Lemma~\ref{lem:growth-block} give
	\begin{align*}
		i(\Delta(\alpha), S)
		&\geq \sum_{j \in \mathcal{J} \cap \{ 1, \dots, q \}} \bigl( \lvert \Sigma_j \rvert - c_0 \bigr) \\
		&\geq \sum_{j \in \mathcal{J} \cap \{ 1, \dots, q \}} \bigl( n_j \ell(j) - K_2 - c_0 \bigr) 
        \geq \sum_{j=1}^{q} n_j \ell(j) - q (K_2 + c_0) ,
	\end{align*}
	using $n_j \ell(j) \le K_2 + c_0$ for $j \notin \mathcal{J}$.
	Taking $K := \max \bigl( K_1, \, q (K_2 + c_0), \, i(\alpha, S) \bigr)$ proves the proposition.
\end{proof}

\subsection{Free abelian subgroups} \label{subsec:free-abelian}

This subsection deduces from the growth estimate that Dehn twists along disjoint core-independent primitive tori generate a free abelian group of the expected rank, and then extends this to classes with pairwise disjoint representatives.

\begin{corollary} \label{cor:disjoint_twists_free_abelian}
	Let $t_1, \dots, t_r$ be pairwise disjoint primitive tori with oriented cores, whose classes $T_1, \dots, T_r$ are pairwise core-independent.
	Then $\delta_{\vec T_1}, \dots, \delta_{\vec T_r}$ generate a free abelian subgroup of $\Out(\FF)$ of rank $r$.
\end{corollary}

\begin{proof}
	The twists commute pairwise, since the homeomorphisms $\delta_{\vec t_v}$ have disjoint supports.
	Suppose that $\delta_{\vec T_1}^{m_1} \dots \delta_{\vec T_r}^{m_r} = 1$ with $m_{v_0} \ne 0$ for some $v_0$.
	Let $\alpha$ be the homotopy class of an essential simple closed curve and $S$ the class of an essential sphere.
	Raising the relation to the $q$-th power and applying Proposition~\ref{prop:twist-growth},
	\[
		i(\alpha, S)
		= i \bigl( \delta_{\vec T_1}^{q m_1} \dots \delta_{\vec T_r}^{q m_r}(\alpha), S \bigr)
		\ge q \sum_{v} \lvert m_v \rvert \, i(\alpha, T_v) \, i(\gamma_v, S) - K
	\]
	for every $q \ge 1$, where $K$ does not depend on $q$.
	Hence $\sum_{v} \lvert m_v \rvert \, i(\alpha, T_v) \, i(\gamma_v, S) = 0$, and in particular $i(\alpha, T_{v_0}) \, i(\gamma_{v_0}, S) = 0$ for every $\alpha$ and $S$.

	We now choose $\alpha$ and $S$ with $i(\alpha, T_{v_0}) \, i(\gamma_{v_0}, S) > 0$, which is a contradiction.
	Some $z \in \FF$ is hyperbolic on the dual tree of $T_{v_0}$.
	Otherwise $\FF$ would fix a vertex by Theorem~\ref{thm:trees}\eqref{item:trees:fg}, that is, stabilize a region $W$, so $C \cap W = C$ by Lemma~\ref{lem:bounded-limit}.
	But the side of a lift of $t_{v_0}$ adjacent to $W$ away from $W$ is nonempty and misses $C \cap W$.

	Represent $z$ by an embedded circle in $\MM$ and let $\alpha$ be its class; then $i(\alpha, T_{v_0}) > 0$ by Proposition~\ref{prop:orbit-count-translation}.
	Since $T_{v_0}$ is primitive, its core $\gamma_{v_0}$ is the class of a primitive curve, so $i(\gamma_{v_0}, S) = 1$ for some class $S$ of essential spheres by the definition in Section~\ref{subsec:primitive}.
	This contradiction shows that $m_v = 0$ for all $v$.
\end{proof}

The corollary needs the tori to be disjoint simultaneously.
The next proposition needs only pairwise disjoint representatives, at the price of two results from the theory of group actions on trees: the existence of a common refinement of pairwise compatible trees, due to Guirardel and Levitt, and Levitt's computation of the group of twists of a graph of groups.

\begin{proposition} \label{prop:pairwise-free-abelian}
	Let $T_1, \dots, T_r$ be pairwise core-independent classes of primitive tori such that any two of these classes have disjoint representatives, and let $\vec T_1, \dots, \vec T_r$ be oriented classes over $T_1, \dots, T_r$.
	Then $\delta_{\vec T_1}, \dots, \delta_{\vec T_r}$ generate a free abelian subgroup of $\Out(\FF)$ of rank $r$.
\end{proposition}

\begin{proof}
	The twists commute pairwise by Theorem~\ref{thm:commuting}.
	For each $a$, fix a representative $t_a \in T_a$ with dual tree $X_a := X_{t_a}$ and core $y_a$, the generator of the stabilizer of a lift representing the oriented core.

	The trees $X_a$ are pairwise compatible.
	Indeed, if $t \in T_a$ and $t' \in T_b$ are disjoint, the tree dual to the lifts of $t \cup t'$ collapses onto the dual trees of $t$ and of $t'$, which are equivariantly isomorphic to $X_a$ and $X_b$ by Section~\ref{subsec:classes}.

	Each $X_a$ is irreducible in the sense of \cite[Section~1.2]{MR3758992}, that is, $\FF$ fixes no point and no end of $X_a$ and preserves no line.
	Indeed, $X_a$ is minimal with edges by Section~\ref{subsec:dehn_prelim}, so $\FF$ fixes no point.
	If $\FF$ fixed an end, every vertex stabilizer would fix a ray, hence an edge, and would be cyclic.
	If $\FF$ preserved a line, that line would be $X_a$ by minimality; the kernel of the action on this line would fix every edge and be cyclic, so every vertex stabilizer would be virtually cyclic, hence cyclic.
	Both contradict Lemma~\ref{lem:no-cyclic-vertex}.
	Moreover $X_a$ has no redundant vertex in the sense of that section: at a vertex $v$ of valence two, the stabilizer of $v$ is not contained in the cyclic stabilizer of its edges by Lemma~\ref{lem:no-cyclic-vertex}, so some element of the stabilizer of $v$ exchanges the two edges and fixes only $v$.

	By Guirardel and Levitt \cite[Proposition~A.17]{MR3758992}, the trees $X_1, \dots, X_r$ have a common refinement $R$: an $\FF$-tree with equivariant collapse maps $R \to X_a$, such that no edge of $R$ is collapsed by all of these maps.

	We describe the edges of $R$.
	A collapse map is injective on the edges it does not collapse, and the preimage of an open edge of $X_a$ is a single open edge of $R$; by equivariance, an edge of $R$ not collapsed in $X_a$ has the same stabilizer as its image.
	So an edge of $R$ not collapsed in $X_a$ and in $X_b$ with $a \ne b$ would have a nontrivial stabilizer conjugate into $\langle y_a \rangle$ and into $\langle y_b \rangle$, contrary to Lemma~\ref{lem:joint-stabilizer}.
	Hence each edge of $R$ is collapsed in all $X_b$ but one, and since $\FF$ is transitive on the edges of $X_a$, the edges not collapsed in $X_a$ form a single orbit $E_a$.
	So the edges of $R$ form the $r$ orbits $E_1, \dots, E_r$, and $R \to X_a$ is injective on $E_a$ and collapses the other orbits.
	Let $\hat e_a \in E_a$ be the edge over the edge $e_a$ of $X_a$ fixed by $y_a$; its stabilizer is $\langle y_a \rangle$.

	No vertex stabilizer of $R$ is cyclic.
	Indeed, let $v$ be a vertex with cyclic stabilizer $\langle z \rangle$, and let $\hat e$ be an edge at $v$ in $E_a$; its stabilizer is conjugate to $\langle y_a \rangle$ and contained in $\langle z \rangle$.
	Since $y_a$ is not a proper power by Lemma~\ref{lem:primitive-basis} and Theorem~\ref{thm:free}\eqref{item:free:primitive}, this stabilizer is $\langle z \rangle$.
	The same holds for every edge at $v$, so by Lemma~\ref{lem:joint-stabilizer} all edges at $v$ lie in $E_a$.
	Collapsing the orbits other than $E_a$ does not identify $v$ with another vertex, since no edge at $v$ is collapsed; so the image of $v$ in $X_a$ has the cyclic stabilizer $\langle z \rangle$, contrary to Lemma~\ref{lem:no-cyclic-vertex}.

	Let $\mathcal{R} := R / \FF$ be the quotient graph of groups.
	Its vertex groups are free and not cyclic, hence have trivial center, and its edge groups are the infinite cyclic groups $\langle y_a \rangle$ up to conjugacy.
	It is minimal and is not a mapping torus, since a valence-one vertex whose group equals the edge group, or a single vertex whose group equals the edge group, would be cyclic.

	By Levitt \cite[Proposition~3.1]{MR2174093}, the group of twists of $\mathcal{R}$ is the quotient of $\prod_{\vec e} Z_{G_{o(\vec e)}}(G_{\vec e})$, over the oriented edges $\vec e$ of $\mathcal{R}$, by the relations coming from the centers of the vertex groups and from the edge groups.
	Here the group of twists is the subgroup of $\Out(\FF)$ generated by the twists around the edges of $\mathcal{R}$ by elements of the centralizers of the edge groups in the adjacent vertex groups.
	In our case the centralizer of $\langle y_a \rangle$ in a vertex group containing $\langle y_a \rangle$ is $\langle y_a \rangle$, by Theorem~\ref{thm:free}\eqref{item:free:commute} and since $y_a$ is not a proper power.
	The vertex relations are trivial, and the edge relation for $e_a$ identifies the twist by $y_a$ at one end of $e_a$ with the twist by $y_a^{-1}$ at the other end.
	So the group of twists of $\mathcal{R}$ is free abelian of rank $r$, freely generated by the twists $D_a$ around the edges $e_a$ by $y_a$.

	It remains to identify $D_a$ with $\delta_{\vec T_a}^{\pm 1}$.
	By definition \cite[Section~2]{MR2174093}, $D_a$ is induced by the automorphism of the path group of $\mathcal{R}$ that fixes the vertex groups and the letters $t_{e_b}$ for $b \ne a$ and sends $t_{e_a}$ to $y_a t_{e_a}$.
	The one-edge graph of groups $\mathcal{R} / (e_b, b \ne a)$ is $X_a / \FF$, since the collapse $R \to X_a$ is equivariant, and its vertex groups are generated by the vertex groups of $\mathcal{R}$ and the letters $t_{e_b}$, $b \ne a$, all fixed by $D_a$.
	So $D_a$ is the twist $D(X_a / \FF, (z_e))$ of Cohen and Lustig with twistor $y_a^{\pm 1}$ at $e_a$, up to inner automorphisms.
	As $o_{t_a}(\vec e_a) = y_a^{\pm 1}$, Lemma~\ref{lem:coincide-twist} gives $D_a = \delta_{\vec T_a}^{\pm 1}$.
	Hence $\delta_{\vec T_1}, \dots, \delta_{\vec T_r}$ generate a free abelian group of rank $r$.
\end{proof}

\section{Compatible collections and fellow travelers} \label{sec:fellow}

The goal of this section is Theorem~\ref{thm:main}: sufficiently large powers of the Dehn twists along a compatible collection of primitive tori generate a right-angled Artin group.
The proof follows Koberda \cite{MR3000498} and Seo \cite{MR4299673}, and uses the lifts of curves to $\widetilde\MM = \mathbb{S}^3 \setminus C$.
We record only the ends of a lift of a curve, defined in Section~\ref{subsec:curve-lifts}, and the two sides into which a lift of a torus cuts $C$, defined in Section~\ref{subsec:ends}; these do not depend on representatives.

\subsection{Compatible collections} \label{subsec:compatible}

This subsection fixes the collection $V$ of tori, its coincidence graph $\Gamma$ and the constant $D(\Gamma)$ that controls the ping-pong argument, and records the count of lifts behind that constant.
It then shows that adjacent members have crossing number zero and that the core of each has intersection number zero with the other, while non-adjacent members have positive crossing number.

Let $V = \{ T_1, \dots, T_k \}$ be pairwise core-independent classes of primitive tori, with cores $\gamma_1, \dots, \gamma_k$.
We say that $T_a$ and $T_b$ \emph{co-penetrate} if $i(\gamma_a, T_b) > 0$ and $i(\gamma_b, T_a) > 0$.
The \emph{coincidence graph} $\Gamma$ of $V$ has vertex set $V$, and distinct $T_a$ and $T_b$ are joined by an edge if and only if they have disjoint representatives.
The collection $V$ is \emph{compatible} if $T_a$ and $T_b$ co-penetrate whenever $a \ne b$ and $T_a$ and $T_b$ are not adjacent in $\Gamma$.
From now on we assume that $V$ is compatible.
We write $\operatorname{lk}(T_a)$ for the link of $T_a$ and $\operatorname{st}(T_a)$ for its star.

Put
\[
	M(\Gamma) := \max_{a, b} i(\gamma_a, T_b) , \qquad
	I(\Gamma) := \max_{a \ne b} N(T_a, T_b) , \qquad
	D(\Gamma) := 2 \bigl( M(\Gamma) + I(\Gamma) \bigr)^2 .
\]
We abbreviate these to $M$, $I$ and $D$; here $N$ is the crossing number of Section~\ref{sec:crossing}.

Fix representatives $t_a \in T_a$ and write $X_a$ for the dual tree of $t_a$.
This is harmless: by Lemma~\ref{lem:tree-sides}, the ends and sides of the lifts of a primitive torus, and whether a lift crosses a lift of another essential sphere or torus with no common end, depend only on the class of the torus.
The stabilizers of two distinct lifts of the tori $t_1, \dots, t_k$ have no common nontrivial element: by Lemma~\ref{lem:joint-stabilizer} for lifts of $t_a$ and $t_b$ with $a \ne b$, and by Lemma~\ref{lem:end-stabilizer} for two distinct lifts of $t_a$.
So their generators do not commute, and the lifts have no common end by Lemma~\ref{lem:common-end}.
In particular a core of a lift of $t_a$ is a core of no lift of $t_b$ for $b \ne a$, and every lift of $\gamma_a$ is a core of exactly one lift of $t_a$.

\begin{lemma} \label{lem:count-both}
	Let $a$, $b$, $c$ be indices with $b \ne a$ and $b \ne c$, and let $\tilde\gamma$ be a lift of $\gamma_a$.
	Let $\tilde t_c$ be a lift of $t_c$ of which $\tilde\gamma$ is not the core.
	Then at most $D$ lifts of $t_b$ are penetrated by $\tilde\gamma$ and cross $\tilde t_c$.
\end{lemma}

\begin{proof}
	Let $g$ represent $\tilde\gamma$; it generates the stabilizer of the lift $\tilde t_a$ of $t_a$ with core $\tilde\gamma$.
	By Proposition~\ref{prop:orbit-count}, the lifts of $t_b$ penetrated by $\tilde\gamma$ meet $i(\gamma_a, T_b) \le M$ orbits of $\langle g \rangle$.
	Let $g_c$ generate the stabilizer of $\tilde t_c$.
	By Lemma~\ref{lem:count-crossing}, the lifts of $t_b$ crossing $\tilde t_c$ meet $N(T_b, T_c) \le I$ orbits of $\langle g_c \rangle$.
	The elements $g$ and $g_c$ do not commute, since $\tilde t_a$ and $\tilde t_c$ are lifts of distinct members of $V$ or distinct lifts of one member.
	So Lemma~\ref{lem:orbits-bound} gives the bound $2 M I \le D$.
\end{proof}

The constant $D$ is larger than the bound $2 M I$ just obtained; it is chosen so as to absorb also the bound $6 M^2 + 8 M I$ that arises in the proof of Proposition~\ref{prop:co-penetrating}.

The next lemma relates the crossing number of two members of $V$ to the intersection numbers of their cores; together with compatibility, it shows below that non-adjacent members of $V$ have positive crossing number.

\begin{lemma} \label{lem:noncrossing-core}
	Let $a \ne b$ with $N(T_a, T_b) = 0$.
	Then no core of a lift of $t_a$ penetrates a lift of $t_b$.
\end{lemma}

\begin{proof}
	If the core of a lift $\tilde t_a$ of $t_a$ penetrated a lift $\tilde t_b$ of $t_b$, then $\tilde t_b$ and $\tilde t_a$ would cross by Lemma~\ref{lem:dual-sphere-crosses}, since $T_a$ and $T_b$ are core-independent; this contradicts $N(T_a, T_b) = 0$.
\end{proof}

\begin{lemma} \label{lem:noncrossing-core-intersection}
	Let $a \ne b$ with $N(T_a, T_b) = 0$.
	Then $i(\gamma_a, T_b) = 0$.
\end{lemma}

\begin{proof}
	Let $\tilde t_a$ be a lift of $t_a$.
	Its core is a lift of $\gamma_a$, penetrates no lift of $t_b$ by Lemma~\ref{lem:noncrossing-core}, and is a core of no lift of $t_b$; so $i(\gamma_a, T_b) = 0$ by Proposition~\ref{prop:orbit-count}.
\end{proof}

Equivalently, if $i(\gamma_a, T_b) > 0$, then $N(T_a, T_b) > 0$; this is the contrapositive of Lemma~\ref{lem:noncrossing-core-intersection}, and is also Corollary~\ref{lem:dual-sphere-crosses-positive}.

Let $T_a$ and $T_b$ be adjacent, with disjoint representatives $t$ and $t'$.
Then $N(T_a, T_b) = 0$ by Lemma~\ref{lem:count-disjoint}, and $i(\gamma_a, T_b) = 0$, since $\gamma_a$ is represented by a curve in a regular neighborhood of $t$ disjoint from $t'$.
For oriented classes $\vec T_a$ and $\vec T_b$ over $T_a$ and $T_b$, the twists $\delta_{\vec T_a}$ and $\delta_{\vec T_b}$ commute by Theorem~\ref{thm:commuting}.
If $a \ne b$ and $T_a$ and $T_b$ are not adjacent, then they co-penetrate since $V$ is compatible, so $N(T_a, T_b) > 0$ by Lemma~\ref{lem:noncrossing-core-intersection}.
Hence distinct $T_a$ and $T_b$ are adjacent if and only if $N(T_a, T_b) = 0$, and if and only if $i(\gamma_a, T_b) = 0$.

\subsection{Common lifts of a curve and a core} \label{subsec:common}

Throughout this subsection, let $T_i \in V$, let $\tilde t_i$ be a lift of $t_i$ with core $\tilde\gamma_i$, represented by $g_i$, and let $\tilde\gamma$ be a lift of the homotopy class $\gamma$ of an essential simple closed curve, represented by $g$.
Let $T_j \in \Gamma \setminus \operatorname{st}(T_i)$, so that $T_i$ and $T_j$ co-penetrate, since $V$ is compatible.
Suppose that $\tilde\gamma$ is a core of no lift of $t_j$, and that $\tilde\gamma \ne \tilde\gamma_i$.
Then $g$ and $g_i$ do not commute: otherwise $g$ would lie in the maximal cyclic subgroup $\langle g_i \rangle$, and $\tilde\gamma$ would be $\tilde\gamma_i$.

The lifts of $t_j$ penetrated by both $\tilde\gamma$ and $\tilde\gamma_i$ are the edges common to the minimal subtrees of $\tilde\gamma$ and of $\tilde\gamma_i$ on $X_j$, by Lemma~\ref{lem:subtree-sides-penetrate}.
By Proposition~\ref{prop:orbit-count} and Lemma~\ref{lem:orbits-bound}, there are at most $2\, i(\gamma, T_j)\, i(\gamma_i, T_j)$ such lifts.
The common edges of two geodesics in a tree form a geodesic, so these lifts are the edges of a finite segment of the axis of $g$, and we write them as
\[
	E_1, \dots, E_J ,
\]
in the order along the axis of $g$, oriented in the direction of translation by $g$.
Reversing the orientation of $\tilde\gamma_i$ if necessary, we assume that the axis of $g_i$ induces the same order.

For each $r$, let $P_r$ be the side of $E_r$ containing $g^{-\infty}$ and $Q_r$ the side containing $g^{+\infty}$.
By Lemma~\ref{lem:subtree-sides}, $g_i^{-\infty} \in P_r$ and $g_i^{+\infty} \in Q_r$.
By Lemmas~\ref{lem:nested-sides-inclusion} and~\ref{lem:nested-sides}, $P_r \subset P_{r'}$ and $Q_{r'} \subset Q_r$ for $r < r'$, and the ends of $E_r$ lie in $P_{r'}$ while the ends of $E_{r'}$ lie in $Q_r$.

Now let $T_\ell \in \operatorname{lk}(T_i)$, and let $\tilde t_\ell$ be a lift of $t_\ell$ penetrated by $\tilde\gamma$.
Let $Z$ and $Z^\ast$ be the sides of $\tilde t_\ell$, and let
\[
	\Lambda_{\tilde t_\ell} := \{\, r \mid \tilde t_\ell \text{ crosses } E_r \,\} .
\]
Since $T_i$ and $T_\ell$ are adjacent, $N(T_i, T_\ell) = 0$ by Section~\ref{subsec:compatible}, and the core $\tilde\gamma_i$ of $\tilde t_i$ penetrates no lift of $t_\ell$, by Lemma~\ref{lem:noncrossing-core}.
Its ends are not ends of $\tilde t_\ell$, so they lie on one side of $\tilde t_\ell$, say $g_i^{\pm\infty} \in Z$.

\begin{lemma} \label{lem:crossing-interval}
	$\Lambda_{\tilde t_\ell}$ is an interval of $\{ 1, \dots, J \}$.
\end{lemma}

\begin{proof}
	Let $a < c < a'$ with $a, a' \in \Lambda_{\tilde t_\ell}$, and suppose that $c \notin \Lambda_{\tilde t_\ell}$.
	By Lemma~\ref{lem:crossing-criterion}, there are $Z' \in \{ Z, Z^\ast \}$ and $P' \in \{ P_c, Q_c \}$ such that $Z' \cap P'$ contains no point other than ends of $\tilde t_\ell$ and of $E_c$.
	We treat $P' = P_c$; the case $P' = Q_c$ is the same with $a'$ in place of $a$ and the sets $Q$ in place of the sets $P$.
	Since $\tilde t_\ell$ crosses $E_a$, Lemma~\ref{lem:crossing-criterion} gives a point $\xi \in Z' \cap P_a$ that is not an end of $\tilde t_\ell$ or of $E_a$.
	Then $\xi \in P_c$, since $P_a \subset P_c$, and $\xi$ is not an end of $E_c$, since the ends of $E_c$ lie in $Q_a$, which is disjoint from $P_a$.
	This contradicts the choice of $Z'$ and $P'$.
\end{proof}

\begin{lemma} \label{lem:crossing-end}
	Suppose that $\Lambda_{\tilde t_\ell} \ne \varnothing$, and that $g^{-\infty} \in Z$.
	Then $J \in \Lambda_{\tilde t_\ell}$.
\end{lemma}

\begin{proof}
	Suppose that $J \notin \Lambda_{\tilde t_\ell}$, and let $b := \max \Lambda_{\tilde t_\ell} < J$.
	The ends $g_i^{-\infty} \in Z \cap P_J$ and $g_i^{+\infty} \in Z \cap Q_J$ are not ends of $\tilde t_\ell$ or of $E_J$, and neither is $g^{+\infty} \in Z^\ast \cap Q_J$.
	Since $\tilde t_\ell$ does not cross $E_J$, Lemma~\ref{lem:crossing-criterion} shows that $Z^\ast \cap P_J$ contains only ends.
	Since $\tilde t_\ell$ crosses $E_b$, there is $\xi \in Z^\ast \cap P_b$ that is not an end of $\tilde t_\ell$ or of $E_b$.
	Then $\xi \in P_J$, and $\xi$ is not an end of $E_J$, since the ends of $E_J$ lie in $Q_b$.
	This is a contradiction.
\end{proof}

\begin{lemma} \label{lem:crossing-bound}
	$\lvert \Lambda_{\tilde t_\ell} \rvert \le D$.
\end{lemma}

\begin{proof}
	Each $E_r$ with $r \in \Lambda_{\tilde t_\ell}$ is a lift of $t_j$ penetrated by $\tilde\gamma_i$ and crossing $\tilde t_\ell$.
	Lemma~\ref{lem:count-both} bounds the number of such lifts by $D$.
\end{proof}

\begin{lemma} \label{lem:middle-lifts}
	Let $D < r \le J - D$.
	Then no lift of $t_\ell$ penetrated by $\tilde\gamma$ crosses $E_r$, for any $T_\ell \in \operatorname{lk}(T_i)$.
\end{lemma}

\begin{proof}
	Let $\tilde t_\ell$ be such a lift with $r \in \Lambda_{\tilde t_\ell}$.
	By Lemmas~\ref{lem:crossing-interval}, \ref{lem:crossing-end} and~\ref{lem:crossing-bound}, $\Lambda_{\tilde t_\ell}$ is an interval with at most $D$ elements, and it contains $J$ if $g^{-\infty} \in Z$.
	If $g^{+\infty} \in Z$, we reverse the orientations of $\tilde\gamma$ and $\tilde\gamma_i$, which reverses the indexing, and conclude that $\Lambda_{\tilde t_\ell}$ contains $1$.
	So $\Lambda_{\tilde t_\ell} \subset \{ 1, \dots, D \} \cup \{ J - D + 1, \dots, J \}$, which does not contain $r$.
\end{proof}

\subsection{Fellow travelers} \label{subsec:fellow}

Let $q$ be a positive integer and $T_i \in V$.
The homotopy class $\gamma$ of an essential simple closed curve \emph{$q$-fellow-travels} with $T_i$ if $i(\gamma, T_i) > 0$ and there are a lift $\tilde\gamma$ of $\gamma$ and a lift $\tilde\gamma_i$ of $\gamma_i$ with the following property.
For every $T_j \in \Gamma \setminus \operatorname{st}(T_i)$, at least $q$ lifts of $t_j$ are penetrated by both $\tilde\gamma$ and $\tilde\gamma_i$.
We call $(\tilde\gamma, \tilde\gamma_i)$ a \emph{witness}.
Let $\mathcal{FT}_q(T_i)$ be the set of such classes $\gamma$, so that $\mathcal{FT}_1(T_i) \supseteq \mathcal{FT}_2(T_i) \supseteq \dots$.

Put
\[
	F(\Gamma) := 2 k D(\Gamma) + 3 .
\]

\begin{lemma} \label{lem:core_not_fellow}
	Let $T_0, T \in V$.
	Then $\gamma_0 \notin \mathcal{FT}_q(T)$ for every $q \ge 1$.
\end{lemma}

\begin{proof}
	Suppose that $\gamma_0 \in \mathcal{FT}_q(T)$.
	Then $i(\gamma_0, T) > 0$, so $T \ne T_0$, since $\gamma_0$ is represented by a curve in the regular neighborhood $N(t_0)$ disjoint from $t_0$.
	Moreover $T_0$ is not adjacent to $T$, since otherwise $i(\gamma_0, T) = 0$ by Section~\ref{subsec:compatible}; that is, $T_0 \in \Gamma \setminus \operatorname{st}(T)$.
	So some lift of $\gamma_0$ penetrates a lift of $t_0$.
	But every lift of $\gamma_0$ is a core of a lift of $t_0$, and penetrates no lift of $t_0$ by Lemma~\ref{lem:core-no-penetrate}.
\end{proof}

\subsection{Twisting along a torus in the link} \label{subsec:twist-link}

The goal of this subsection is Proposition~\ref{prop:link-twist}: twisting along tori in the link of $T_i$ does not destroy fellow-traveling with $T_i$.
For each $a$, fix one of the two oriented classes $\vec T_a$ over $T_a$ and write $\delta_a := \delta_{\vec T_a}$.
Orient the core of the fixed representative $t_a$ so that $\vec t_a \in \vec T_a$, which is possible by Section~\ref{subsec:classes}; when a lift of $\delta_a^{\,n}$ is taken, $\delta_a^{\,n}$ also denotes the homeomorphism $\delta_{\vec t_a}^{\,n}$, which represents it by Section~\ref{subsec:twist-homeo}.
% We use the lifts of powers of Dehn twists of Section~\ref{subsec:twist-lift}, in particular Lemmas~\ref{lem:twist-adjacent}, \ref{lem:twist-adjacent-sides}, \ref{lem:twist-adjacent-image} and~\ref{lem:twist-curve}.

Throughout this subsection, let $T_i \in V$, $T_\ell \in \operatorname{lk}(T_i)$, $n \in \ZZ$ and $q \ge 2D + 2$, and let $\gamma \in \mathcal{FT}_q(T_i)$ with witness $(\tilde\gamma, \tilde\gamma_i)$, where $g$ represents $\tilde\gamma$ and $g_i$ represents $\tilde\gamma_i$.
Since $T_i$ and $T_\ell$ are adjacent, $i(\gamma_i, T_\ell) = 0$ and $\tilde\gamma_i$ is a core of no lift of $t_\ell$, by Section~\ref{subsec:compatible}.
So $g_i$ fixes no edge of $X_\ell$ and has translation length $0$ on $X_\ell$ by Proposition~\ref{prop:orbit-count-translation}, and by Theorem~\ref{thm:trees}\eqref{item:trees:dichotomy} it fixes exactly one vertex $W$ of $X_\ell$, the minimal subtree of $\tilde\gamma_i$.
The vertex $W$ is a region of $t_\ell$, stabilized by $g_i$, so $g_i^{\pm\infty} \in C \cap W$ by Lemma~\ref{lem:end-in-side}.
Let $\tilde\delta$ be the lift of $\delta_\ell^{\,n}$ that is the identity on $W$.
It fixes the ends of $\tilde\gamma_i$ by the remark after Lemma~\ref{lem:twist-fixes-positive}, and $\tilde\delta(\tilde\gamma)$ is a lift of $\delta_\ell^{\,n}(\gamma)$ by Lemma~\ref{lem:twist-curve}.

Since $i(\gamma, T_i) > 0$, the class $\gamma$ is not a power of $\gamma_i$, so $\tilde\gamma \ne \tilde\gamma_i$.
If $\Gamma \setminus \operatorname{st}(T_i) = \varnothing$, then $\mathcal{FT}_q(T_i)$ is the set of classes $\gamma$ with $i(\gamma, T_i) > 0$ for every $q$.
Otherwise fix $T_j \in \Gamma \setminus \operatorname{st}(T_i)$; the statements below that involve $T_j$ are made under this assumption.
Since $\tilde\gamma$ penetrates at least $q \ge 1$ lifts of $t_j$, it is a core of no lift of $t_j$ by Lemma~\ref{lem:core-no-penetrate}, so Section~\ref{subsec:common} applies.
Let $E_1, \dots, E_J$ be the lifts of $t_j$ penetrated by both $\tilde\gamma$ and $\tilde\gamma_i$, in the order along the axis of $g$ fixed in Section~\ref{subsec:common}, where $\tilde\gamma_i$ is oriented so that the axis of $g_i$ induces the same order.
For each $r$, let $P_r$ be the side of $E_r$ containing $g^{-\infty}$ and $Q_r$ the side containing $g^{+\infty}$.
Then $J \ge q$, and $g^{-\infty}, g_i^{-\infty} \in P_r$ and $g^{+\infty}, g_i^{+\infty} \in Q_r$ for every $r$.

\begin{lemma} \label{lem:link-subtree}
	$\tilde\gamma$ is a core of no lift of $t_\ell$.
\end{lemma}

\begin{proof}
	Suppose that $\tilde\gamma$ were a core of a lift $\tilde t_\ell$ of $t_\ell$, so that the ends of $\tilde\gamma$ are those of the core of $\tilde t_\ell$.
	By Lemma~\ref{lem:noncrossing-core}, with the roles of $T_i$ and $T_\ell$ exchanged, which applies since $N(T_\ell, T_i) = 0$ by Section~\ref{subsec:compatible}, the core of $\tilde t_\ell$ penetrates no lift of $t_i$, and neither does $\tilde\gamma$, since penetration depends only on the ends.
	Moreover $\tilde\gamma$ is a core of no lift of $t_i$, since the lifts of $t_i$ and of $t_\ell$ have no common end.
	So $i(\gamma, T_i) = 0$ by Proposition~\ref{prop:orbit-count}, a contradiction.
\end{proof}

By Lemma~\ref{lem:link-subtree}, $\tilde\gamma$ has a minimal subtree $A$ on $X_\ell$.

\begin{lemma} \label{lem:link-subtree-vertex}
	If $\Gamma \setminus \operatorname{st}(T_i) \ne \varnothing$, then the minimal subtree $A$ of $\tilde\gamma$ on $X_\ell$ contains $W$.
\end{lemma}

\begin{proof}
	Let $\Gamma \setminus \operatorname{st}(T_i) \ne \varnothing$, fix $T_j \in \Gamma \setminus \operatorname{st}(T_i)$, and let $E_1, \dots, E_J$ be the lifts of $t_j$ penetrated by both $\tilde\gamma$ and $\tilde\gamma_i$, where $J \ge q$.
	For each $r$, let $P_r$ be the side of $E_r$ containing $g^{-\infty}$ and $g_i^{-\infty}$, and $Q_r$ the side of $E_r$ containing $g^{+\infty}$ and $g_i^{+\infty}$.
	Suppose that $W \notin A$.
	Let $\tilde t_\ell$ be the first edge of the geodesic from $W$ to $A$.
	Let $Z := C(Y_\kappa)$ and $Z^\ast := C(Y_{\kappa^\ast})$ be the sides of $\tilde t_\ell$, where $\kappa$ is the component of $X_\ell \setminus \tilde t_\ell$ containing $A$ and $\kappa^\ast$ the component containing $W$.
	By Lemma~\ref{lem:subtree-sides}, $g^{\pm\infty} \in Z$.
	The minimal subtree of $\tilde\gamma_i$ is the vertex $W \in \kappa^\ast$, so $g_i^{\pm\infty} \in Z^\ast$ by Lemma~\ref{lem:subtree-sides}.

	For every $r$, each of the four sets $Z \cap P_r$, $Z \cap Q_r$, $Z^\ast \cap P_r$, $Z^\ast \cap Q_r$ contains an end of $\tilde\gamma$ or of $\tilde\gamma_i$, which is not an end of $\tilde t_\ell$ or of $E_r$.
	By Lemma~\ref{lem:crossing-criterion}, $\tilde t_\ell$ crosses every $E_r$.
	But the $E_r$ are penetrated by $\tilde\gamma_i$, so at most $D$ of the lifts $E_r$ cross $\tilde t_\ell$ by Lemma~\ref{lem:count-both}, while $J \ge q > D$.
\end{proof}

Fix $r$ with $D < r \le J - D$, and write $E := E_r$, $P := P_r$, $Q := Q_r$.
By Lemma~\ref{lem:middle-lifts}, no lift of $t_\ell$ penetrated by $\tilde\gamma$ crosses $E$.

Suppose that $A$ is the axis of $g$, let $\tilde t_\ell^+$ be the edge of $A$ at $W$ on the positive ray, and let $Z^+$ be the side of $\tilde t_\ell^+$ not containing $g_i^{+\infty}$.
Let $\kappa^+$ be the component of $X_\ell \setminus \tilde t_\ell^+$ not containing $W$ and $\kappa$ the other one.
The minimal subtree of $\tilde\gamma_i$ is the vertex $W \in \kappa$, so $g_i^{+\infty} \in C(Y_\kappa)$ by Lemma~\ref{lem:subtree-sides}; hence $C(Y_\kappa)$ is the side containing $g_i^{+\infty}$, and $Z^+ = C(Y_{\kappa^+})$.

\begin{lemma} \label{lem:link-sides}
	$g^{+\infty} \in Z^+$.
\end{lemma}

\begin{proof}
	By Lemma~\ref{lem:subtree-sides}, $g^{+\infty} \in C(Y_{\kappa^+}) = Z^+$.
\end{proof}

\begin{lemma} \label{lem:link-sides-contained}
	Every point of $Z^+$ lies in $Q$ or is an end of $E$.
\end{lemma}

\begin{proof}
	Since $\tilde t_\ell^+$ is penetrated by $\tilde\gamma$, it does not cross $E$.
	By Lemmas~\ref{lem:link-sides} and~\ref{lem:subtree-sides}, $g^{+\infty} \in Z^+$ and $g^{-\infty} \in C(Y_\kappa)$.
	The point $g^{+\infty} \in Z^+ \cap Q$, the point $g^{-\infty} \in C(Y_\kappa) \cap P$ and the point $g_i^{+\infty} \in C(Y_\kappa) \cap Q$ are not ends of $\tilde t_\ell^+$ or of $E$.
	So three of the four sets of Lemma~\ref{lem:crossing-criterion} for the pair $(\tilde t_\ell^+, E)$ contain points that are not ends, and the fourth, $Z^+ \cap P$, contains only ends.
	A point of $Z^+$ is not an end of $\tilde t_\ell^+$, so it lies in $Q$ or is an end of $E$.
\end{proof}

\begin{lemma} \label{lem:link-penetrates}
	The end $\tilde\delta(g^{+\infty})$ lies in $Q$ or is an end of $E$.
\end{lemma}

\begin{proof}
	By Lemma~\ref{lem:link-subtree-vertex}, $W \in A$.
	If $A = \{ W \}$, then $g$ stabilizes $W$, so $g^{\pm\infty} \in C \cap W$ by Lemma~\ref{lem:end-in-side}, and $\tilde\delta$ fixes them by the remark after Lemma~\ref{lem:twist-fixes-positive}; so $\tilde\delta(g^{+\infty}) = g^{+\infty}$ lies in $Q$.
	Otherwise $A$ is the axis of $g$ through $W$.
	Let $\tilde t_\ell^+$ be the edge of $A$ at $W$ on the positive ray, and $Z^+$ the side of $\tilde t_\ell^+$ not containing $g_i^{+\infty}$.
	Then $g^{+\infty} \in Z^+$ by Lemma~\ref{lem:link-sides}, and $\tilde\delta$ preserves $Z^+$ by Lemma~\ref{lem:twist-adjacent-sides}, so $\tilde\delta(g^{+\infty}) \in Z^+$ lies in $Q$ or is an end of $E$ by Lemma~\ref{lem:link-sides-contained}.
\end{proof}

\begin{lemma} \label{lem:link-penetrates-negative}
	The end $\tilde\delta(g^{-\infty})$ lies in $P$ or is an end of $E$.
\end{lemma}

\begin{proof}
	Reversing the orientations of $\tilde\gamma$ and $\tilde\gamma_i$ exchanges $P$ with $Q$ and the positive ray with the negative one, so this is Lemma~\ref{lem:link-penetrates} for the reversed orientations.
\end{proof}

\begin{lemma} \label{lem:twist-link}
	Let $T_i \in V$, $T_\ell \in \operatorname{lk}(T_i)$, $n \in \ZZ$ and $q \ge 2D + 2$, and let $\gamma \in \mathcal{FT}_q(T_i)$ with witness $(\tilde\gamma, \tilde\gamma_i)$.
	Let $\tilde\delta$ be the lift of $\delta_\ell^{\,n}$ that is the identity on the region $W$ of $t_\ell$ stabilized by the element representing $\tilde\gamma_i$.
	Then $\delta_\ell^{\,n}(\gamma) \in \mathcal{FT}_{q - 2D}(T_i)$, with witness $(\tilde\delta(\tilde\gamma), \tilde\gamma_i)$.
\end{lemma}

\begin{proof}
	First, since $T_\ell$ and $T_i$ are adjacent, there are disjoint tori $t \in T_i$ and $t' \in T_\ell$; choose a regular neighborhood $N(t')$ disjoint from $t$, and orient the core of $t'$ so that $\vec t'$ lies in $\vec T_\ell$, which is possible by Section~\ref{subsec:classes}.
	Then $\delta_{\vec t'}^{\,-n}$ represents $\delta_\ell^{\,-n}$ by Section~\ref{subsec:twist-homeo} and fixes $t$, so $\delta_\ell^{\,-n}(T_i)$ is the class of $t$, which is $T_i$; hence $i(\delta_\ell^{\,n}(\gamma), T_i) = i(\gamma, \delta_\ell^{\,-n}(T_i)) = i(\gamma, T_i) > 0$.
	This proves the claim if $\Gamma \setminus \operatorname{st}(T_i) = \varnothing$.

	Otherwise, since $T_j \in \Gamma \setminus \operatorname{st}(T_i)$ was arbitrary, it suffices to show that $\tilde\delta(\tilde\gamma)$ penetrates $E_r$ for each $r$ with $D < r \le J - D$, where $E_1, \dots, E_J$ are the lifts of $t_j$ penetrated by both $\tilde\gamma$ and $\tilde\gamma_i$, as there are at least $q - 2D \ge 2$ such indices $r$.
	If $\tilde\delta(\tilde\gamma)$ is not a core of $E_r$, then its ends are not ends of $E_r$, so it penetrates $E_r$ by Lemmas~\ref{lem:link-penetrates} and~\ref{lem:link-penetrates-negative}.
	A lift of a curve is a core of at most one lift of $t_j$.
	If $\tilde\delta(\tilde\gamma)$ were a core of $E_r$, it would penetrate $E_{r'}$ for another such $r'$, contradicting Lemma~\ref{lem:core-no-penetrate}.
\end{proof}

\begin{proposition} \label{prop:link-twist}
	Let $T_i \in V$, let $T_{j_1}, \dots, T_{j_m} \in \operatorname{lk}(T_i)$ be distinct, and let $n_1, \dots, n_m \in \ZZ$.
	Then
	\[
		\delta_{j_m}^{\,n_m} \dots \delta_{j_1}^{\,n_1} \bigl( \mathcal{FT}_{F(\Gamma)}(T_i) \bigr) \subseteq \mathcal{FT}_3(T_i) .
	\]
\end{proposition}

\begin{proof}
	Since the indices are distinct, $m \le k$.
	So $F(\Gamma) - 2D(s-1) \ge 2D + 3$ for $s \le m$, and Lemma~\ref{lem:twist-link} applied $m$ times gives $\mathcal{FT}_{F(\Gamma)}(T_i) \to \mathcal{FT}_{F(\Gamma) - 2Dm}(T_i) \subseteq \mathcal{FT}_3(T_i)$.
\end{proof}

\subsection{Twisting along a co-penetrating torus} \label{subsec:co-penetrating}

The goal of this subsection is Proposition~\ref{prop:co-penetrating}: a large power of the twist along $T_j$ carries the curves that fellow-travel a co-penetrating torus $T_i$, and the core of $T_i$, to curves that $F(\Gamma)$-fellow-travel $T_j$.

\begin{proposition} \label{prop:co-penetrating}
	Let $T_i, T_j \in V$ co-penetrate, and let $n \in \ZZ$ with $\lvert n \rvert \ge (2k + 7) D(\Gamma) + 4$.
	Then
	\[
		\delta_j^{\,n} \bigl( \mathcal{FT}_3(T_i) \cup \{ \gamma_i \} \bigr) \subseteq \mathcal{FT}_{F(\Gamma)}(T_j) .
	\]
\end{proposition}

\begin{proof}
	Let $\gamma \in \mathcal{FT}_3(T_i) \cup \{ \gamma_i \}$.
	Since $i(\gamma_i, T_j) > 0$, the classes $T_i$ and $T_j$ are distinct and not adjacent by Section~\ref{subsec:compatible}, and $N(T_i, T_j) > 0$ by Lemma~\ref{lem:noncrossing-core-intersection}; so $T_j \in \Gamma \setminus \operatorname{st}(T_i)$.

	We first note that $i(\delta_j^{\,n}(\gamma), T_j) > 0$.
	Since the twist along $t_j$ preserves $t_j$, $i(\delta_j^{\,n}(\gamma), T_j) = i(\gamma, \delta_j^{-n}(T_j)) = i(\gamma, T_j)$.
	If $\gamma = \gamma_i$, then $i(\gamma, T_j) > 0$ is co-penetration.
	Otherwise some lift of $\gamma$ penetrates a lift of $t_j$, so $i(\gamma, T_j) > 0$ by Proposition~\ref{prop:orbit-count}.

	Next we choose three consecutive lifts of $t_j$ along a lift of $\gamma$.
	If $\gamma \in \mathcal{FT}_3(T_i)$, let $(\tilde c, \tilde\gamma_i)$ be a witness.
	Then $\tilde c$ is a core of no lift of $t_j$ by Lemma~\ref{lem:core-no-penetrate}, and $\tilde c \ne \tilde\gamma_i$ since $i(\gamma, T_i) > 0$.
	So Section~\ref{subsec:common} applies, and there are three consecutive lifts $\tilde t_j, \tilde t_j', \tilde t_j''$ of $t_j$ penetrated by both $\tilde c$ and $\tilde\gamma_i$, consecutive edges of the axis of the element $g$ representing $\tilde c$.
	If $\gamma = \gamma_i$, let $\tilde c = \tilde\gamma_i$ be a lift of $\gamma_i$, represented by $g$.
	Its minimal subtree on $X_j$ is an axis by Proposition~\ref{prop:orbit-count-translation}, and we let $\tilde t_j, \tilde t_j', \tilde t_j''$ be the lifts at three consecutive edges of this axis.

	Let $\tilde\gamma'$ be the core of $\tilde t_j'$.
	Let $W$ be the region of $t_j$ between $\tilde t_j$ and $\tilde t_j'$, and let $\tilde\delta$ be the lift of $\delta_j^{\,n}$ that is the identity on $W$.
	By Lemmas~\ref{lem:twist-adjacent-sides} and~\ref{lem:twist-adjacent-image}, $\tilde\delta$ preserves $\tilde t_j$ and $\tilde t_j'$ with their sides, and $\tilde\delta(\tilde t_j'') = w^{\lvert n \rvert} \tilde t_j''$ for a suitable generator $w$ of the stabilizer of $\tilde t_j'$.
	By Lemma~\ref{lem:twist-curve}, $\tilde\delta(\tilde c)$ is a lift of $\delta_j^{\,n}(\gamma)$.
	We show that $(\tilde\delta(\tilde c), \tilde\gamma')$ is a witness for $\delta_j^{\,n}(\gamma) \in \mathcal{FT}_{F(\Gamma)}(T_j)$.

	We locate the ends of $\tilde\delta(\tilde c)$.
	The lifts $\tilde t_j, \tilde t_j', \tilde t_j''$ are distinct lifts of $t_j$, so they are pairwise disjoint and their stabilizers pairwise have no common nontrivial element, by Lemma~\ref{lem:end-stabilizer}.
	Let $P$ be the side of $\tilde t_j$, and $R$ the side of $\tilde t_j''$, not containing the ends of $\tilde t_j'$.
	Orient $\tilde c$ so that its axis passes through $\tilde t_j, \tilde t_j', \tilde t_j''$ in this order.
	By Lemmas~\ref{lem:subtree-sides} and~\ref{lem:nested-sides}, $g^{-\infty} \in P$ and $g^{+\infty} \in R$, and likewise $\tilde\gamma_i$ has one end in $P$ and one in $R$.
	By Corollary~\ref{lem:twist-sides-curve}, applied to the lifts $\tilde t_j$, $\tilde t_j'$, $\tilde t_j''$ and $f = \tilde\delta$, the lift $\tilde\delta(\tilde c)$ has one end in $P$ and one end in $w^{\lvert n \rvert} R$.

	Fix $T_\ell \in \Gamma \setminus \operatorname{st}(T_j)$; since $V$ is compatible, $T_j$ and $T_\ell$ co-penetrate, and $T_\ell = T_i$ is allowed.
	We apply Lemma~\ref{lem:block-count} with $u = t_\ell$, with $\tilde t_1 = \tilde t_j$, $\tilde t_2 = \tilde t_j'$ and $\tilde t_3 = \tilde t_j''$, with the generator $w$ of the stabilizer of $\tilde t_j'$ and the core $\tilde\gamma'$ of $\tilde t_j'$, with the sides $P$ and $R$, with $\lvert n \rvert$ in the role of $n$, and with $\tilde\gamma_i$ in the role of $\tilde a$ and $\tilde\delta(\tilde c)$ in the role of $\tilde c$.
	The hypotheses of Lemma~\ref{lem:block-count} hold: $T_j$ and $T_\ell$ are core-independent, the lifts $\tilde t_j, \tilde t_j', \tilde t_j''$ are pairwise disjoint with stabilizers pairwise having no common nontrivial element, the lift $\tilde\gamma_i$ has one end in $P$ and one end in $R$, the lift $\tilde\delta(\tilde c)$ has one end in $P$ and one end in $w^{\lvert n \rvert} R$, and $\tau = i(\gamma_j, T_\ell) \ge 1$ by co-penetration.

	It remains to compare the resulting bound with $F(\Gamma)$.
	Here $\tau \le M$, $\theta \le M$, $\iota \le I$ and $i(\gamma_i, T_\ell) \le M$, so
	\[
		4 \tau \theta + 2 \tau \, i(\gamma_i, T_\ell) + 8 \tau \iota + 2
		\le 6 M^2 + 8 M I + 2
		\le 3 D + 2 ,
	\]
	using $D = 2 (M + I)^2$.
	So at least $\lvert n \rvert \tau - 3D - 2 \ge \lvert n \rvert - 3D - 2$ lifts of $t_\ell$ are penetrated by both $\tilde\delta(\tilde c)$ and $\tilde\gamma'$.
	Since $T_i$ and $T_j$ co-penetrate, $M \ge 1$, so $D \ge 2$.
	Finally $\lvert n \rvert - 3D - 2 \ge (2k + 4) D + 2 \ge 2kD + 3 = F(\Gamma)$, since $\lvert n \rvert \ge (2k + 7) D + 4$ and $D \ge 2$.
\end{proof}

\subsection{The right-angled Artin group} \label{subsec:raag}

Let $A(\Gamma)$ be the right-angled Artin group of $\Gamma$: the group generated by $T_1, \dots, T_k$ subject to the relations $[T_a, T_b] = 1$ for the edges $\{ T_a, T_b \}$ of $\Gamma$.

\begin{theorem} \label{thm:main}
	Let $V = \{ T_1, \dots, T_k \}$ be a compatible collection with coincidence graph $\Gamma$.
	Let $n_1, \dots, n_k$ be integers with $\lvert n_a \rvert \ge (2k + 7) D(\Gamma) + 4$ for all $a$.
	Then the homomorphism
	\[
		\Phi \colon A(\Gamma) \to \Out(\FF) , \qquad T_a \mapsto \delta_a^{\,n_a}
	\]
	is injective.
\end{theorem}

We note that the homomorphism $\Phi$ is well-defined by Theorem~\ref{thm:commuting}, since adjacent members of $V$ have disjoint representatives.
% From Sections~\ref{sec:crossing}--\ref{sec:growth}, the proof uses chiefly Theorem~\ref{thm:commuting}, Lemma~\ref{lem:count-crossing}, Lemma~\ref{lem:block-count} and Proposition~\ref{prop:pairwise-free-abelian}, together with Lemmas~\ref{lem:count-disjoint} and~\ref{lem:dual-sphere-crosses}, Proposition~\ref{prop:orbit-count} and the lifts of Dehn twists of Section~\ref{subsec:twist-lift}.
% Proposition~\ref{prop:pairwise-free-abelian} is used only for the members of $V$ adjacent to all other members; when these have simultaneously disjoint representatives, Corollary~\ref{cor:disjoint_twists_free_abelian} replaces it and the proof is self-contained.

We recall the normal form of Koberda \cite[Section~6.8]{MR3000498}.
For a nontrivial $x \in A(\Gamma)$, its \emph{support} $\operatorname{supp}(x)$ is the set of vertices occurring in a reduced word representing $x$; it does not depend on the reduced word.
An element is \emph{clique-supported} (also known as \emph{central}) if its support spans a complete subgraph of $\Gamma$.
A product $x = x_r \dots x_1$ of clique-supported elements, with $x_r$ the leftmost factor, is a \emph{left greedy normal form} if for each $s < r$ and each $v \in \operatorname{supp}(x_s)$ there is $u \in \operatorname{supp}(x_{s+1})$ with $[u, v] \ne 1$, that is, $u \ne v$ and $u$ not adjacent to $v$.
Every nontrivial element of $A(\Gamma)$ has a left greedy normal form, which can be obtained by repeatedly splitting off a maximal central factor on the left.

In the proof below, the factor $x_1$ acts first on the curve $\gamma_0$, and $x_s$ acts after $x_{s-1}$.
A graph is \emph{centerless} if no vertex is adjacent to all other vertices; then $A(\Gamma)$ has trivial center.

\begin{lemma} \label{lem:step}
	Let $x$ be a clique-supported element of $A(\Gamma)$ and $T_v \in \operatorname{supp}(x)$.
	Let $T_u \in V$ be distinct from $T_v$ and not adjacent to $T_v$, and let $\gamma \in \mathcal{FT}_3(T_u) \cup \{ \gamma_u \}$.
	Then $\Phi(x)(\gamma) \in \mathcal{FT}_3(T_v)$.
\end{lemma}

\begin{proof}
	Since $x$ is clique-supported, $x = y \, T_v^{\,m}$ with $m \ne 0$ and $\operatorname{supp}(y) \subseteq \operatorname{lk}(T_v)$.
	Since $V$ is compatible, $T_u$ and $T_v$ co-penetrate.
	As $\lvert m n_v \rvert \ge (2k + 7) D(\Gamma) + 4$, Proposition~\ref{prop:co-penetrating} gives $\Phi(T_v^{\,m})(\gamma) = \delta_v^{\,m n_v}(\gamma) \in \mathcal{FT}_{F(\Gamma)}(T_v)$.
	Since $\operatorname{supp}(y) \subseteq \operatorname{supp}(x)$ spans a complete subgraph, $\Phi(y)$ is a product of one power of the twist along each member of $\operatorname{supp}(y) \subseteq \operatorname{lk}(T_v)$.
	So Proposition~\ref{prop:link-twist} gives $\Phi(x)(\gamma) = \Phi(y) \bigl( \Phi(T_v^{\,m})(\gamma) \bigr) \in \mathcal{FT}_3(T_v)$.
\end{proof}

\begin{lemma} \label{lem:main-centerless}
	Let $V = \{ T_1, \dots, T_k \}$ be a compatible collection with coincidence graph $\Gamma$, let $n_1, \dots, n_k$ be integers with $\lvert n_a \rvert \ge (2k + 7) D(\Gamma) + 4$ for all $a$, and let $\Phi \colon A(\Gamma) \to \Out(\FF)$ be the homomorphism with $\Phi(T_a) = \delta_a^{\,n_a}$ of Theorem~\ref{thm:main}.
	Suppose that $\Gamma$ is centerless.
	Then $\Phi$ is injective.
\end{lemma}

\begin{proof}
	Let $x \in A(\Gamma)$ be nontrivial, with left greedy normal form $x = x_r \dots x_1$.
	Some vertex $T_0 \notin \operatorname{supp}(x_1)$ is not adjacent to some $T_{i_1} \in \operatorname{supp}(x_1)$: otherwise a vertex of the complete subgraph spanned by $\operatorname{supp}(x_1)$ would be adjacent to all other vertices.
	By Lemma~\ref{lem:step}, $\Phi(x_1)(\gamma_0) \in \mathcal{FT}_3(T_{i_1})$.

	Suppose inductively that $\Phi(x_{s-1} \dots x_1)(\gamma_0) \in \mathcal{FT}_3(T_{i_{s-1}})$ with $T_{i_{s-1}} \in \operatorname{supp}(x_{s-1})$, for some $2 \le s \le r$.
	By the left greedy condition, some $T_{i_s} \in \operatorname{supp}(x_s)$ is different from and not adjacent to $T_{i_{s-1}}$.
	By Lemma~\ref{lem:step}, $\Phi(x_s \dots x_1)(\gamma_0) \in \mathcal{FT}_3(T_{i_s})$.

	By induction, $\Phi(x)(\gamma_0) \in \mathcal{FT}_3(T_{i_r})$.
	By Lemma~\ref{lem:core_not_fellow}, $\gamma_0 \notin \mathcal{FT}_3(T_{i_r})$.
	So $\Phi(x)(\gamma_0) \ne \gamma_0$ and $\Phi(x) \ne 1$.
\end{proof}

\begin{proof}[Proof of Theorem~\ref{thm:main}]
	Let $V_0$ be the set of members of $V$ adjacent to all other members.
	Let $\Gamma'$ be the subgraph of $\Gamma$ spanned by the vertices not in $V_0$; then $\Gamma'$ is centerless and $A(\Gamma) = A(\Gamma') \times \ZZ^{V_0}$, with center $\ZZ^{V_0}$.
	The members of $V$ in $\Gamma'$ form a compatible collection $V'$ with coincidence graph $\Gamma'$, and $(2k' + 7) D(\Gamma') + 4 \le (2k + 7) D(\Gamma) + 4$ for its size $k'$, so the restriction of $\Phi$ to $A(\Gamma')$ is injective by Lemma~\ref{lem:main-centerless} applied to $V'$.

	It remains to treat the center.
	Let $x = z y \in \ker \Phi$ with $z \in \ZZ^{V_0}$ and $y \in A(\Gamma')$.
	Then $\Phi(y) = \Phi(z)^{-1}$ commutes with $\Phi(A(\Gamma'))$, which is isomorphic to $A(\Gamma')$ and has trivial center.
	So $\Phi(y) = 1$, hence $y = 1$ and $x = z$.
	The members of $V_0$ are pairwise adjacent, so any two members of $V_0$ have disjoint representatives, and they are pairwise core-independent.
	So by Proposition~\ref{prop:pairwise-free-abelian} the twists $\{\delta_a \mid T_a \in V_0\}$ generate a free abelian group of rank $\lvert V_0 \rvert$.
	Hence $\Phi$ is injective on $\ZZ^{V_0}$, and $x = 1$.
\end{proof}

\section{The space of primitive curves} \label{sec:curve-space}

The goal of this section is Propositions~\ref{prop:C_compact} and~\ref{prop:C_compact-action}: the intersection numbers with essential spheres embed the set of primitive curves into a projective space, and the closure $\overline{\mathcal{C}}$ of the image is a compact Hausdorff $\Out(\FF)$-space.
Its dynamics is studied in Section~\ref{sec:dynamics}.

\subsection{The projective embedding} \label{subsec:embedding}

The goal of this subsection is Proposition~\ref{prop:injective}: the intersection numbers with essential spheres determine a primitive curve, so that the set of primitive curves embeds in a projective space, which is Hausdorff by Lemma~\ref{lem:hausdorff}.

Let $\mathcal{C}$ be the set of homotopy classes of unoriented primitive curves, and $\mathcal{S}$ the set of classes of essential spheres, which are also their isotopy classes by Theorem~\ref{thm:laudenbach}\eqref{item:laudenbach:spheres}.
By Lemma~\ref{lem:primitive-basis}, a class lies in $\mathcal{C}$ if and only if it is represented by an element of a free basis of $\FF$; so $\Out(\FF)$ acts on $\mathcal{C}$.
Let $\mathbb{R}_+ := [0, \infty)$, and define
\[
	i_\ast \colon \mathcal{C} \to \mathbb{R}_+^{\mathcal{S}} , \qquad
	i_\ast(\gamma)(S) := i(\gamma, S) .
\]
Let $P(\mathbb{R}_+^{\mathcal{S}})$ be the quotient of $\mathbb{R}_+^{\mathcal{S}} \setminus \{ 0 \}$ by multiplication by positive reals, with the quotient of the product topology, and $p$ the quotient map.
By the definition of a primitive curve, every $\gamma \in \mathcal{C}$ has $i(\gamma, S) = 1$ for some $S \in \mathcal{S}$, so $i_\ast(\gamma) \ne 0$ and $p \circ i_\ast$ is defined.
Since $i(\phi(\gamma), S) = i(\gamma, \phi^{-1}(S))$ for $\phi \in \Out(\FF)$, the map $p \circ i_\ast$ is equivariant for the action of $\Out(\FF)$ on $P(\mathbb{R}_+^{\mathcal{S}})$ by permutation of the coordinates.

\begin{proposition} \label{prop:injective}
	The map $p \circ i_\ast \colon \mathcal{C} \to P(\mathbb{R}_+^{\mathcal{S}})$ is injective.
\end{proposition}

\begin{proof}
	Let $\gamma, \gamma' \in \mathcal{C}$ with $(p \circ i_\ast)(\gamma) = (p \circ i_\ast)(\gamma')$, so that $i_\ast(\gamma) = \lambda \, i_\ast(\gamma')$ for some $\lambda > 0$.
	Fix a summand $\mathbb{S}^2 \times \mathbb{S}^1$ of $\MM$, a simple closed curve $c_1$ of the form $\{ y \} \times \mathbb{S}^1$ in this summand, and a separating essential sphere $s_1$ disjoint from $c_1$ cutting off this summand.
	The curve $c_1$ meets the essential sphere $\mathbb{S}^2 \times \{ z \}$ of the summand in one point, so its class $a_1$ is primitive.
	Under the identification $\pi_1(\MM) = \FF$ of Section~\ref{subsec:doubled_handlebody}, $a_1$ belongs to a free basis $\{ a_1, \dots, a_g \}$ by Lemma~\ref{lem:primitive-basis}, which we fix.
	By Lemma~\ref{lem:primitive-basis}, $\gamma$ is represented by an element $w$ of a free basis $\{ w, b_2, \dots, b_g \}$.
	Let $\phi$ be the automorphism of $\FF$ with $\phi(a_1) = w$ and $\phi(a_i) = b_i$ for $i \ge 2$.
	By Theorem~\ref{thm:laudenbach}\eqref{item:laudenbach:mcg}, $\phi$ is induced by a homeomorphism $h$ of $\MM$, so $h(c_1) \in \gamma$; put $c := h(c_1)$ and $s := h(s_1)$.

	The component $Q$ of $\MM \setminus s$ containing $c$ is homeomorphic to $(\mathbb{S}^2 \times \mathbb{S}^1) \setminus \{ x \}$ for a point $x$, and $\pi_1(Q)$ is generated by $c$; since $s$ is separating, $\pi_1(Q)$ is a free factor of $\FF$, so $\pi_1(Q) \to \pi_1(\MM)$ is injective.
	Let $s' \subset Q$ be an essential sphere meeting $c$ in one point, and $S'$ its class.
	Then $i(\gamma, [s]) = 0$, and $i(\gamma, S') \le 1$ since $c$ meets $s'$ in one point.
	The algebraic intersection number of $c$ with $s'$ is $\pm 1$ and depends only on the classes, so $i(\gamma, S') \ge 1$, and hence $i(\gamma, S') = 1$.
	By proportionality, $i(\gamma', [s]) = 0$ and $i(\gamma', S') = 1 / \lambda > 0$.

	So $\gamma'$ has a representative disjoint from a sphere homotopic to $s$, which is isotopic to $s$ by Theorem~\ref{thm:laudenbach}\eqref{item:laudenbach:spheres}.
	Applying the inverse of an ambient isotopy, which exists by the isotopy extension theorem, $\gamma'$ has a representative $c'$ disjoint from $s$.
	The curve $c'$ lies in $Q$, since a curve in the other component misses $s'$, contradicting $i(\gamma', S') > 0$.
	Then $c'$ represents a power of $c$ in $\pi_1(Q) \cong \ZZ$, and the power is $\pm 1$ since $\gamma'$ is not a proper power in $\FF$ by Theorem~\ref{thm:free}\eqref{item:free:primitive} and $\pi_1(Q) \to \FF$ is injective.
	Hence $\gamma = \gamma'$.
\end{proof}

By Proposition~\ref{prop:injective}, we regard $\mathcal{C}$ as a subset of $P(\mathbb{R}_+^{\mathcal{S}})$.

\begin{lemma} \label{lem:hausdorff}
	$P(\mathbb{R}_+^{\mathcal{S}})$ is Hausdorff.
\end{lemma}

\begin{proof}
	For $S_0, S_1 \in \mathcal{S}$ and $r > 0$, the sets $\{ [z] \mid z(S_0) > r \, z(S_1) \}$ and $\{ [z] \mid z(S_0) < r \, z(S_1) \}$ are disjoint and open, since their preimages in $\mathbb{R}_+^{\mathcal{S}} \setminus \{ 0 \}$ are open.
	Let $x$, $y$ be non-proportional, and choose $S_0$ with $x(S_0) > 0$.
	If $y(S_0) = 0$, choose $S_1$ with $y(S_1) > 0$ and $r > 0$ with $x(S_0) > r \, x(S_1)$; then the two sets separate $[x]$ from $[y]$.
	If $y(S_0) > 0$, scale so that $x(S_0) = y(S_0) = 1$; then $x(S_1) \ne y(S_1)$ for some $S_1$, say $x(S_1) < r < y(S_1)$, and the sets with $1/r$ in place of $r$ separate $[x]$ from $[y]$.
\end{proof}

\subsection{Compactness} \label{subsec:compactness}

For the isotopy class $\Sigma$ of a sphere system $\sigma$, we write $i(\gamma, \Sigma)$ for the minimum of $\lvert c \cap \sigma' \rvert$ over representatives $c \in \gamma$ and $\sigma' \in \Sigma$ transverse to each other.
Likewise $i(S, \Sigma)$ is the minimum of $\lvert \pi_0(s \cap \sigma') \rvert$ over transverse representatives $s \in S$, $\sigma' \in \Sigma$.
On a surface, the corresponding step is contained in the work of Fathi, Laudenbach and Po{\'e}naru \cite[Expos\'e~4]{MR3053012}, whose proof puts curves and a pants decomposition in minimal position simultaneously.
This is not available in $\MM$: we cannot put a curve, a sphere and a sphere system in minimal position at once.
Lemmas~\ref{lem:sum} and~\ref{lem:bounded} replace that step, using only one surface in minimal position with $\Sigma$ at a time and the simple connectivity of the components of $\MM$ cut along a representative of $\Sigma$.

\begin{lemma} \label{lem:sum}
	Let $\sigma$ be a sphere system with spheres $s_1, \dots, s_q$ of classes $S_1, \dots, S_q$, and $\Sigma$ its isotopy class.
	Then $i(\gamma, \Sigma) = \sum_{j} i(\gamma, S_j)$ for every homotopy class $\gamma$ of an essential simple closed curve.
\end{lemma}

\begin{proof}
	For $c \in \gamma$ and $\sigma' \in \Sigma$ transverse, with spheres $s'_j \in S_j$, we have
    \[
        \lvert c \cap \sigma' \rvert = \sum_{j = 1}^q \lvert c \cap s'_j \rvert \ge \sum_{j = 1}^q i(\gamma, S_j),
    \]
    which gives one inequality.

	For the other, let $X_\Sigma$ be the dual tree of $\sigma$, whose vertices are the components of $\widetilde\MM \setminus \rho^{-1}(\sigma)$ and whose edges are the lifts of the $s_j$.
	Let $g$ represent $\gamma$, with lift $\tilde\gamma$ and minimal subtree $A$ on $X_\Sigma$.
	Penetration is defined through the sides in $C$, so it does not depend on the other spheres of $\sigma$.
	Lemmas~\ref{lem:subtree-sides} and~\ref{lem:subtree-sides-penetrate}, stated for the dual tree of a single sphere or torus, hold for the dual tree $X_\Sigma$ of the system $\sigma$ with the same proofs, which use only the components of $\widetilde\MM$ minus the lifts and their halfspaces.
	So by Lemma~\ref{lem:subtree-sides-penetrate}, applied to $X_\Sigma$, the edges of $A$ that are lifts of $s_j$ are the lifts of $s_j$ penetrated by $\tilde\gamma$, and by Proposition~\ref{prop:orbit-count} they form $i(\gamma, S_j)$ orbits under $\langle g \rangle$.
	Hence $g$ is hyperbolic on $X_\Sigma$ with translation length $\sum_j i(\gamma, S_j)$, or fixes a vertex and $\sum_j i(\gamma, S_j) = 0$.

	Let $x_0$ be a point in a component $v_0$ of $\widetilde\MM \setminus \rho^{-1}(\sigma)$ lying on $A$.
	Let $\beta$ be a path from $x_0$ to $g x_0$ through the components along $A$, crossing each lift on $A$ between $v_0$ and $g v_0$ exactly once and transversely, and meeting no other lift.
	If $g$ fixes the vertex $v_0$, then $\beta$ is a path inside the single component $v_0$ and meets no lift.
	The translates $g^n \beta$ form a $\langle g \rangle$-invariant line, whose projection is a closed curve representing $\gamma$ and meeting $\sigma$ in $\sum_j i(\gamma, S_j)$ points.
	By general position in dimension $3$, it is homotopic to a simple closed curve with the same intersection with $\sigma$.
	So $i(\gamma, \Sigma) \le \sum_j i(\gamma, S_j)$.
\end{proof}

Let $\sigma$ be a maximal sphere system and $\Sigma$ its isotopy class.
By Theorem~\ref{thm:systems}, each component $N$ of $\MM$ cut along $\sigma$ is a holed sphere with three boundary spheres; in particular $N$ is simply connected and $H_1(N; \ZZ/2) = 0$.
Let $S \in \mathcal{S}$ with $i(S, \Sigma) > 0$, and let $s \in S$ be transverse to $\sigma$ with $\lvert \pi_0(s \cap \sigma) \rvert = i(S, \Sigma)$.
Such an $s$ exists: an isotopy of the finite disjoint union $\sigma$ of embedded spheres extends to an ambient isotopy by the isotopy extension theorem.
This ambient isotopy carries any representative of $\Sigma$, the isotopy class of the system $\sigma$, to $\sigma$.
Call the components of $s \cap N$ the \emph{pieces} in $N$.

\begin{lemma} \label{lem:pieces}
	Each $N$ contains at most $2 \, i(S, \Sigma)$ pieces.
\end{lemma}

\begin{proof}
	Each piece is a compact connected surface with nonempty boundary in $\partial N$, since $s$ is connected and meets $\sigma$.
	Each circle of $s \cap \sigma$ is a boundary circle of at most two pieces in $N$, one on each side of the sphere of $\sigma$ containing the circle.
	So $N$ contains at most $2 \, i(S, \Sigma)$ pieces.
\end{proof}

\begin{lemma} \label{lem:pieces-separate}
	Each piece separates $N$.
\end{lemma}

\begin{proof}
	Each piece $P$ is two-sided, so the number of intersection points with $P$ modulo $2$ is a homomorphism $H_1(N; \ZZ/2) \to \ZZ/2$, which is trivial.
	Hence $P$ separates $N$.
\end{proof}

\begin{lemma} \label{lem:arc-pieces}
	Each arc $a$ in $N$ with endpoints in $\partial N$, transverse to the pieces, is homotopic in $N$ relative to its endpoints to an arc meeting each piece transversely in at most one point.
\end{lemma}

\begin{proof}
	Suppose that $a$ meets a piece $P$ in two points $x_1$, $x_2$ that are consecutive on $a$ among the points of $a \cap P$.
	Since $N$ is simply connected, the subarc of $a$ from $x_1$ to $x_2$ is homotopic relative to its endpoints to a path $b$ in $P$.
	The arc $a$ crosses $P$ at $x_1$ from one side, which we call the first side, to the other, stays on the other side until $x_2$, and crosses back to the first side at $x_2$.
	Replace the subarc of $a$ from a point just before $x_1$ to a point just after $x_2$ by a path parallel to $b$ on the first side of $P$.

	The new arc is homotopic to $a$ relative to its endpoints, and it meets $P$ in two fewer points.
	It meets the other pieces only at points of $a$, since a path parallel to $b$ lies in a thin regular neighborhood of $P$, which misses the other pieces.
	Repeating this, we obtain an arc meeting each piece at most once.
\end{proof}

\begin{lemma} \label{lem:bounded-disjoint}
	Let $\sigma$ be a maximal sphere system and $s$ an essential sphere disjoint from $\sigma$.
	Then $s$ is isotopic to a sphere of $\sigma$.
\end{lemma}

\begin{proof}
	The sphere $s$ lies in a holed sphere $N$ with three boundary spheres, by Theorem~\ref{thm:systems}.
	Regard $N$ as $\mathbb{S}^3$ minus three open balls; by Alexander's theorem, $s$ bounds balls on both sides in $\mathbb{S}^3$, and each removed ball lies in one of these two balls.
	Since $s$ is essential in $\MM$, neither of the two balls is disjoint from all three removed balls, so one of the two balls, $B$, contains exactly one removed ball, and $B \cap N$ is homeomorphic to $\mathbb{S}^2 \times [0, 1]$.
	Hence $s$ is isotopic in $N$ to a boundary sphere of $N$, that is, to a sphere of $\sigma$.
\end{proof}

\begin{lemma} \label{lem:bounded}
	Let $\Sigma$ be the isotopy class of a maximal sphere system.
	Then $i(\gamma, S) \le i(\gamma, \Sigma) \, (2 \, i(S, \Sigma) + 1)$ for all $\gamma \in \mathcal{C}$ and $S \in \mathcal{S}$.
\end{lemma}

\begin{proof}
	Let $\sigma$ be a maximal sphere system representing $\Sigma$.
	Let $s \in S$ be transverse to $\sigma$ with $\lvert \pi_0(s \cap \sigma) \rvert = i(S, \Sigma)$, and let $c \in \gamma$ be transverse to $\sigma$ with $\lvert c \cap \sigma \rvert = i(\gamma, \Sigma)$.
	Since an ambient isotopy, given by the isotopy extension theorem, carries any representative of $\Sigma$ to $\sigma$, these minima are attained with $\sigma$ fixed.
	After a small perturbation of $c$ fixed near $\sigma$, we may assume that $c$ is transverse to $s$.

	Suppose first that $i(S, \Sigma) = 0$, so that $s \cap \sigma = \varnothing$.
	Then $s$ is isotopic to a sphere of $\sigma$ by Lemma~\ref{lem:bounded-disjoint}, so $i(\gamma, S) \le \lvert c \cap \sigma \rvert = i(\gamma, \Sigma)$.

	Suppose now that $i(S, \Sigma) > 0$.
	Let $N$ range over the components of $\MM$ cut along $\sigma$, each a simply connected holed sphere with three boundary spheres by Theorem~\ref{thm:systems}, and call the components of $s \cap N$ the pieces in $N$.
	We construct a curve $c'$ homotopic to $c$ as follows.
	The arcs of $c \cap N$, over all $N$, number $\lvert c \cap \sigma \rvert = i(\gamma, \Sigma)$.
	By Lemma~\ref{lem:arc-pieces}, every arc of $c \cap N$, for every $N$, is homotopic in $N$ relative to its endpoints to an arc meeting each piece in $N$ transversely in at most one point.
	Replace every arc of $c \cap N$ by such an arc, pushed into the interior of $N$ except near its endpoints, so that the new curve meets $\sigma$ only at the points of $c \cap \sigma$.
	Make the new arcs embedded and pairwise disjoint by general position in dimension $3$, without changing their intersections with $s$.
	The union $c'$ of the new arcs is a simple closed curve homotopic to $c$.

	We count the intersections of $c'$ with $s$.
	In each $N$, the number of points of $c' \cap s$ is at most the number of arcs times the number of pieces, which is at most $2 \, i(S, \Sigma)$ by Lemma~\ref{lem:pieces}.
	Summing over $N$, $i(\gamma, S) \le \lvert c' \cap s \rvert \le i(\gamma, \Sigma) \cdot 2 \, i(S, \Sigma)$.
\end{proof}

\begin{proposition} \label{prop:C_compact}
	The closure $\overline{\mathcal{C}}$ of $\mathcal{C}$ in $P(\mathbb{R}_+^{\mathcal{S}})$ is compact.
\end{proposition}

\begin{proof}
	Fix a maximal sphere system $\sigma$ with spheres of classes $S_1, \dots, S_q$, and let $\Sigma$ be its isotopy class.
	Every primitive curve meets every representative of $\Sigma$, since the components of its complement are simply connected; so $i(\gamma, \Sigma) > 0$ for $\gamma \in \mathcal{C}$.
	By Lemmas~\ref{lem:bounded} and~\ref{lem:sum}, the functions $f_\gamma := i_\ast(\gamma) / i(\gamma, \Sigma)$, $\gamma \in \mathcal{C}$, lie in
	\[
		F := \left\{ f \in \prod_{S \in \mathcal{S}} \bigl[ 0, \, 2 \, i(S, \Sigma) + 1 \bigr] \;\Big|\; \sum_{j=1}^{q} f(S_j) = 1 \right\} .
	\]

	The product is compact by Tychonoff's theorem, and $F$ is closed in the product and misses $0$, so $p(F)$ is compact.
	By Lemma~\ref{lem:hausdorff}, $p(F)$ is closed, and it contains $\mathcal{C}$, hence $\overline{\mathcal{C}}$.
	So $\overline{\mathcal{C}}$ is compact.
\end{proof}

Being a subspace of $P(\mathbb{R}_+^{\mathcal{S}})$, the closure $\overline{\mathcal{C}}$ is Hausdorff by Lemma~\ref{lem:hausdorff}.

\begin{proposition} \label{prop:C_compact-action}
	$\Out(\FF)$ acts on $\overline{\mathcal{C}}$ by homeomorphisms.
\end{proposition}

\begin{proof}
	$\Out(\FF)$ acts on $P(\mathbb{R}_+^{\mathcal{S}})$ by homeomorphisms permuting the coordinates and preserves $\mathcal{C}$, so it preserves $\overline{\mathcal{C}}$.
\end{proof}

\section{Dynamics on the space of primitive curves} \label{sec:dynamics}

The goal of this section is Theorem~\ref{thm:loxodromic}: a fully irreducible outer automorphism acts with north--south dynamics on the compact space $\overline{\mathcal{C}}$ of Section~\ref{sec:curve-space}.
On the way, Theorem~\ref{thm:multi-converge} describes the limits of the orbits of multitwists, which are weighted combinations of the cores of the tori.
These are the two model behaviors of the action of a mapping class group on the space of projective measured laminations, transferred to $\Out(\FF)$: multitwists act parabolically and fully irreducible elements act with north--south dynamics.

\subsection{Convergence under multitwists} \label{subsec:convergence}

We keep the notation of Section~\ref{subsec:growth}: $t_1, \dots, t_r$ are pairwise disjoint primitive tori with oriented cores, whose classes $T_1, \dots, T_r$ are pairwise core-independent, and $\gamma_v$ is the core of $T_v$.
For $\vec m = (m_1, \dots, m_r) \in \ZZ^r$ we write $\Delta_{\vec m} := \delta_{\vec T_1}^{\,m_1} \dots \delta_{\vec T_r}^{\,m_r}$ for the multitwist $\Delta$ of Section~\ref{subsec:growth} with these exponents; the factors commute, since the twist homeomorphisms along the disjoint tori $t_v$ have disjoint supports.

Let $\alpha \in \mathcal{C}$, and let $\vec m^{(n)}$ be a sequence in $\ZZ^r$ such that
\[
	w_n := \sum_{v=1}^{r} \lvert m_v^{(n)} \rvert \, i(\alpha, T_v) \to \infty
	\quad\text{and}\quad
	\frac{\lvert m_v^{(n)} \rvert \, i(\alpha, T_v)}{w_n} \to c_v \quad \text{for each}~ 1 \le v \le r .
\]
Then $c_v \ge 0$ and $\sum_v c_v = 1$, since for each $n$ the ratios are nonnegative and sum to $1$.

\begin{theorem} \label{thm:multi-converge}
	Let $t_1, \dots, t_r$ be pairwise disjoint primitive tori with oriented cores, whose classes $T_1, \dots, T_r$ are pairwise core-independent, with cores $\gamma_1, \dots, \gamma_r$, and write $\Delta_{\vec m} := \delta_{\vec T_1}^{\,m_1} \dots \delta_{\vec T_r}^{\,m_r}$ for $\vec m \in \ZZ^r$.
	Let $\alpha \in \mathcal{C}$, and let $\vec m^{(n)}$ be a sequence in $\ZZ^r$ such that $w_n := \sum_{v} \lvert m_v^{(n)} \rvert \, i(\alpha, T_v) \to \infty$ and $\lvert m_v^{(n)} \rvert \, i(\alpha, T_v) / w_n \to c_v$ for each $v$.
	Then
	\[
		(p \circ i_\ast) \bigl( \Delta_{\vec m^{(n)}}(\alpha) \bigr)
		\to \left[ \, \sum_{v=1}^{r} c_v \, i_\ast(\gamma_v) \, \right]
		\qquad \text{in } P(\mathbb{R}_+^{\mathcal{S}}) .
	\]
\end{theorem}

\begin{proof}
	Outer automorphisms carry free bases to free bases, so $\Delta_{\vec m^{(n)}}(\alpha) \in \mathcal{C}$.
	Put $\eta := \sum_v c_v \, i_\ast(\gamma_v) \in \mathbb{R}_+^{\mathcal{S}}$.
	By Proposition~\ref{prop:twist-growth}, for every $S \in \mathcal{S}$ there is $K_S$ with
	\[
		\left\lvert \, i \bigl( \Delta_{\vec m^{(n)}}(\alpha), S \bigr) - \sum_{v=1}^{r} \lvert m_v^{(n)} \rvert \, i(\alpha, T_v) \, i(\gamma_v, S) \, \right\rvert \le K_S
	\]
	for all $n$.
	Dividing by $w_n$ and letting $n \to \infty$, the error term tends to $0$ and the sum tends to $\eta(S)$.
	So $w_n^{-1} \, i_\ast(\Delta_{\vec m^{(n)}}(\alpha)) \to \eta$ in the product topology.

	The limit $\eta$ is nonzero.
	Choose $v_0$ with $c_{v_0} > 0$; since $\gamma_{v_0}$ is the class of a primitive curve, there is $S_0 \in \mathcal{S}$ with $i(\gamma_{v_0}, S_0) = 1$.
	Then $\eta(S_0) \ge c_{v_0} \, i(\gamma_{v_0}, S_0) > 0$.
	Since $p$ is continuous and invariant under scaling, $(p \circ i_\ast)(\Delta_{\vec m^{(n)}}(\alpha)) \to p(\eta)$.
\end{proof}

For the next three corollaries, let $\vec m \in \ZZ^r$ and put $\Delta := \Delta_{\vec m}$.

\begin{corollary} \label{cor:multitwist-parabolic}
	Let $\alpha \in \mathcal{C}$ with $i(\alpha, T_v) = 0$ for every $v$ with $m_v \ne 0$.
	Then $\Delta^n(\alpha) = \alpha$ for all $n \in \ZZ$.
\end{corollary}

\begin{proof}
	Let $v$ with $m_v \ne 0$.
	Since $i(\alpha, T_v) = i(\alpha, [t_v]) = 0$, some torus $t'$ homotopic to $t_v$ is disjoint from a representative $c$ of $\alpha$.
	By Section~\ref{subsec:classes}, $t'$ lies in the class $T_v$.
	We may orient the core of $t'$ so that $\vec t'$ lies in the oriented class $\vec T_v$, since every torus of $T_v$ admits a core orientation placing it in either oriented class over $T_v$.
	Alternatively, the core and transverse orientations carried from $\vec t_v$ by the homotopy place $\vec t'$ in $\vec T_v$ by Lemma~\ref{lem:oriented-homotopic}.

	Choose a regular neighborhood of $t'$ disjoint from $c$.
	Then the twist homeomorphism along $\vec t'$ with this regular neighborhood represents $\delta_{\vec T_v}$ by Section~\ref{subsec:twist-homeo} and fixes $c$.
	So $\delta_{\vec T_v}(\alpha) = \alpha$ for every $v$ with $m_v \ne 0$, and $\Delta^n(\alpha) = \alpha$.
\end{proof}

\begin{corollary} \label{cor:multitwist-parabolic-converge}
	Let $\alpha \in \mathcal{C}$ with $i(\alpha, T_v) > 0$ for some $v$ with $m_v \ne 0$.
	Then
	\[
		(p \circ i_\ast) \bigl( \Delta^n(\alpha) \bigr)
		\to \left[ \, \sum_{v=1}^{r} \lvert m_v \rvert \, i(\alpha, T_v) \, i_\ast(\gamma_v) \, \right]
		\quad \text{as} \quad \lvert n \rvert \to \infty .
	\]
\end{corollary}

\begin{proof}
	We have $\Delta^n = \Delta_{n \vec m}$ since the factors commute, and we apply Theorem~\ref{thm:multi-converge} separately for $n \to +\infty$ and for $n \to -\infty$, with $\vec m^{(n)} := n \vec m$: here $w_n = \lvert n \rvert \sum_v \lvert m_v \rvert \, i(\alpha, T_v) \to \infty$ and $c_v = \lvert m_v \rvert \, i(\alpha, T_v) / \sum_u \lvert m_u \rvert \, i(\alpha, T_u)$.
\end{proof}

\begin{corollary} \label{cor:multitwist-parabolic-fixed}
	Every point $\bigl[ \sum_v s_v \, i_\ast(\gamma_v) \bigr]$ with $s_1, \dots, s_r \ge 0$ not all zero is fixed by $\Delta$.
\end{corollary}

\begin{proof}
	We first show that $\delta_{\vec T_u}(\gamma_v) = \gamma_v$ for all $u$, $v$.
	Indeed, a core of $t_v$ lies on $t_v$, which is disjoint from the regular neighborhoods $N(t_u)$ with $u \ne v$, so the twist along $t_u$ fixes it.
	For $u = v$, a core of $t_v$ is homotopic to a curve on $\partial N(t_v)$, which the twist along $t_v$ fixes.
	Hence $\Delta(\gamma_v) = \gamma_v$, and since $i(\phi(\gamma), S) = i(\gamma, \phi^{-1}(S))$, the function $i_\ast(\gamma_v)$ is fixed by the coordinate action of $\Delta$; so is every point $\bigl[ \sum_v s_v \, i_\ast(\gamma_v) \bigr]$.
\end{proof}

\begin{corollary} \label{cor:core_simplex}
	Let $\alpha \in \mathcal{C}$ and $P_\alpha := \{ v \mid i(\alpha, T_v) > 0 \}$.
	Then for all $s_v \ge 0$, $v \in P_\alpha$, not all zero,
	\[
		\left[ \, \sum_{v \in P_\alpha} s_v \, i_\ast(\gamma_v) \, \right] \in \overline{\mathcal{C}} .
	\]
\end{corollary}

\begin{proof}
	For $v \in P_\alpha$ put $m_v^{(n)} := \lfloor n s_v / i(\alpha, T_v) \rfloor$, and put $m_v^{(n)} := 0$ for $v \notin P_\alpha$.
	Then $\lvert m_v^{(n)} \rvert \, i(\alpha, T_v) = n s_v + O(1)$, so $w_n = n \sum_v s_v + O(1) \to \infty$ and $c_v = s_v / \sum_u s_u$.
	By Theorem~\ref{thm:multi-converge}, the point $\bigl[ \sum_{v \in P_\alpha} s_v \, i_\ast(\gamma_v) \bigr]$ is a limit of points of $\mathcal{C}$.
\end{proof}

\subsection{Fully irreducible elements} \label{subsec:iwip}

The goal of this subsection is Theorem~\ref{thm:loxodromic}, which we deduce from the north--south dynamics of a fully irreducible element on the minimal set of projective currents through the continuous equivariant surjection onto $\overline{\mathcal{C}}$ of Lemmas~\ref{prop:surjection} and~\ref{prop:surjection-equivariant} and Proposition~\ref{prop:surjection-onto}.

A \emph{geodesic current} on $\FF$ is an $\FF$-invariant, flip-invariant, positive Radon measure on the space $\partial^2 \FF$ of pairs of distinct ends of $\FF$.
Let $\operatorname{Curr}(\FF)$ be the space of geodesic currents and $\mathbb{P}\operatorname{Curr}(\FF)$ its projectivization, both with the weak topology, on which $\Out(\FF)$ acts.
Let $\overline{cv(\FF)}$ be the closure of unprojectivized Outer space.
By Kapovich and Lustig \cite{MR2496058}, there is a unique continuous $\Out(\FF)$-equivariant intersection form
\[
	\langle \, \cdot \, , \, \cdot \, \rangle \colon \overline{cv(\FF)} \times \operatorname{Curr}(\FF) \to \mathbb{R}_+ ,
\]
$\mathbb{R}_+$-homogeneous in the first argument and $\mathbb{R}_+$-linear in the second, with $\langle Y, \eta_g \rangle = \lVert g \rVert_Y$ for all $Y \in \overline{cv(\FF)}$ and $g \in \FF$, where $\lVert g \rVert_Y$ is the translation length of $g$ on $Y$ and $\eta_g$ is the counting current of $g$.

Each $S \in \mathcal{S}$ determines a one-edge free splitting of $\FF$ whose Bass--Serre tree is the dual tree $X_S$; it is a very small simplicial $\FF$-tree, hence lies in $\overline{cv(\FF)}$ by Cohen and Lustig \cite{MR1341810}.
Each $\gamma \in \mathcal{C}$ determines the counting current $\eta_\gamma$ of a representative.
By Proposition~\ref{prop:orbit-count-translation}, $i(\gamma, S)$ is the translation length of $\gamma$ on $X_S$ with unit edge lengths, so
\[
	i(\gamma, S) = \langle X_S, \eta_\gamma \rangle .
\]

Let $\mathcal{M} \subseteq \mathbb{P}\operatorname{Curr}(\FF)$ be the closure of the set of projective classes $[\eta_\gamma]$ of the counting currents of the primitive elements $\gamma$ of $\FF$.
It is nonempty, closed and $\Out(\FF)$-invariant, since outer automorphisms carry primitive elements to primitive elements, and the primitive classes are dense in $\mathcal{M}$ by definition.
We claim that $\mathcal{M}$ is a minimal nonempty closed $\Out(\FF)$-invariant subset of $\mathbb{P}\operatorname{Curr}(\FF)$.
The set $\mathcal{M}$ was introduced by Martin \cite{MR2693216}, who proved that every $\Out(\FF)$-orbit in $\mathcal{M}$ is dense, in every rank $g \ge 2$; see \cite[Proposition~4.5]{MR2966693}.
For $g \ge 3$, Kapovich and Lustig \cite[Theorem~B]{MR2322181} show moreover that $\mathcal{M}$ is the unique minimal nonempty closed $\Out(\FF)$-invariant subset of $\mathbb{P}\operatorname{Curr}(\FF)$.

For $g = 2$, uniqueness fails, since the class $[\eta_{[a,b]}]$ of the commutator of a free basis is fixed by $\Out(\FF_2)$, but minimality holds.
Indeed, $\FF_2$ is the fundamental group of the once-punctured torus $S_{1,1}$, with $\Out(\FF_2) \cong \mathrm{GL}(2, \ZZ)$ its mapping class group, and the primitive elements of $\FF_2$ are the classes of the simple closed curves on $S_{1,1}$.
The space of measured laminations of $S_{1,1}$ maps continuously and equivariantly into $\operatorname{Curr}(\FF_2)$, carrying a weighted simple closed curve to the multiple of its counting current \cite[Section~2.2]{MR3370027}.
The projective classes of essential simple closed curves are dense in $\mathbb{P}\mathcal{ML}(S_{1,1}) \cong \mathbb{RP}^1$, their slopes being the points of $\mathbb{Q} \cup \{ \infty \}$, so $\mathcal{M}$ is the image of the compact space $\mathbb{P}\mathcal{ML}(S_{1,1})$.
On this space $\mathrm{GL}(2, \ZZ)$ acts through $\mathrm{PGL}(2, \ZZ)$ on $\mathbb{R} \cup \{ \infty \}$ by fractional linear transformations, and every orbit of this action is dense, so $\mathcal{M}$ is minimal.

Call a current $\mu$ \emph{non-degenerate} if $\langle X_S, \mu \rangle > 0$ for some $S \in \mathcal{S}$.

\begin{lemma} \label{lem:nondegeneracy}
	Every class in $\mathcal{M}$ is non-degenerate.
\end{lemma}

\begin{proof}
	Let $\mathcal{Z}$ be the set of classes $[\mu] \in \mathbb{P}\operatorname{Curr}(\FF)$ with $\langle X_S, \mu \rangle = 0$ for all $S \in \mathcal{S}$.
	Each pairing $\langle X_S, \cdot \rangle$ is continuous and homogeneous, so $\mathcal{Z}$ is closed.
	The family $\{ X_S \}$ is $\Out(\FF)$-invariant, with $\phi \, X_S = X_{\phi(S)}$, and the intersection form is equivariant; so $\mathcal{Z}$ is $\Out(\FF)$-invariant.
	Let $\gamma \in \mathcal{C}$ and let $S$ be the class of an essential sphere with $i(\gamma, S) = 1$, which exists by the definition of a primitive curve.
	Then $\langle X_S, \eta_\gamma \rangle = 1$, so $[\eta_\gamma] \in \mathcal{M} \setminus \mathcal{Z}$.
	Hence $\mathcal{Z} \cap \mathcal{M}$ is a closed invariant proper subset of $\mathcal{M}$, which is empty by minimality.
\end{proof}

By Lemma~\ref{lem:nondegeneracy}, the coordinate map
\[
	j \colon \mathcal{M} \to P(\mathbb{R}_+^{\mathcal{S}}) , \qquad
	[\mu] \mapsto \bigl[ \, (\langle X_S, \mu \rangle)_{S \in \mathcal{S}} \, \bigr]
\]
is defined, and it agrees with $p \circ i_\ast$ on the dense set of primitive classes.

The space $\mathcal{M}$ is closed in the compact metrizable space $\mathbb{P}\operatorname{Curr}(\FF)$, hence compact and Hausdorff, and $\overline{\mathcal{C}}$ is compact and Hausdorff by Proposition~\ref{prop:C_compact} and Lemma~\ref{lem:hausdorff}.

\begin{lemma} \label{prop:surjection}
	The map $j$ is continuous.
\end{lemma}

\begin{proof}
	Continuity follows from the continuity of each pairing $\langle X_S, \cdot \rangle$.
\end{proof}

\begin{lemma} \label{prop:surjection-equivariant}
	The map $j$ is $\Out(\FF)$-equivariant.
\end{lemma}

\begin{proof}
	For $\phi \in \Out(\FF)$, equivariance of the intersection form gives $\langle X_S, \phi \mu \rangle = \langle \phi^{-1} X_S, \mu \rangle = \langle X_{\phi^{-1}(S)}, \mu \rangle$, which is the coordinate action of $\phi$ on $P(\mathbb{R}_+^{\mathcal{S}})$; so $j$ is equivariant.
\end{proof}

\begin{proposition} \label{prop:surjection-onto}
	$j(\mathcal{M}) = \overline{\mathcal{C}}$.
\end{proposition}

\begin{proof}
	The map $j$ is continuous by Lemma~\ref{prop:surjection}.
	The image $j(\mathcal{M})$ is compact, hence closed in the Hausdorff space $P(\mathbb{R}_+^{\mathcal{S}})$, and contains the dense subset $\mathcal{C}$ of $\overline{\mathcal{C}}$; and $j(\mathcal{M}) \subseteq \overline{\mathcal{C}}$ since $j$ is continuous and the primitive classes are dense in $\mathcal{M}$.
	So $j(\mathcal{M}) = \overline{\mathcal{C}}$.
\end{proof}

A homeomorphism $\phi$ of a compact Hausdorff space $X$ acts with \emph{north--south dynamics} with attracting point $x_+$ and repelling point $x_-$ if the following holds.
We have $\phi(x_\pm) = x_\pm$ and $x_+ \ne x_-$, and for every compact $K \subseteq X \setminus \{ x_- \}$ and every neighborhood $V$ of $x_+$ we have $\phi^n(K) \subseteq V$ for all large $n$.
The same is required of $\phi^{-1}$ with $x_-$ and $x_+$ exchanged.

Let $\phi \in \Out(\FF)$ be fully irreducible, with expansion factors $\lambda_+ := \lambda(\phi) > 1$ and $\lambda_- := \lambda(\phi^{-1})^{-1} < 1$, and let $\mu_\pm$ be its stable and unstable currents, introduced by Martin \cite{MR2693216}; see Uyanik \cite[Lemma~3.7 and Remark~3.8]{MR3273532}, where $\mu_+$ is constructed from a train track representative with $\phi \, \mu_+ = \lambda_+ \mu_+$, and $\mu_- := \mu_+(\phi^{-1})$, so that $\phi \, \mu_- = \lambda_- \mu_-$.
These equations hold as currents, not merely projectively, $[\mu_+] \ne [\mu_-]$, and $\mu_\pm(\phi^m) = \mu_\pm(\phi)$ up to scaling for $m \ge 1$.
When $\phi$ is induced by a pseudo-Anosov homeomorphism of a surface, $\mu_\pm$ are its stable and unstable measured laminations regarded as currents \cite[Section~3]{MR3370027}.

The inverse $\phi^{-1}$ is fully irreducible, since the definition is symmetric under inversion, and by the definition of $\mu_-$ its stable and unstable currents are $\mu_+(\phi^{-1}) = \mu_-$ and $\mu_-(\phi^{-1}) = \mu_+$.
So every statement below about $\phi$, $\mu_+$ and $\mu_-$ applies to $\phi^{-1}$, $\mu_-$ and $\mu_+$.

\begin{theorem}[Martin \cite{MR2693216}, Uyanik \cite{MR3273532, MR3370027}] \label{thm:NS-currents}
	Let $\phi \in \Out(\FF)$ be fully irreducible.
	Then $[\mu_+], [\mu_-] \in \mathcal{M}$.
\end{theorem}

\begin{theorem}[Martin \cite{MR2693216}, Uyanik \cite{MR3273532, MR3370027}] \label{thm:NS-currents-dynamics}
	Let $\phi \in \Out(\FF)$ be fully irreducible.
	Then $\phi$ acts on $\mathcal{M}$ with north--south dynamics with attracting point $[\mu_+]$ and repelling point $[\mu_-]$.
\end{theorem}

We make the attribution precise.
In the atoroidal case, Martin \cite{MR2693216} and Uyanik \cite{MR3273532} prove north--south dynamics on all of $\mathbb{P}\operatorname{Curr}(\FF)$, and the restriction to $\mathcal{M}$ is immediate.
This case does not occur for $g = 2$, since every element of $\Out(\FF_2)$ preserves the conjugacy class of the commutator of a free basis up to inversion.
Uyanik states the non-atoroidal case for every fully irreducible $\phi$ \cite[Corollary~4.5]{MR3370027}; his proof uses neither the uniqueness of $\mathcal{M}$ nor a bound on the rank, but it is written for orientable surfaces and for homeomorphisms fixing the boundary pointwise.
We give the proof with the changes needed for nonorientable surfaces and orientation-reversing homeomorphisms.
Both theorems rest on the following lemma.

\begin{lemma} \label{thm:NS-currents-attracting}
	Let $\phi \in \Out(\FF)$ be fully irreducible, let $K \subseteq \mathcal{M} \setminus \{ [\mu_-] \}$ be compact, and let $V$ be a neighborhood of $[\mu_+]$ in $\mathbb{P}\operatorname{Curr}(\FF)$.
	Then $\phi^n(K) \subseteq V$ for all large $n$.
\end{lemma}

\begin{proof}
	Suppose first that $\phi$ is atoroidal.
	Then $\phi$ acts on all of $\mathbb{P}\operatorname{Curr}(\FF)$ with north--south dynamics with attracting point $[\mu_+]$ and repelling point $[\mu_-]$, by Martin \cite{MR2693216} and Uyanik \cite{MR3273532}.
	The set $K$ is a compact subset of $\mathbb{P}\operatorname{Curr}(\FF) \setminus \{ [\mu_-] \}$, so $\phi^n(K) \subseteq V$ for all large $n$.

		Suppose now that $\phi$ is not atoroidal.
	By Bestvina and Handel \cite{MR1147956}, $\phi$ is induced by a pseudo-Anosov homeomorphism $f$ of a compact surface $S_\phi$ with one boundary component $c$ and $\pi_1(S_\phi) = \FF$.

	We first reduce to the case where $f$ preserves $c$ with its orientation, by replacing $\phi$ by $\phi^2$.
	Then $f$ is orientation-preserving and isotopic to a homeomorphism fixing the boundary pointwise, which is the setting of \cite[Section~3]{MR3370027}.
	The same holds for the lift $\hat f$ used below when $S_\phi$ is nonorientable.

	The reduction is harmless.
	Indeed, $\mu_\pm(\phi^2) = \mu_\pm(\phi)$ up to scaling, and $\phi$ fixes $[\mu_\pm]$.
	Let $K \subseteq \mathcal{M} \setminus \{ [\mu_-] \}$ be compact and $V$ a neighborhood of $[\mu_+]$.
	The set $K \cup \phi(K)$ is again compact and misses $[\mu_-]$.
	So if the claim holds for $\phi^2$, then $\phi^{2n}(K \cup \phi(K)) \subseteq V$ for all large $n$, that is, $\phi^m(K) \subseteq V$ for all large $m$.

	Let $H_0$ be the segment of $\mathbb{P}\operatorname{Curr}(\FF)$ spanned by $[\mu_-]$ and $[\eta_c]$, and $H_1$ the segment spanned by $[\mu_+]$ and $[\eta_c]$.
	Suppose that $S_\phi$ is orientable.
	Then \cite[Theorem~3.5]{MR3370027}, which allows any number of boundary components, gives the following.
	For every compact $K$ disjoint from $H_0$ and every neighborhood $V$ of $[\mu_+]$, we have $\phi^n(K) \subseteq V$ for all large $n$.
	The symmetric statement holds for $\phi^{-1}$, $H_1$ and $[\mu_-]$.

	We extend this statement to nonorientable $S_\phi$ by passing to the orientation double cover.
	If $S_\phi$ is nonorientable, let $\hat S \to S_\phi$ be its orientation double cover, an orientable surface of genus $g - 1 \ge 1$ with two boundary components $c_1$, $c_2$ over $c$.
	Let $\hat f$ be the lift of $f$ preserving $c_1$ and $c_2$, which is pseudo-Anosov.
	The subgroup $\hat\FF := \pi_1(\hat S)$ has index two in $\FF$, so $\partial^2 \hat\FF = \partial^2 \FF$, and $\operatorname{Curr}(\FF)$ is the closed subspace of $\FF$-invariant currents in $\operatorname{Curr}(\hat\FF)$.
	On this subspace $\hat f$ acts as $\phi$, with stable and unstable currents $\mu_\pm$, and $\eta_c = \eta_{c_1} + \eta_{c_2}$.

	The elements of $\FF \setminus \hat\FF$ exchange $\eta_{c_1}$ and $\eta_{c_2}$ and fix $\mu_-$.
	The $\FF$-invariant points of the simplex spanned by $[\mu_-]$, $[\eta_{c_1}]$, $[\eta_{c_2}]$ are those with equal coefficients on $\eta_{c_1}$ and $\eta_{c_2}$, which form $H_0$.
	So the exceptional set of \cite[Theorem~3.5]{MR3370027} for $\hat f$, which is this simplex, meets $\mathbb{P}\operatorname{Curr}(\FF)$ in $H_0$, and the same statement holds for $\phi$.

		In both cases, the segments $H_0$ and $H_1$ meet $\mathcal{M}$ at most in $[\mu_-]$ and $[\mu_+]$.
	To see this, we use Whitehead graphs.

	Fix a free basis of $\FF$.
	The Whitehead graph of the support of a current $\mu$ has the letters of the basis and their inverses as vertices, and an edge for each reduced two-letter word whose cylinder set has positive $\mu$-measure.
	For the counting current of an element, this is the Whitehead graph of the element.

	We claim that the Whitehead graph of the support of a current in $\mathcal{M}$ is disconnected or has a cut vertex.
	This is the criterion used in the proof of \cite[Corollary~4.5]{MR3370027}, and we include the argument.
	By Whitehead's lemma \cite{MR1575455}, the Whitehead graph of a primitive element is disconnected or has a cut vertex.

	Let $\mu$ be a current with $[\mu] \in \mathcal{M}$.
	It is a limit of currents $\lambda_n \eta_{\gamma_n}$ with $\gamma_n$ primitive and $\lambda_n > 0$, since $\mathbb{P}\operatorname{Curr}(\FF)$ is metrizable and the primitive classes are dense in $\mathcal{M}$.
	Cylinder sets are compact and open, so the measure of a cylinder is continuous on $\operatorname{Curr}(\FF)$.
	Hence a cylinder of positive $\mu$-measure has positive $\eta_{\gamma_n}$-measure for all large $n$.

	Since there are finitely many reduced two-letter words, the Whitehead graph of the support of $\mu$ is therefore a subgraph, with the same vertices, of the Whitehead graph of some $\gamma_n$.
	Removing edges from a graph that is disconnected or has a cut vertex gives a graph with the same property.
	This proves the claim.

	In a suitable basis, the boundary word of $S_\phi$ is $[a_1, b_1] \dots [a_h, b_h]$ with $g = 2h$ if $S_\phi$ is orientable, and $a_1^2 \dots a_g^2$ otherwise.
	Its Whitehead graph is a cycle through all vertices; in the second case it is the cycle $a_1^{-1}, a_1, a_2^{-1}, a_2, \dots, a_g^{-1}, a_g, a_1^{-1}$.
	The support of a point of $H_0 \cup H_1$ other than $[\mu_\pm]$ contains that of $\eta_c$, so its Whitehead graph contains this cycle.
	Such a graph is connected and has no cut vertex, so the point does not lie in $\mathcal{M}$ by the claim.

	Hence the compact set $K \subseteq \mathcal{M} \setminus \{ [\mu_-] \}$ misses $H_0$, and $\phi^n(K) \subseteq V$ for all large $n$.
\end{proof}

\begin{proof}[Proof of Theorem~\ref{thm:NS-currents}]
	The set $\mathcal{M}$ is infinite, since it contains the classes of the counting currents of the infinitely many conjugacy classes of primitive elements.
	Choose $x \in \mathcal{M}$ with $x \ne [\mu_-]$.
	By Lemma~\ref{thm:NS-currents-attracting} applied to $K = \{ x \}$, $\phi^n(x) \to [\mu_+]$, so $[\mu_+] \in \mathcal{M}$ since $\mathcal{M}$ is closed and invariant.
	The same argument for $\phi^{-1}$, whose stable current is $\mu_-$, gives $[\mu_-] \in \mathcal{M}$.
\end{proof}

\begin{proof}[Proof of Theorem~\ref{thm:NS-currents-dynamics}]
	The points $[\mu_\pm]$ lie in $\mathcal{M}$ by Theorem~\ref{thm:NS-currents}, are fixed by $\phi$ since $\phi \, \mu_\pm = \lambda_\pm \mu_\pm$, and are distinct.
	Let $K \subseteq \mathcal{M} \setminus \{ [\mu_-] \}$ be compact and $V$ a neighborhood of $[\mu_+]$ in $\mathcal{M}$, so that $V = V' \cap \mathcal{M}$ for a neighborhood $V'$ of $[\mu_+]$ in $\mathbb{P}\operatorname{Curr}(\FF)$.
	By Lemma~\ref{thm:NS-currents-attracting}, $\phi^n(K) \subseteq V'$ for all large $n$, and $\phi^n(K) \subseteq \mathcal{M}$ since $\mathcal{M}$ is invariant; so $\phi^n(K) \subseteq V$.
	The same argument for $\phi^{-1}$, whose stable and unstable currents are $\mu_-$ and $\mu_+$, gives the statement required of $\phi^{-1}$.
\end{proof}

\begin{lemma} \label{lem:descent}
	Let $q \colon X \to Y$ be a continuous surjection of compact Hausdorff spaces, equivariant for a homeomorphism $\phi$ of $X$ and a homeomorphism $\phi$ of $Y$.
	Suppose that $\phi$ acts on $X$ with north--south dynamics, with attracting fixed point $x_+$ and repelling fixed point $x_-$, and that $q(x_+) \ne q(x_-)$.
	Then $\phi$ acts on $Y$ with north--south dynamics, with fixed points $y_\pm := q(x_\pm)$.
\end{lemma}

\begin{proof}
	Equivariance gives $\phi(y_\pm) = q(\phi(x_\pm)) = y_\pm$, and $y_+ \ne y_-$ by hypothesis.
	Let $K \subseteq Y \setminus \{ y_- \}$ be compact and $V$ an open neighborhood of $y_+$.
	Then $q^{-1}(K)$ is compact and misses $x_-$, and $U := q^{-1}(V)$ is an open neighborhood of $x_+$.
	So there is $N$ with $\phi^n(q^{-1}(K)) \subseteq U$ for all $n \ge N$.
	By surjectivity and equivariance, $\phi^n(K) = q(\phi^n(q^{-1}(K))) \subseteq q(U) \subseteq V$ for $n \ge N$.
	The statement for $\phi^{-1}$ and $y_-$ follows in the same way.
\end{proof}

\begin{lemma} \label{lem:reciprocity}
	$j([\mu_+]) \ne j([\mu_-])$.
\end{lemma}

\begin{proof}
	Write $f_\pm(S) := \langle X_S, \mu_\pm \rangle$.
	By equivariance and $\phi \, \mu_\pm = \lambda_\pm \, \mu_\pm$,
	\[
		f_\pm(\phi^{-1}(S)) = \langle X_{\phi^{-1}(S)}, \mu_\pm \rangle = \langle X_S, \phi \, \mu_\pm \rangle = \lambda_\pm \, f_\pm(S)
	\]
	for all $S \in \mathcal{S}$.
	Suppose that $j([\mu_+]) = j([\mu_-])$, so that $f_+ = \kappa f_-$ for some $\kappa > 0$.
	Evaluating at $\phi^{-1}(S)$,
	\[
		\lambda_+ \, f_+(S) = f_+(\phi^{-1}(S)) = \kappa \, f_-(\phi^{-1}(S)) = \kappa \, \lambda_- \, f_-(S) = \lambda_- \, f_+(S) .
	\]
	Since $\lambda_+ \ne \lambda_-$, this gives $f_+ = 0$, contradicting Lemma~\ref{lem:nondegeneracy}, since $[\mu_+] \in \mathcal{M}$ by Theorem~\ref{thm:NS-currents}.
\end{proof}

\begin{theorem} \label{thm:loxodromic}
	Let $\phi \in \Out(\FF)$ be fully irreducible.
	Then $\phi$ acts on the compact space $\overline{\mathcal{C}}$ with north--south dynamics, with attracting fixed point $j([\mu_+])$ and repelling fixed point $j([\mu_-])$.
\end{theorem}

\begin{proof}
	By Theorem~\ref{thm:NS-currents-dynamics}, $\phi$ acts on $\mathcal{M}$ with north--south dynamics.
	By Lemmas~\ref{prop:surjection} and~\ref{prop:surjection-equivariant} and Proposition~\ref{prop:surjection-onto}, $j \colon \mathcal{M} \to \overline{\mathcal{C}}$ is a continuous equivariant surjection of compact Hausdorff spaces, and $j([\mu_+]) \ne j([\mu_-])$ by Lemma~\ref{lem:reciprocity}.
	Lemma~\ref{lem:descent} gives the conclusion.
\end{proof}

\bibliographystyle{amsalpha}
\bibliography{references}

\end{document}